[1994/12/01]% LaTeX date must December 1994 or later
\documentclass[10pt,reqno]{amsart}
\usepackage{a4wide}
\usepackage[english, activeacute]{babel}
\usepackage{amsmath,amsthm,amsxtra}
\usepackage{epsfig}
\usepackage{amssymb}
\usepackage{latexsym}
\usepackage{amsfonts}
\usepackage{hyperref}
\title{On the inelastic 2-soliton collision for gKdV equations with general nonlinearity}
\author{Claudio Mu\~noz C.}
\address{Universit\'e de Versailles Saint-Quentin-en-Yvelines \\ LMV-UMR 8100, 45 av. des Etats-Unis, 78035 Versailles cedex, France}
\email{Claudio.Munoz@math.uvsq.fr}
\date{February, 2009}
\thanks{This research was supported in part by a CONICYT-Chile and an \emph{Allocation de Recherche} grants}
\chardef\bslash=`\\ % p. 424, TeXbook
\newtheorem{thm}{Theorem}[section]

\newtheorem{lem}[thm]{Lemma}
\newtheorem{prop}[thm]{Proposition}

\theoremstyle{definition}
\newtheorem{defn}{Definition}[section]

\theoremstyle{remark}
\newtheorem{rem}{Remark}[section]

\newtheorem{Cl}{Claim}

\newcommand{\R}{\mathbb{R}}

\newcommand{\wqs}{{Q}_c}
\newcommand{\ys}{y_c}
\newcommand{\al}{\alpha}
\newcommand{\ga}{\gamma}

\def\bm{\left( \begin{array}{cc}}
\def\endm{\end{array}\right)}

 \providecommand{\abs}[1]{\lvert#1 \rvert}
 \providecommand{\norm}[1]{\lVert#1 \rVert}

\newcommand{\ve}{\varepsilon}

\newcommand{\be}{\begin{equation}}
\newcommand{\ee}{\end{equation}}
\newcommand{\ba}{\begin{array}}
\newcommand{\ea}{\end{array}}
\newcommand{\bea}{\begin{eqnarray}}
\newcommand{\eea}{\end{eqnarray}}
\newcommand{\bee}{\begin{eqnarray*}}
\newcommand{\eee}{\end{eqnarray*}}

\newcommand{\eval}[2][\right]{\relax
  \ifx#1\right\relax \left.\fi#2#1\rvert}

\let\abs=\envert

\let\norm=\enVert

\begin{document}
\begin{abstract}
We study the problem of 2-soliton collision for the generalized Korteweg-de Vries equations, completing some recent works of Y. Martel and F. Merle \cite{MMcol1, MMcol2}. We classify the nonlinearities  for which collisions are elastic or inelastic. Our main result states that in the case of small solitons, with one soliton smaller than the other one, the unique nonlinearities allowing a perfectly elastic collision are precisely the integrable cases, namely the quadratic (KdV), cubic (mKdV) and Gardner nonlinearities. 
\end{abstract}
\maketitle \markboth{2-soliton collision for gKdV equations} {Claudio Mu\~noz}
\renewcommand{\sectionmark}[1]{}

%\tableofcontents

%%%%%%%%%%%%%%%%%%%%%%%%%%%%%%%%%%%%%%%%%%%%%%%%%%%%%%%%%%
%%%%%%%%%%%%%%%%%%%%%%%%%%%%%%%%%%%%%%%%%%%%%%%%%%%%%%%%%%

\section{Introduction and Main Results}

\medskip

In this work we consider the \emph{generalized Korteweg-de Vries equation} (gKdV) on the real line
\begin{equation}\label{gKdVm1}
u_t + (u_{xx} + f(u))_x =0,\quad \hbox{ in }  \R_t\times \R_x. 
\end{equation}

\medskip

Here $u=u(t,x)$  is a real-valued function, and $f: \R\to \R$ a nonlinear function, often refered as the \emph{nonlinearity} of (\ref{gKdVm1}). %$C^{m+2}(\R) function, where we assume, to be explained later, that $m=2$ or $3$. 
This equation represents a mathematical generalization of the \emph{Korteweg-de Vries} equation (KdV), namely the case $f(s) = s^2$,
\be\label{KdV}
u_t + (u_{xx} +u^2)_x =0, \quad \hbox{ in }\R_t\times \R_x;
\ee
other physically important cases are the cubic one $f(s)=s^3$, and the \emph{quadratic-cubic} nonlinearity, namely $f(s) = s^2 - \mu s^3$, $\mu\in \R$. In the former case, the equation (\ref{gKdVm1}) is often refered as the (focusing) \emph{modified} KdV equation (mKdV), and in the latter, it is known as the \emph{Gardner} equation. 

\smallskip

Concerning the KdV equation, it arises in Physics as a model of propagation of dispersive long waves, as was pointed out by J. S. Russel in 1834 \cite{Miura}. The exact formulation of the KdV equation comes from Korteweg and de Vries (1895) \cite{KdV}. This equation was re-discovered in a numerical work by N. Zabusky and M. Kruskal in 1965 \cite{KZ}.

\smallskip

After this work, a great amount of literature has emerged, physical, numerical and mathematical, for the study of this equation, see for example \cite{Bo3, Bo14, Shih34, Li18, HaSa12, Mi30, Miura}. Although under different points of view, among the main topics treated are the following: existence of explicit solutions and their stability, local and global well posedness, long time behavior properties and, of course, related generalized models, \emph{hierarchies}  and their properties.

\smallskip

This continuous, focused research on the KdV equation can be in part explained by some striking algebraic properties. One of the first properties is the existence of localized, rapidly decaying, stable and smooth solutions called \emph{solitons}. Given three real numbers $t_0, x_0$ and $c>0$, solitons are solutions of (\ref{KdV}) of the form
\be\label{(3)}
u(t,x):= Q_c(x-x_0-c(t-t_0)),  \quad Q_c(s):=c Q(c^{1/2} s),
\ee  
and where $Q$ satisfies the second order nonlinear differential equation
$$
Q'' -Q + Q^2 =0, \quad Q(x) =\frac{3}{2\cosh^2(\frac 12 x)}.
$$

The 3-parameter family of solitons (\ref{(3)}) contains three important symmetries of the equation, namely \emph{scaling} and \emph{translation} in space and time invariances. From the No\"ether theorem, these two last symmetries are related to \emph{conserved quantities}, invariant under the KdV flow, usually called \emph{Mass} and \emph{Energy}, represented below in (\ref{M})-(\ref{E}) (in a general form). Moreover, due to the mass and energy conservation, the Sobolev space $H^1(\R)$ appears as an ideal space to study long time properties of KdV.    

Even more striking is the fact that KdV, mKdV and the Gardner equation, being infinite dimensional dynamical systems, possess an \emph{infinity number of conserved quantities}, a consequence of the so-called \emph{complete integrability} property. This one is closely related to the existence of a \emph{Lax pair} for these equations  (see Lax, \cite{LAX1}). Another important property is the following well known fact: given any Schwartz initial data, the corresponding solution to the Cauchy problem for (\ref{KdV}) exists globally in time and decouples, as $t\to +\infty$, into a \emph{radiation} part going leftward plus a nonlinear \emph{multisoliton} component going to the right, see \cite{S}. %The radiative term disperses like a solution of the linear equation associated to (\ref{KdV}), the \emph{Airy} equation.  

\smallskip

%Many of these properties were discovered looking at 
The dynamical problem of 2-soliton collision is a classical problem in nonlinear wave propagation (see \cite{MMcol1} for a review and references therein). By 2-soliton collision we mean the following problem: given two solitons, solutions of (\ref{gKdVm1}), largely separated at some early time and having different velocities, we expect that they have to collide at some finite time. The resulting solution after the collision is precisely the object of study. In particular, one considers if any change in size, position, or shape, even destruction of the solitons, after some large time, may be present.    

\smallskip

Let us review some relevant works in this direction. First, the works of Fermi, Pasta and Ulam \cite{FPU} and Zabusky and Kruskal \cite{KZ} exhibited numerical results showing a remarkable phenomena related to solitons collision. More precisely, they put in evidence the \emph{elastic} character of the collision between two solitons. By \emph{elastic} we mean that collision keeps the solitons unchanged and does not produce any residual term of positive mass for large times. The unique consequence of the collision is a shift translation on each soliton, depending on their sizes. Next, the work of Lax  \cite{LAX1} developed  a mathematical framework to study these problems. After this, the \emph{inverse scattering method} (we refer e.g. to \cite{AC} and \cite{Miura} for a review) provided explicit formulas for $N$-soliton solutions (Hirota \cite{HIROTA}). Indeed, let $c_1>c_2>0$ and  $\delta_1,\delta_2\in \mathbb{R}$ be arbitrary given numbers. There exists an explicit  solution $U= U_{c_1,c_2}(t,x)$ of \eqref{KdV} which satisfies
\be\label{intcases}
\Big\|U(t,\cdot)- \sum_{j=1}^2 Q_{c_j}(\cdot -c_jt-\delta_j)\Big\|_{H^1(\R)} \mathop{\longrightarrow}_{t\to -\infty} 0,\quad
\Bigl\|U(t,\cdot)- \sum_{j=1}^2 Q_{c_j}(\cdot -c_jt-\delta_j')\Big\|_{H^1(\R)} \mathop{\longrightarrow}_{t\to +\infty} 0,
\ee
for some $\delta_j'$ such that the shifts $\Delta_j=\delta_j'-\delta_j$ depend only on $c_1,c_2$. This solution, called 2-soliton, represents the pure collision of two solitons, with no residual terms before and after the collision. In other words, the collision is \emph{elastic}.

These properties are also valid for the \emph{cubic} mKdV, %where \emph{Miura transform} allows to reduce the problem to the quadratic KdV equation, 
(see \cite{AC}, p. 390) and for the \emph{Gardner} equation (see \cite{Ga, Wa} and references there in). In particular, complete integrability and elastic collisions are still present. Let us recall that for the Gardner equation
\be\label{Geq}
u_t + (u_{xx} + u^2 - \mu u^3)_x =0,
\ee
given $\mu\in \R$, soliton solutions exist for all $c>0$ in the case $\mu<0$, and provided $c<\frac 2{9\mu}$ if $\mu>0$. These solutions are explicit and given by $u(t,x)=Q_{\mu, c}(x-ct)$, where $Q_{\mu, c}$ is the Schwartz function \cite{Wa}
\be\label{SolG}
Q_{\mu, c}(x) :=  \frac{3c}{1+ \rho \cosh (\sqrt{c}x)}; \quad \rho := (1-\frac 92 \mu c )^{1/2}.
\ee
In particular, no soliton-solution exists provided $\mu>0$, and $c>0$ large enough, where the character of the equation becomes \emph{defocusing}.

\medskip

We point out that these techniques are known to be too rigid to be applied to more general models, and have no equivalent for the case of the gKdV equation (\ref{gKdVm1}) with a general nonlinearity. The first purpose of this paper is to confirm this belief under reasonable hypothesis on the nonlinearity:  the collision of two solitons is not elastic in general, except by KdV, mKdV and the Gardner equations. Before establishing our main result we explain the framework where the problem must be posed.

The complete integrability property has been studied in many other differential equations, as NLS, KPI, Benjamin-Ono, etc.; see for example \cite{AC}. In particular, when complete integrability is lost, very little is known. We mention the recent works of Perelman \cite{P}, Holmer, Marzuola and Zworski \cite{HZ,HMZ0, HMZ} and Abou Salem, Fr\"ohlich and Sigal \cite{F} on the problem of 2-soliton collision for the nonlinear Schr\"odinger equation (NLS) under the action of a potential and considering higher velocities.

\smallskip

\subsection{Setting and hypothesis} Let us come back to the general equation (\ref{gKdVm1}). Assume that the nonlinearity $f\in C^3(\R)$.
The Cauchy problem for equation \eqref{gKdVm1} (namely, adding the initial condition $u(t=0)=u_0$) is \emph{locally well-posed} for $u_0\in H^1(\R)$ (see Kenig, Ponce and Vega \cite{KPV}). %In particular, the solution $u(t)\in H^1(\R)$ in a small interval of time containing $t=0$. 

For $H^1(\R)$ solutions, in the general case, unlike the integrable cases, only the following two quantities are conserved by the flow:
\be\label{M}
M(t):= \int_\R u^2(t,x)\,dx = \int_\R u_0^2(x)\, dx =M(0), \quad \hbox{(Mass),}
\ee
and
\bea
E(t) &:= & \frac 12 \int_\R u_x^2(t,x)\,dx - \int_\R F(u(t,x))\,dx  \nonumber \\
& & \quad =  \frac 12 \int_\R (u_0)_x^2(x)\,dx - \int_\R F(u_0(x))\,dx =  E(0), \quad \hbox{(Energy)\label{E}}
\eea
where we have denoted 
\be\label{F}
 F(s) := \int_0^s f(\sigma)\ d\sigma.
\ee

In the case of a pure power $f(s) = s^m$, $m<5$, any $H^1(\R)$ solution is global in time thanks to the conservation of energy (\ref{E}). For $m=5$, solitons are shown to be \emph{unstable} and the Cauchy problem for the corresponding gKdV equation has finite-time blow-up solutions, and see \cite{MMblow} and references there in. It is believed that for $m> 5$ the situation is the same. The origin \emph{grosso modo} of this instability comes from the lack of control for the injection $H^1(\R) \to L^p(\R)$ for $p\geq 5$. Indeed, from the Galiardo-Nirenberg inequality
$$
\int_\R |v|^{p+1} \leq C(p) \big(\int_\R v^2 \big)^{\frac{p+3}{4}}\big(\int_\R v_x^2\big)^{\frac{p-1}{4}},
$$
valid for any $v\in H^1(\R)$, one can see that the energy (\ref{E}) cannot be controlled by the usual $H^1$-norm. Consequently, in this work, we will discard high-order nonlinearities at leading order. Indeed, we will consider nonlinearities $f$ of the form
\be\label{surf}
f\in C^{m+2}(\R), \quad f(u) := u^m + f_1(u), \quad m=2,3,4, \quad \hbox{ with } \quad \lim_{s\to 0}\frac{\abs{f_1(s)}}{\abs{s}^m}=0.  
\ee
Moreover, using stability properties of the solitons, we will have only \emph{global in time} solutions, namely $u(t)\in H^1(\R)$ for all time $t\in \R$.   
  
\smallskip
 
The positive sign leading in front of $f$ (see (\ref{surf})) allows the existence of solitons for  (\ref{gKdV0}) of the form
$$
u(t,x) := Q_c(x-x_0-ct), 
$$
with $c>0$ small enough and $x_0\in \R$, where the function $Q_c$ satisfies the elliptic equation
\be\label{ellipticf}
Q_c'' + f(Q_c) = cQ_c, \quad Q_c\in H^1(\R).
\ee
From Berestycki and Lions \cite{BL} and \eqref{gKdVm1}, it follows that there exists $c_*(f)>0$ (possibly $+\infty$) defined by
$$
c_*(f):=\sup\{ c>0 \text{ such that for all $c'\in (0,c)$, exists $Q_{c'}$ positive
solution of \eqref{ellipticf}} \}.
$$
For all $c>0$, if a solution $Q_c>0$ of \eqref{ellipticf} exists then it can be chosen \emph{even} on $\mathbb{R}$ and \emph{exponentially decreasing} on $\mathbb{R}^+$ (and similarly if $Q_c<0$). Moreover, in \cite{MMas1}, the authors have showed that $0<c<c_*(f)$ is a sufficient condition for \emph{asymptotic stability} in the energy space $H^1$ around the soliton $Q_c$, see also Proposition \ref{ASYMPTOTIC} for details. 

Finally, in this paper, we consider only \emph{nonlinear stable} solitons in the sense of Weinstein
\cite{We2}, i.e. such that
\begin{equation}\label{stab}
\frac d{dc'} \int {Q_{c'}^2(x) dx } \Big|_{c'=c}>0.
\end{equation}
Note that since $m=2,3,4$ in \eqref{surf}, this condition is automatically satisfied for $c>0$ small enough (in the pure power case $f(s)=s^m$, this condition is satisfied for any $c>0$ provided $m<5$, see \cite{We2}).

\subsection{Previous analytic results on 2-soliton collision in non-integrable cases}

As pointed out in \cite{MMcol1}, the problem of describing the collision of two traveling waves or solitons is a general problem for nonlinear PDEs, which is almost completely open, except in the integrable cases described above. On the other hand, these problems have been studied since the 60's from both experimental and numerical points of view. 

\smallskip

We deal with these questions for \eqref{gKdVm1} with a general nonlinearity
$f(u)$ in a particular setting: we consider two positive solitons $Q_{c_1}$, $Q_{c_2}$, $0<c_2<c_1<c_*(f)$, and we assume $c_2$ small compared with $c_1$. 

Under these assumptions, Martel and Merle \cite{MMcol1} considered the collision problem for (\ref{gKdV}) in the quartic  case, $f(s)=s^4$, with one soliton small with respect to the other. They showed that the collision is almost elastic, but inelastic, by showing the nonexistence of pure 2-soliton solution. 

\begin{thm}[Non-existence of a pure 2-soliton solution, quartic case \cite{MMcol1}]\label{PURE4}

Let $f(s):=s^4$ and $0<c:=\frac{c_2}{c_1}<1$. There exists a constant $c_0>0$ such  that if $c<c_0$ then the following holds. Let $u(t)\in H^1(\R)$ be the unique solution of (\ref{gKdVm1}) such that 
\be\label{minusinfty}
\lim_{t\to -\infty} \norm{u(t) - Q_{c_1} (\cdot -c_1 t) -Q_{c_2}(\cdot -c_2 t)}_{H^1(\R)}=0. 
\ee
Then there exist $x_1^+, x_2^+$, $c_1^+ >c_2^+ >0$ and constants $T_0, K>0$ large enough such that
$$
w^+(t,\cdot) := u(t,\cdot) - Q_{c_1^+} (\cdot -x_1^+ -c_1^+ t) -Q_{c_2^+}(\cdot -x_2^+ -c_2^+t)
$$ 
satisfies
\begin{enumerate}
\item Support on the left of solitons. 
$$
\lim_{t\to +\infty} \norm{w^+(t)}_{H^1(x>\frac 1{10} c_2 t)} =0.
$$
\item Parameters perturbation. The limit scaling parameters $c_1^+$ and $c_2^+$ satisfy
$$
\frac 1K c^{\frac{17}6} \leq \frac{c^+_1}{c_1}-1\leq K c^{\frac{11}{6}}, \quad \hbox{ and }
\quad \frac 1K c^{\frac 83} \leq 1 - \frac{c^+_2}{c} \leq K c^{\frac 13}.
$$
In particular, $c_1^+>c_1$ and $c_2^+ <c_2$.
\item Non zero residual term.  For every $t\geq T_0$, the adapted $H^1$-norm of $w^+(t)$ satisfies 
$$
\frac 1K  c_1^{\frac{7}{12}} c^{\frac{17}{12}} \leq \norm{w^+_x (t)}_{L^2(\R)} + \sqrt{ c_1c} \norm{w^+(t)}_{L^2(\R)} \leq K  c_1^{\frac{7}{12}}c^{\frac {11}{12}}.
$$
\end{enumerate}   
\end{thm} 

\begin{rem}
The existence and uniqueness of the solution of (\ref{gKdVm1}) satisfying (\ref{minusinfty}) was proved in \cite{Martel}.
\end{rem}
\begin{rem}
Note that $\norm{Q_{c_2}}_{H^1(\R)} \sim c^{\frac 1{12}} \gg Kc^{\frac{5}{12}} \geq \norm{w^+(t)}_{L^2(\R)}$ for $c$ small. In other words the defect $w^+$ is really small compared with $Q_{c_2}$.
\end{rem}
\medskip

The next question arising from this result is to generalize these results to (\ref{gKdVm1}) under assumption (\ref{surf}). In this case, Martel and Merle \cite{MMcol2} proved that the collision is still stable, giving upper bounds on the residual terms appearing after the collision.  In particular, their result extends the positive part of Theorem \ref{PURE4}.

\begin{thm}[Behavior after collision of a pure 2-soliton solution, \cite{MMcol2}]\label{PURE}~\\
Let $f$ satisfying \eqref{surf}. Let $0<c_2< c_1< c_*(f)$ be such that the positive solution $Q_{c_1}$ of
    \eqref{ellipticf} 
    satisfies \eqref{stab}.  Then there exists $c_0=c_0(c_1)\in (0,c_1)$ such that if $c_2<c_0(c_1)$ then the following holds.
    Let $u(t)$ be the solution of \eqref{gKdV0} satisfying
    \begin{equation}\label{th01}
    \lim_{t\to -\infty} \|u(t)-Q_{c_1}(\cdot - c_1t)-Q_{c_2}(\cdot -c_2 t)\|_{H^1(\R)} =0.
    \end{equation}
    Then, there exist  $\rho_1(t),\rho_2(t)$, $c_1^+>c_2^+>0$ and  $K>0$ such that
    \begin{equation*}
    w^+(t,x):=u(t,x)-Q_{c_1^+}(x-\rho_1(t))-Q_{c_2^+}(x-\rho_2(t)) \quad
        \end{equation*}
        satisfies  $\sup_{t\in \mathbb{R}} \|w^+(t)\|_{H^1(\R)}\leq K c_2^{\frac 1{m-1}}$ and for $q=q_m := \frac 2{m-1}+\frac 14$,
       \begin{equation}\label{th02}
        \lim_{t\to +\infty} \|w^+(t)\|_{H^1(x>\frac 1{10} c_2 t)}=0,\quad
        \limsup_{t\to +\infty} \|w^+(t)\|_{H^1(\R)} \leq K c_2^{q - \frac 12 -\frac 1{100}},
   \end{equation}
   \begin{equation}\label{th02tri}
    \lim_{t\to +\infty}|\rho_1'(t)-c_1^+ | + | \rho_2'(t)-c_2^+|=0.
   \end{equation}
   Moreover, $\lim_{t\to +\infty} E(w^+(t))=:E^+$ and  $\lim_{t\to +\infty} \int_\R (w^+)^2(t)=:M^+$ exist and the following bounds hold
   \begin{equation}\label{th02bis}
   \frac 12 \limsup_{t\to +\infty}  \int_\R ((w^+_x)^2  +   c_2  (w^+)^2)(t)\leq 2 E^+ + c_2 M^+ \leq \liminf_{t\to +\infty} \int_\R ((w^+_x)^2 + 2 c_2 (w^+)^2)(t).
      \end{equation}
Finally, the limit parameters $c^+_1$ and $c^+_2$ satisfy the following bounds      
     \begin{equation}\label{th03}
            \frac 1{K} (2 E^+ + c_2 M^+) \leq \frac{c_1^+}{c_1} -1
            \leq {K}   (2 E^+ + c_2 M^+),
     \end{equation}
     and
        \begin{equation}\label{th04} 
        \frac 1K c_2^{q-\frac 34} (2 E^+ + c_1 M^+) 
          \leq 1-\frac {c_2^+}{c_2} 
            \leq {K} c_2^{q-\frac 34} (2 E^+ + c_1 M^+) .
   \end{equation}
\end{thm} 

\begin{rem}
In Theorem \ref{PURE}, if $c_1^+=c_1$ and $c_{2}^+=c_{2}$ (or equivalently $E^+ =M^+ =0$),  then the solution $u(t)$ is a \emph{pure 2-soliton solution and the collision is elastic}.  
\end{rem}

In \cite{MMcol2}, the question of whether the collision is elastic or inelastic in the general case --and thus the nonexistence of pure 2-soliton solutions-- was left open, see \cite{MMcol2}, Remark 1. More precisely, the authors conjectured a \emph{classification result} concerning the nonlinearities $f(s)$ allowing small stable solitons. This affirmation asserts that under reasonable stability properties, \emph{the unique nonlinearities for which any 2-soliton collision is pure are the integrable cases, $f(s)=s^2$, $f(s)=s^3$ and a linear combination of both nonlinearities.} Theorem \ref{PURE4} from \cite{MMcol1} was the first step in this direction. By extending some techniques from \cite{MMcol1}, \cite{MMcol2} and developping new computations, we are able to provide a satisfactory answer to this open question. 

\subsection{Main results}

Consider the framework introduced in Theorem \ref{PURE}. In addition to this result, we have the following

\begin{thm}[Non-existence of pure 2-soliton solution, general case]\label{PUREgen0}
Let $f$ be as in (\ref{surf}), with $m=2$ or $3$, and 
\be\label{fp}
f\in C^{p+1}(\R), \quad f^{(p)}(0)\neq 0 \quad \hbox{for some}\quad  p\geq 4.
\ee
For $0<c_2\ll c_1\ll1$ equation (\ref{gKdVm1}) has no pure 2-soliton solution of sizes $c_1, c_2$. In particular Theorem \ref{PURE} holds with $c_1^+>c_1$ and $c_2^+<c_2$.
\end{thm} 

\begin{rem}
The nonzero condition $f^{(p)}(0)\neq 0$ for some $p\geq 4$ rules out the integrable cases $f(s)=s^m$, $m=2$ or $3$ and the Gardner nonlinearity $f(s) = s^2 - \mu s^3$.
\end{rem}
\begin{rem}\label{1.3}
We do not treat the degenerate cases $f(s) = s^m + f_1(s)$, for $f_1(s)\neq 0$ but $f_1^{(p)}(0) = 0$ for all $p\geq 4$. These cases seem to be not physically relevant.
\end{rem}
\begin{rem}
The result of Theorem \ref{PUREgen0} in the quartic case $m=4$ follows directly from the proof of Theorem \ref{PURE4} in \cite{MMcol1} together with the techniques used in the present paper. This remark and Theorem \ref{PUREgen0} allow to classify the nonlinearities for which 2-soliton collision is elastic. In particular, with the restriction mentioned in Remark \ref{1.3}, we obtain that  \emph{pure 2-soliton solutions are present for any pair of solitons with different velocities if and only if $f$ corresponds to the integrable cases, $f(s)=s^2, s^3$ and linear combinations}. We recall that for $m\geq 5$ in (\ref{surf}), solitons have been shown to be unstable (see \cite{BSS}).  It is believed that collision may produce blow-up solutions in finite time. 
\end{rem}

\medskip

Theorem \ref{PUREgen0} is a consequence of the following reduction of the problem. Let $p$ be the smallest integer greater or equal than $4$ satisfying (\ref{fp}). Let $c_1>0$ small. Consider the transformation
$$
\tilde u (t,x) :=c_1^{-\frac 1{m-1}} u(c_1^{-\frac 32} t, c_1^{-\frac 12} x),
$$
which maps $Q_{c_1}(x-c_1t)$ to $Q(x-t)$ and $Q_{c_2}(x-c_2 t)$ to $Q_c(x-ct)$, with $c:=\frac{c_2}{c_1}$. If $u=u(t,x)$ is solution of (\ref{gKdVm1}) then $\tilde u$ is solution of the equation
$$
\begin{cases}
\tilde u_t + (\tilde u_{xx} + \tilde f (\tilde u))_x=0,\\
\hbox{with }Ê\quad \tilde f(\tilde u):= \tilde u^m + \tilde f_1(\tilde u), \quad \tilde f_1(\tilde u) := c_1^{-\frac m{m-1}} f_1(c_1^{\frac 1{m-1}} \tilde u). 
\end{cases}
$$
Note that $ \tilde f_1$ satisfies (\ref{surf}). Then, for the case $m=3$, $\tilde f_1$ can be expanded as
\bee
\tilde f_1 (\tilde u) & = & c_1^{-\frac 3{2}}\Big[ \frac 1{p!} f_1^{(p)}(0) (c_1^{\frac 1{2}}\tilde u)^p + O((c_1^{\frac 1{2}}\tilde u)^{p+1})\Big] \\
& =: & \ve \tilde u^p + |\ve|^{1+\frac{1}{p-3}} \hat f_1(\tilde u),
\eee
where
$\ve := \frac{1}{p!} f_1^{(p)}(0) c_1^{\frac{p-3}{2}}$ is small and $\hat f_1$ satisfies the decay relation $\lim_{s\to 0}|s|^{-p}\hat f_1(s) =0$. 

For the quadratic case, we need more care \emph{because of the Gardner nonlinearity}. We have
\bee
\tilde f_1 (\tilde u) & = & c_1^{-2}\Big[  \frac 16 f_1^{(3)}(0)(c_1 \tilde u)^3 +  \frac 1{p!} f_1^{(p)}(0) (c_1 \tilde u)^p + O((c_1\tilde u)^{p+1})\Big] \\
& =&  \frac 16 f_1^{(3)}(0) c_1 \tilde u^3 +  \frac 1{p!} f_1^{(p)}(0) c_1^{p-2} \tilde u^p + O(c_1^{p-1}\tilde u^{p+1})\\
& =: & \mu(\ve) \tilde u^3 + \ve \tilde u^p + |\ve|^{1+ \frac 1{p-2}} \hat f_1(\tilde u),
\eee
where $\ve := \frac{1}{p!} f_1^{(p)}(0) c_1^{p-2}\neq 0$ by hypothesis, and $\mu(\ve) := \frac{1}{6} f_1^{(3)}(0)c_1 = \hat \mu \ve^{\frac 1{p-2}}$, $\hat \mu \in \R$. Here both $\ve$ and $\mu$ are small (depending on $c_1$) and $\hat f_1$ satisfies the decay relation $\lim_{s\to 0}|s|^{-p}\hat f_1(s) =0$. Note that in this framework, the quadratic case can be seen as a particular case of the Gardner nonlinearity, for which $\hat \mu=0$.

Finally, we drop the tilde on $\tilde u$ and $\tilde f$ and the hat on $\hat f_1$.  We are now reduced to the \emph{$\ve$-dependent} equation
\be\label{gKdV0}
u_t + (u_{xx} + f(u))_x =0, 
\ee
where $\mu (\ve) = \hat \mu \ve^{\frac 1{p-2}}$,  $\hat \mu \in \R$,
\be\label{gKdV0b}
f=f_\ve \in C^{p+1}(\R), \;  f(u)= \begin{cases} u^2+ \mu(\ve) u^3 + \ve u^p  + \ve^{1+ \frac 1{p-2}}f_1(u), \; m=2,\\
u^3+\ve u^p  + \ve^{1+\frac 1{p-3}}f_1(u), \; m=3, \end{cases} \; \lim_{s\to 0} \frac{f_1(s)}{\abs{s}^p} =0,
\ee
and $\ve$ is small, and $p\geq 4$. For notational commodity we will skip the $\ve$-dependence on the functions considered along this work, except in some computations performed in Appendix \ref{D}. Lastly, note that for $\ve$ small $Q$ and $Q_c$ satisfy (\ref{stab}), see also Remark \ref{stabeps}. 

In this framework, we now claim the main result of this paper: 

\begin{thm}[Non-existence of pure 2-soliton solution, general case]\label{PUREgen}~\\
Suppose $m=2,3$ and $f$ satisfying (\ref{gKdV0b}) for $p\geq 4$. There exists a  constant $\ve_{0}>0$ such that if 
\be\label{Cond0} %0< c^{\frac 2{m-1} + \frac 14}<\ve
0<\abs{\ve}< \ve_0, \quad \hbox{ and } \quad  0 < c \leq  \abs{\ve}^{m-1 + \frac 1{25}},  
\ee  
then the following holds. Let $u(t)$ be solution of (\ref{gKdV0})  satisfying 
\be\label{minusinfty2}
    \lim_{t\to -\infty} \|u(t)-Q (\cdot - t)-Q_{c}(\cdot -c t)\|_{H^1(\R)} =0,
\ee 
there exist $K,T_0>0$ such that
\begin{enumerate}

\item Non zero residual term.  There exist $\rho_1(t)$, $\rho_2(t)$, $c_1^+>c_2^+ >0$ such that %the residual term 
$$
w^+(t) := u(t,x) -Q_{c_1^+}(x-\rho_1(t)) -Q_{c_2^+}(x-\rho_2(t))
$$
 satisfies, for every $t\geq T_0$,  
\begin{equation}\label{th023}
        \lim_{t\to +\infty} \|w^+(t)\|_{H^1(x>\frac 1{10} c t)}=0, %\quad \hbox{(support on the left of solitons),}
\end{equation}
and for $q = \frac 2{m-1} + \frac 14$,%  the adapted $H^1$-norm of $w^+(t)$ satisfies,
\be\label{restgen}
\frac 1K \abs{\ve} c^{q +\frac 12}  \leq \norm{w^+_x (t)}_{L^2(\R)} + \sqrt{c} \norm{w^+(t)}_{L^2(\R)} \leq K \abs{\ve} c^{q}.
\ee
\item The asymptotic scaling parameters $c_1^+$ and $c_2^+$ satisfy
\bea\label{scagen}
& & \frac 1K \ve^2 c^{2q + 1} \leq c^+_1-1\leq K \ve^2 c^{2q}, \, \hbox{ and } \nonumber \\
& &  \frac 1K \ve^2 c^{3q +\frac 14} \leq 1 - \frac{c^+_2}{c} \leq K \ve^2 c^{3q -\frac 74}.
\eea
%In particular, $c_1^+>1$, $c_2^+ <c$ and $\norm{w^+ (t)}_{L^2(\R)} \leq K \abs{\ve} c^{q-\frac 12} \ll \norm{Q_c}_{L^2(\R)} \sim c^{\frac 12q -\frac 38}$. 
\end{enumerate}        
\end{thm} 

\medskip

Before sketching the proof of this Theorem, some remarks are in order.

\begin{rem}[Comments on the assumptions]
In the present paper, as in \cite{MMcol1,MMcol2}, we study the collision of two solitons $Q_{c_1}$ and $Q_{c_2}$. The assumption $c_2$ small allows to linearize in $c_2$, and then to reduce the non existence of a pure 2-soliton solution  to the computation of a \emph{coefficient} depending only on $c_1$. For general $f$, as in (\ref{surf}), and general $c_1>0$, it is an open question to compute this coefficient, see Remark \ref{smallep} (for $p=4$, a special algebraic structure allowed to compute this coefficient, see \cite{MMcol1}).   

According to this, we compute the asymptotics of this coefficient as $c_1$ is small (or equivalent, $\ve$ is small), see Appendix \ref{D}. This is the only place where $\ve$ small is needed. This asymptotic allows us to conclude under the additional restriction $0<c<\abs{\ve}^{m-1+ \frac 1{25}}$ (see Proposition \ref{cor:1} and estimate (\ref{cond1})). The exponent $\frac 1{25}$ has no special meaning, and can be taken as small as we want, as long as $c$ is taken even smaller.

Two open questions then arise:
\begin{enumerate}
\item Can we relax in (\ref{Cond0}) the second condition on $c$? 
\item For general $f$, do there exist special values of $c_1$ for which the coefficient is zero? The residue from the collision would then be of smaller order in $c_2$.
\end{enumerate} 
%For further discussion on these non-perturbative regimes see also \cite{Tao}.  
\end{rem}

\begin{rem}
For $m=2$, the smoothness condition (\ref{gKdV0b}) allows nonlinearities of type $f(s) =s^2 + \nu \ve^{p-2} s^{p}$, with $p\geq 4$ (possibly non integer). For $m=3$, the same conclusion follows for nonlinearities of the type $f(s) = s^3 + \ve s^p$, $p=4$ and $p\geq 5$ (possibly non integer). See also Appendix \ref{sec:B} and final remarks in Appendix \ref{D}.
\end{rem}

\begin{rem}
Although this theorem asserts that collision is indeed inelastic, near-elastic, the appearance of smaller solitons on the left of the solitons is not discarded by our proof, and (\ref{th023}). However, we believe that, at least under the condition of Theorem \ref{PUREgen}, there are no such small solitons. 
\end{rem}

\subsection{Sketch of the proof}\label{sop}

Our proof will follow closely the approach described by Martel, Merle and Mizumachi \cite{MMcol1, MMcol2, MMMcol}. The argument is as follows: we consider the solution (unique, see \cite{Martel}) $u(t)$ of (\ref{gKdV0}) satisfying (\ref{th01}) at time $t\sim -\infty$. Then, we separate the analysis among three different time intervals: $t\ll -c^{-\frac 12}$, $\abs{t} \leq c^{-\frac 12}$ and $c^{-\frac 12} \ll t$. On each interval the solution possesses a specific behavior which we briefly describe:      

\begin{enumerate}
\item  ($t\ll -c^{-\frac 12}$). In this interval of time we prove that $u(t)$ remains close to a 2-soliton solution with no changes on scaling and shift parameters. This result is possible for negative long enough times, such that both solitons are still far from each other, and is a consequence of  \cite{Martel}. 

\item ($\abs{t}\leq c^{-\frac 12}$). This is the interval where solitons collision leads the dynamic of $u(t)$. The novelty in the method is the construction of an \emph{approximate solution} of (\ref{gKdV0}) with high order of accuracy such that $(a)$ at time $t\sim -c^{-\frac 12}$ this solution is close to a 2-soliton solution and therefore to $u(t)$, $(b)$ it describes the 2-soliton collision in this interval, $(c)$ at time $t\sim c^{-\frac 12}$, when solitons are sufficiently separated, it possesses an extra, nonzero, residual term product of the collision, and characterized by a number $d(\ve)\neq 0$ (cf. (\ref{diff0})-(\ref{diff2})), and $(d)$ it is possible to extend the solution $u(t)$ to the whole interval $[-c^{\frac 12}, c^{\frac 12}]$ being still close to our approximate solution, uniformly on time, modulo modulation on a translation parameter. This property confirms that our Ansatz is indeed the correct approximate solution describing the collision. 

\item ($t\gg c^{-\frac 12}$)  Here some \emph{stability} properties (see Proposition \ref{ASYMPTOTIC}) will be used to establish the convergence of the solution $u(t)$ to a 2-soliton solution with modified parameters. Moreover, by using a \emph{monotony} argument, it will be possible to show that the residue appearing after the collision at time $t\sim c^{-\frac 12}$ is still present at infinity. This gives the conclusion of the Theorem. 
\end{enumerate}

\medskip

The plan of this paper is as follows. In Section 2 we construct the aforementioned approximate solution and compute the error term produced in terms of a set of linear problems. Then we solve such linear systems and finally we give the first basic estimates concerning this solution. We finally prove that it is indeed close to a 2-soliton solution. In section 3 we construct an actual solution $u$ close to the approximate solution for small times, and state some stability results to study the long time behavior of the solution $u$. Finally, in section 4, we prove Theorem using above results.

%%%%%%%%%%%%%%%%%%%%%%%%%%%%%%
%%%%%%%%%%%%%%%%%%%%%%%%%%%%%%

\bigskip

\section{Construction of an approximate 2-soliton solution}\label{sec:2}

%\subsection{Warnings}
The objective of this section is to construct an approximate solution of
the gKdV equation (\ref{gKdV0}), which will precisely describe the collision of two solitons. 
Hereafter, we assume the hypothesis of Theorem \ref{PUREgen}. We suppose both solitons are positive (the negative case, for $m=3$, can be treated in the same way).

Secondly, note that $Q$ and $Q_c$ have \emph{velocity} (and size) 1 and $c$ respectively;  so that working with $u(t,x+t)$ instead of $u(t,x)$ we can assume that the great soliton $Q$ is fixed at $x=0$ and the small soliton has velocity $c-1<0$. Of course, $v(t,x):=u(t, x+t)$ satisfies now the \emph{translated} equation
\be\label{gKdV}
v_t + (v_{xx}-v +f(v))_x =0 \quad \hbox{ on }\R_t\times \R_x.
\ee

Finally, denote
\be\label{Tc}
T_c:=c^{-\frac 12 -\frac 1{100}}>0.
\ee
This quantity can be understood as the time of interaction between the two solitons. The exponent $\frac 1{100}$ can be replaced by any small positive number without relevant modifications. 

\medskip

The following result deals with the problem of describing the collision in the interval of time $[-T_c, T_c]$:

\begin{prop}[Construction of an approximate solution of the gKdV equation]\label{Constr}

Let $m=2,3$ and $f$ as in (\ref{gKdV0b}). There exist constants $c_0=c_0(f)>0$ and $K_0=K_0(f)$ such that for all $0<c <c_0$ there exists a function $\tilde u= \tilde u_{1,c}(t,x)$  such that the following hold:

\begin{enumerate}
\item Approximate solution on $[-T_c, T_c]$. For all $t\in [-T_c, T_c]$, 
$$
\norm{\tilde u_t + (\tilde u_{xx} -\tilde u + f(\tilde u))_x}_{H^2(\R)} \leq K_0 c^{\frac 3{m-1} +\frac 34}.
$$

\item Closeness to the sum of two solitons: For all time $t\in [-T_c, T_c]$, the function $\tilde u$ belongs to $H^1(\R)$ and satisfies
$$
\norm{\tilde u(t) -Q(x-\al)- Q_c(x+(1-c)t) }_{H^1(\R)} \leq K_0c^{\frac 1{m-1}},
$$
where $\alpha=\al(t,x)$ is a smooth bounded function, to be defined below, see (\ref{defALPHA}).
\end{enumerate}
  
\end{prop}

\begin{rem}
The proof of this proposition requires several steps, starting in Subsection \ref{sec:2-1} to finally ending in Subsection \ref{sec:3.7}, Proposition \ref{lem:z}. However, the proof is intuitively clear to describe: our approximate solution will consists of a \emph{linear combination} of a \emph{nonlinear basis} well behaved under the gKdV flow, together a variable decomposition resembling the classical separation of variables from second order linear PDEs. 
This description was first introduced by Martel and Merle \cite{MMcol1}, \cite{MMcol2}. 
\end{rem}

First of all we explain how the approximate solution is composed. We follow \cite{MMcol1}.  

\subsection{Decomposition of the approximate solution}\label{sec:2-1}

We look for $\tilde u(t,x)$, the approximate solution for (\ref{gKdV}),  carring out a specific structure. We first introduce a set of indices, depending on the cases we deal with. Let
$$
\Sigma_2 := \{ (k,l) = (1,0),\, (1,1),\, (2,0),\, (2,1),\, (1,2),\, (3,0)\}, 
$$
for the quadratic case ($m=2$), and
$$
\Sigma_3 := \{ (k,l) = (1,0),\, (1,1),\, (2,0),\, (2,1),\, (3,0),\, (4,0)\},
$$
for the cubic one ($m=3$).% We will see in (\ref{defv})-(\ref{defW}) that these sets represent the size of the correction term to a 2-soliton solution.

We recall now an order relation for indices $(k,l), (k',l')\in \Sigma_m$ introduced in \cite{MMcol1}. We say that 
\be\label{order}
(k',l')< (k,l) \quad \hbox{ if and only if} \quad \begin{cases} k'< k \hbox{ and } l'\leq l, \hbox{ or } \\  k'\leq k \hbox{ and } l'<l .\end{cases} 
\ee

\smallskip

We set two variables denoting the position of each soliton. For the small soliton, let
\begin{equation*}
    y_{c}:=x+(1-c)t  \quad \text{and} \quad R_c(t,x):= Q_c(y_c),
\end{equation*} 
and for the great soliton,
\begin{equation*}
    y:=x-\alpha(y_{c}) \quad \text{and} \quad     R(t,x): =Q(y),
\end{equation*}
where for $(a_{k,l})_{(k,l)\in \Sigma_m}$,  
\begin{equation}\label{defALPHA} 
    \alpha(s):=\int_{0}^s \beta(s') ds',\quad     
    \beta(s):=\sum_{(k,l)\in \Sigma_m} a_{k,l} \, c^l Q_{c}^k(s).
\end{equation}

The correction term $\alpha$ is intended to describe the shift on the position of the great soliton. Note that $\alpha$ might be nonzero even in the integrable case, see (\ref{intcases}). Moreover, in the quartic case $m=4$, $\ve=0$, one has $\abs{\al}\to +\infty$ as $c\to 0$, see \cite{MMas1}. Along this work $\alpha$ will be a bounded function, uniformly on $c$.  

\smallskip

The form of $\tilde u(t,x)$ is, as it should be expected, the sum of the two soliton plus a correction term:
\begin{equation}\label{defv} 
    \tilde u(t,x) :=Q(y)+Q_{c}(y_{c})+W(t,x),
\end{equation}
\begin{equation}\label{defW}
    W(t,x):=\sum_{(k,l)\in \Sigma_m} 
        c^l\left(Q_{c}^k(y_{c}) A_{k,l}(y)+(Q_{c}^k)'(y_{c}) B_{k,l}(y)\right),
\end{equation}
where $ a_{k,l} $, $A_{k,l}$, $B_{k,l}$ are unknowns to be determined.

The motivation in \cite{MMcol1} 
for choosing $W$ of the form \eqref{defW} is precisely the closeness of the family of functions
\begin{equation}\label{ily}
    \left\{c^l Q_c^k,\ c^l (Q_c^k)', \ k\geq 1,\ l \geq 0 \right\}
\end{equation}
under multiplication and differentiation, due to the specific form of the equation  of $Q_c$ (see Lemma 2.1 in \cite{MMcol1}). In the case of equation \eqref{gKdV},
for a general nonlinearity this structure is preserved up to a lower order term (see Lemma \ref{surQc2}). 

We want to measure the size of the error produced by inserting $\tilde u$ as defined in (\ref{defv})-(\ref{defW}) in the equation (\ref{gKdV}). For this, let 
\begin{equation}\label{2.2bis}
S[\tilde u](t,x) :=\tilde u_t + (\tilde u_{xx} - \tilde u + f(\tilde u))_x.
\end{equation}
%It is clear that for an actual  solution $u=u(t,x)$ of (\ref{gKdV}) we have $S[u]\equiv 0$.

Our first result in the above direction is the following

\begin{prop}[Decomposition of $S(\tilde u)$]\label{prop:decomp}
%Assume that $f$ is of class $C^{m+2}(\R)$.  
Let
\begin{equation}\label{defLy}
    \mathcal{L} w := - w_{yy} + w -f'(Q) w.
\end{equation}
 Then,
    \begin{align*}
            S[\tilde u](t,x) & = 
        \sum_{(k,l)\in \Sigma_m}
        c^l Q_c^k(y_c)    \Big[a_{k,l} (-3 Q+2 f(Q))'(y)    -(\mathcal{L} A_{k,l})'(y) +F_{k,l}(y)  \Big]
        \\& \quad
        + \sum_{(k,l)\in \Sigma_m}
        c^l (Q_c^k)'(y_c)    
        \Big[a_{k,l} (-3 Q'')(y) + \left(3A_{k,l}'' +f'(Q) A_{k,l}\right)(y)    - (\mathcal{L} B_{k,l})'(y) + G_{k,l}(y) \Big]
        \\& \quad +  \mathcal{E}(t,x)
    \end{align*}
    where 
    $F_{k,l} $, $G_{k,l} $ and $\mathcal{E}$ satisfy, for any
    $(k,l)\in \Sigma_m$,
    \begin{itemize}
        \item[{\rm (i)}] Dependence property of $F_{k,l}$ and $G_{k,l}$:
        The expressions of $F_{k,l}$ and $G_{k,l}$ depend only on $(a_{k',l'})$, $(A_{k',l'})$, $(B_{k',l'})$ for $(k',l') < (k,l)$.
        \item[{\rm (ii)}] Parity property of $F_{k,l}$ and $G_{k,l}$:   Assume that for any $(k',l')$ such that $(k',l')< (k,l)$
         $A_{k',l'}$ is even and    $B_{k',l'}$ is odd, then 
        $F_{k,l} $ is odd and  $G_{k,l}$ is  even.

 Moreover,
$F_{1,0}=(f'(Q))'$ and $G_{1,0} = f'(Q)$, and higher order terms are given in Appendix \ref{sec:B}.
	\item[{\rm (iii)}] Estimate on $\mathcal{E}$: Assume both $(A_{k,l})$ and $(B_{k,l})$ bounded, and $(A_{k,l}'), (B_{k,l}') \in \mathcal Y$ for $(k,l)\in \Sigma_m$. Then there exists $\kappa>0$ such that for all $j=0,1,2$, and for every $(t,x)\in [-T_c,T_c]\times \R$,  
\begin{equation*}%\label{onE}
	|\partial_x^j \mathcal{E}(t,x)|\leq
\kappa c^{m-1} Q_c(y_c).
\end{equation*}
    \end{itemize}
\end{prop}

\medskip

\begin{rem}
Note that $(\mathcal L w)_y$, as defined in (\ref{defLy}), represents the linear operator associated to the gKdV equation (\ref{gKdV}). Thus, the expression for $S[\tilde u ]$ above stated can be seen as a generalization of the linearized gKdV equation, with the addition of some correction terms.
\end{rem}

\begin{proof}
We postpone the proof of the Proposition \ref{prop:decomp}, merely calculative, to Appendix \ref{sec:B}. We note that this Proposition has been already stated in \cite{MMcol2}, but here we will need an improved version, describing explicitly every term $F_{k,l}, G_{k,l}$ up to a fixed high order. For the details, see Appendix \ref{sec:B}.
\end{proof}

Note that if we want to improve the approximation $\tilde u$, the unknown functions $A_{k,l}$ and $B_{k,l}$ for a fixed $(k,l)$ must be chosen satisfying a sort of modified linear gKdV system where the source terms are composed of preceding, well-known, $A_{k',l'}$ and $B_{k',l'}$ functions. Indeed, if we choose (formally) $A_{k,l}$ and $B_{k,l}$ such that for any $(k,l)\in \Sigma_m$

\begin{equation*}
  (\Omega_{k,l})\quad 
  \left\{
    \begin{aligned}
& (\mathcal L A_{k,l})'+a_{k,l}(3Q-2f(Q))'=F_{k,l},\\
& (\mathcal LB_{k,l})'+ 3a_{k,l}Q'' -3A_{k,l}'' - f'(Q)A_{k,l}=G_{k,l},
    \end{aligned}
\right.
\end{equation*}
then the error term will be reduced to the quantity
$$
S[\tilde u ] = \mathcal{E}(t,x).
$$

Of course the solvability theory for the linear systems $(\Omega_{k,l})$ and the measure of this error term must be stated in a rigorous form. This will be established in the following section.

%%%%%%%%%%%%%%%%%%%%%%%%%%%%%%%%%%%%
%%%%%%%%%%%%%%%%%%%%%%%%%%%%%%%%%%%%

\subsection{Resolution of linear systems $(\Omega_{k,l})$ }

First, we recall some preliminary notation and results from \cite{MMcol1}.
We denote by $\mathcal{Y}$ the set of $C^\infty$ functions $f$ such that
\begin{equation}\label{eq:Y}
\forall j\in \mathbb{N},~\exists K_j,~r_j>0,~Ê\forall x \in \mathbb{R},\quad
|f^{(j)}(x)|\leq K_j (1+|x|)^{r_j} e^{-|x|}.
\end{equation}

We recall some well-known results concerning a \emph{resonance} function and the operator $\mathcal{L}$.

\begin{Cl}[\cite{MMcol2}]\label{surphi}
    The function $\varphi(x)=-\frac {Q'(x)}{Q(x)}$ is odd and    satisfies:
    \begin{itemize}
        \item[{\rm (i)}] $\lim_{x\rightarrow -\infty} \varphi(x)=-1$; $\lim_{x\rightarrow +\infty} \varphi(x)=1$;
        \item[{\rm (ii)}] $\forall x\in \mathbb{R}$, $|\varphi'(x)|+|\varphi''(x)|+|\varphi^{(3)}(x)|\leq C e^{-|x|}$.
        \item[{\rm (iii)}] $\varphi'\in \mathcal{Y}$,  $(1-\varphi^2) \in \mathcal{Y}$.
               \end{itemize} 
\end{Cl}

\begin{lem}[Properties of $\mathcal{L}$, see \cite{MMcol2}]\label{surL}  The operator $\mathcal{L}$ defined in $L^2(\mathbb{R})$ by \eqref{defLy}  has domain $H^2(\R)$, is self-adjoint and satisfies the following properties:
    \begin{itemize}
        \item[{\rm (i)}] There exist a unique $\lambda_0>0$, $  \chi_0 \in H^1(\mathbb{R})$,
  $ \chi_0 >0$ such that
  $\mathcal{L}   \chi_0 =-\lambda_0  \chi_0 $.
                \item[{\rm (ii)}] The kernel of $\mathcal{L}$ is 
        $\{\lambda Q', \lambda \in \mathbb{R}\}$.
        Let  $\Lambda Q :=\frac d{dc} {Q_{c}}_{| c=1}$, then
         $\mathcal{L} (\Lambda Q)=-Q $.
        \item[{\rm (iii)}] (Inverse) For all   $h \in L^2(\mathbb{R})$ such that $\int_\R h Q'=0$, 
        there exists a unique $\widetilde h \in H^2(\mathbb{R})$  such that $\int_\R \widetilde hQ'=0$ and $\mathcal{L} \widetilde h=h$; moreover,
        if $h$ is even (resp. odd), then $\widetilde h$ is even (resp. odd).
        \item[{\rm (iv)}] For $h\in H^2(\mathbb{R})$,  $\mathcal{L} h \in \mathcal{Y}$ implies $h\in \mathcal{Y}$.
	\item[{\rm (v)}] (Coercivity) If $\frac d{d\widetilde c} {\int_\R Q_{\widetilde c}^2}_{|\widetilde c=c}>0$
then there exists  $\lambda_c>0$ such that if
$$
\int_\R  w Q_c=\int_\R w Q_c'=0\quad \hbox{ then } \quad
\int_\R (w_x^2+c w^2 -f'(Q_c) w^2) \geq \lambda_c \int_\R w^2.
$$
  \item[{\rm (vi)}]  There exist unique even solutions $P$ and  $\bar P$ of the ordinary differential equations
\be\label{P}
\mathcal{L} P = 3Q'' + f'(Q)Q, \quad P\in \mathcal{Y},
\ee
\be\label{barp}
\mathcal{L} {\bar P} = f'(Q), \quad \bar P\in \mathcal{Y}.
\ee
Moreover, $P:=-(xQ' +\Lambda Q + Q)$.

    \end{itemize}
\end{lem}

\begin{rem}
Item (vi) from above Lemma is new; the proof follows directly from (ii), (iii) and (iv). On the other hand, for general nonlinearities $\bar P$ is not explicit. 
\end{rem}

%%%%%%%%%%%%%%%%%%%%%%%%
%%%%%%%%%%%%%%%%%%%%%%%%

\subsubsection{Existence theory for a model problem}

We recall that linear systems $(\Omega_{k,l})$ are very similar and then proving existence reduces to prove the result for a model problem. This idea comes from \cite{MMcol1}, but we will need a simplified version, from \cite{MMMcol}. 

\begin{prop}[Existence for a model problem, see \cite{MMcol2}]\label{prop:ex}
Let $F\in \mathcal{Y}$, odd, and $G\in \mathcal{Y}$, even.
Let $\gamma, \kappa \in \R$. Then, there exist $a,b\in \R$,   
$\widetilde A\in \mathcal{Y}$ even, and $\widetilde B\in \mathcal{Y}$ odd,
such that
$$
A=\widetilde A+ \gamma,\quad\hbox{ andÊ}\quad  B=\widetilde B+ b \varphi + \kappa Q'
$$
satisfy
\begin{equation*}
  (\Omega)\quad 
  \left\{
    \begin{aligned}
& (\mathcal L A)'+a(3Q-2f(Q))'=F,\\
& (\mathcal LB)'+ 3aQ'' -3A'' - f'(Q)A=G
    \end{aligned}
\right.
\end{equation*}
Moreover, %in this case we can recover $a$ and $b$ by the formulae
\be\label{a}
a = \frac{-1}{\int_\R \Lambda Q Q} \Big\{  \ga\int_\R P +\int_\R GQ - \int_\R  F\int_0^x P \Big\} 
\ee
and
\be\label{b}
b = \frac 12\Big[ \ga\int_\R \bar P + a\int_\R \Lambda Q -\int_\R F\int_0^x \bar P + \int_\R G \Big] .
\ee
\end{prop}

\begin{proof}
We give a sketch of the proof for the sake of completeness. The original result comes from \cite{MMcol2}, 
and here it is even simpler since we deal only with $F,G\in \mathcal{Y}$.

Set $
A:=\widetilde A+ \gamma,$ $ B:=\widetilde B+ b \varphi
$, where $\gamma$ is given, while $b$ is a parameter to be found.
Since $(\mathcal{L} 1)'=(1-f(Q))'=-(f(Q))' $, we obtain the following system for $\widetilde A$, $\widetilde B$:
\begin{equation*}
    \left\{
    \begin{aligned}
& (\mathcal L\widetilde A)'+a(3Q-2f(Q))'=F+\gamma (f(Q))',\\
& (\mathcal L\widetilde B)'+ 3aQ'' -3\widetilde A''-  f'(Q) \widetilde A=G +\gamma f'(Q) -b (\mathcal L\varphi)'.
    \end{aligned}
\right.
\end{equation*}
Note that $F\in \mathcal{Y}$ is odd, therefore $\mathcal{H}(x)=\int_{-\infty}^x F(z) dz
+  \gamma f(Q)$ belong to $\mathcal{Y}$ and is even.
By integration of the first line, we are reduced to solve
\begin{equation*}
    \left\{
    \begin{aligned}
&  \mathcal L\widetilde A +a(3Q-2f(Q)) =\mathcal{H},\\
& (\mathcal L\widetilde B)'+ 3aQ'' -3\widetilde A''- f'(Q)\widetilde A=G +\gamma f'(Q) -b (\mathcal L\varphi)'.
    \end{aligned}
\right.
\end{equation*}
Since $\int_\R \mathcal{H} Q'=0$ (by parity) and $\mathcal{H}\in \mathcal{Y}$, 
by Lemma \ref{surL}, there exists $\overline H \in \mathcal{Y}$, even, such that
$\mathcal L\overline H= \mathcal{H}$. 

Define $\hat P$ to be the unique even solution of
$$
\mathcal{L} \hat P = 3Q-2f(Q), \quad \hat P \in \mathcal{Y}.
$$
Indeed, $\hat P$ has an explicit formula
\be\label{hatp}
\hat P = -(xQ' +\Lambda Q), \quad \hbox{ with } \mathcal L (\Lambda Q) = -Q.
\ee
It follows that $\widetilde A:=-a \hat P + \overline H$
is even, belongs to $\mathcal{Y}$ and solves the first line of the previous system.
Note that at this stage, the parameters $a$ and $b$ are still free.

Now, we only need to find $\widetilde  B\in \mathcal{Y}$, odd, such that
$(\mathcal L\widetilde  B)'=-a Z_0 + D - b (\mathcal L\varphi)',$
where 
$$
D := 3 \overline{H}'' + f'(Q) \overline H + G + \gamma f'(Q) \in \mathcal{Y}, \hbox{ even, }
Z_0 := 3Q'' + 3\hat P'' + f'(Q)\hat P \in \mathcal{Y}, \hbox{ even}.
$$
Let 
$$
E:=\int_0^x (D-aZ_0)(z) dz - b \mathcal L \varphi.
$$
This function \emph{a priori} is in $L^{\infty}(\R)$, independent of $a,b$.
\begin{Cl}\label{cl:ndg}
There exist numbers $a$ and $b$ such that $E\in \mathcal{Y}$ and  $\displaystyle{\int_\R  E Q'=0}$.
\end{Cl}
Assuming Claim \ref{cl:ndg}, we fix  $a,b$ so that $E\in \mathcal{Y}$ and  $\int_\R  E Q'=0$.
It follows from Lemma \ref{surL} that there exists $\widetilde B\in \mathcal{Y}$, odd,
such that $\mathcal L\widetilde  B=E$. The final solution is then given by $A:=\widetilde A + \gamma$
and $B:=\widetilde B + b \varphi + \kappa Q'$, where $\kappa$ is a free parameter, because $\mathcal{L} Q' =0$ (see Lemma \ref{surL}, (ii)).

\begin{proof}[Proof of Claim \ref{cl:ndg}]

First, we  check a sort of non-degeneracy condition, namely that $\int_\R Z_0Q\neq 0$. Indeed,
by (\ref{P})
\begin{align*}
\int_\R Z_0 Q &= -3 \int_\R Q'^2 + \int_\R \mathcal{L} P \hat P = -3 \int_\R Q'^2 + \int_\R P(3Q-2f(Q)).
\end{align*}
We recall now the following auxiliary result. 
\begin{Cl}[\cite{MMcol2}, Claim 2.2]\label{LQQ}
We have
$$
3 \int_\R Q'^2 - \int_\R (3Q-2f(Q))P= \int_\R \Lambda Q Q \neq 0. 
$$
\end{Cl}

\begin{rem}
Indeed,
$$
\int_\R \Lambda Q Q = \frac 12\partial_c \int_\R Q_c^2 \Big|_{c=1} >0,
$$
thanks to (\ref{stab}) provided $\ve$ small enough (independent of $c$). 
\end{rem}

Let us continue with the proof of Claim \ref{cl:ndg}. By the preceding result, it suffices to choose $\displaystyle{a:=\frac{\int_\R DQ}{\int_\R Z_0 Q}},$
and $\displaystyle{b:=\int_0^{+\infty} (D-aZ_0)(z) dz }$ (note that
$\lim_{\pm \infty} \mathcal L \varphi=\lim_{\pm \infty} \varphi=\pm 1$).  This finishes the proof of Claim \ref{cl:ndg}.
\end{proof}

\smallskip      

We return to the proof of Proposition \ref{prop:ex}. Now we find the constants $a$ and $b$ in terms of known quantities in $(\Omega)$. First, we multiply the equation of $B$ by $Q$ and use $\mathcal L Q'=0$. We get
\begin{align*}
-3a \int_\R Q'^2 & =  
\int_\R ( 3Q'' + f'(Q)Q) A + \int_\R G Q\\
&=  \int_\R (\mathcal L A ) P + \int_\R G Q.
\end{align*}
Second, we multiply the equation of $A$ by $\int_0^x P(s) \ ds$. We obtain
\begin{align*}
\int_\R (\mathcal L A)' \int_0^x P &= -\int_\R (\mathcal L A) P
+\ga \int_\R P\\
& =  a \int_\R (3 Q - 2f(Q)) P + \int_\R F \int_0^x P.
\end{align*}
Thus, combining the two identities, we get:
\begin{align*}
& -a \left\{
3 \int_\R Q'^2 - \int_\R (3 Q - 2f(Q)) P \right\} =  \gamma \int_\R P + \int_\R G Q - \int_\R F \int_0^x P. 
\end{align*}
and the expression for $a$ follows from Claim \ref{LQQ}.

\smallskip

To find out $b$, we integrate the equation for $B$ in $(\Omega)$ over $\R$ to obtain
\be\label{b10}
2 b = \int_\R f'(Q)A + \int_\R G.
\ee
Now we consider $\bar P$ the function defined in (\ref{barp}). We multiply the equation for $A$ by $\int_0^x \bar P(s) \ ds$ and then we integrate. We get
$$
\int_\R f'(Q)A = \ga\int_\R \bar P - a\int_\R \bar P(3Q-2f(Q)) -\int_\R F\int_0^x \bar P. 
$$
Now, note that
$$
\int_\R \bar P(3Q-2f(Q)) = \int_\R \mathcal L \hat P \bar P = \int_\R \hat P f'(Q) =\int_\R \hat P(1-\mathcal L1) = \int_\R \hat P -\int_\R (3Q-2f(Q)).
$$
From (\ref{hatp}) we replace the explicit value of $\hat P$ and we use the equation satisfied by $Q$, namely $Q'' -Q + f(Q)=0$, to obtain 
$$
\int_\R \bar P(3Q-2f(Q)) = -\int_\R \Lambda Q.
$$
With $a$ previously known we replace this quantity in (\ref{b10}) to obtain (\ref{b}). 
This finishes the proof.

\end{proof}

We have now a good solvability theory for the linear systems $(\Omega_{k,l})$, that avoids the emergency of linearly growing solutions at this order. As an example, the general theory constructed in \cite{MMcol1} for the quartic KdV equation deals with possibly growing solutions, see \cite{MMcol1} Proposition 2.3.

\medskip

Here, for each system $(\Omega_{k,l})$, $(k,l)\in \Sigma_m$, we will look for solutions such that
\begin{equation}\label{eq:st}
A_{k,l}=\widetilde A_{k,l} + \gamma_{k,l},\quad
B_{k,l}=\widetilde B_{k,l} + b_{k,l}\varphi + \kappa_{k,l}Q',\quad a_{k,l}, b_{k,l}, \kappa_{k,l}\in \R;
\end{equation}
where $\widetilde A_{k,l}\in \mathcal{Y}$ is even and 
$\widetilde B_{k,l}\in \mathcal{Y}$ is odd.
(see Proposition \ref{lem:z}  for a justification of this choice).

\medskip

This election will have several good properties, but we will emphasize a crucial one. Let $(k,l)\in \Sigma_m$ fixed. We say that $(k,l)$ satisfies the {\bf (IP)} property ({\bf IP} = important property) if and only if
\begin{equation*}
{\bf(IP)} \begin{cases}\hbox{ Any derivative of $A_{k,l}$ or $B_{k,l}$ is a localized $\mathcal{Y}$-function.} \\ \hbox{ Moreover, for $(k,l)=(1,0)$ we have $A_{1,0}\in \mathcal{Y}$.}
\end{cases}
\end{equation*}

This property, although depending on the specific pair $(k,l)$, will be useful to quickly discard localized terms composing $F_{k,l}, G_{k,l}$, and seeing essentially the bounded but non localized terms. Indeed, note that thanks to Claim \ref{surphi} any solution as in (\ref{eq:st}) satisfies this property. For the details, see Appendix \ref{sec:B}.  

We start by solving the first system.

\subsubsection{Resolution of the system $(\Omega_{1,0})$}

From Proposition \ref{prop:decomp} (ii) the system $(\Omega_{1,0})$ is given by
  \begin{align}
\label{A10}
& (\mathcal{L}A_{1,0})'=  -a_{1,0}(3Q - 2f(Q))' + (f'(Q))', \\
\label{B10}
& (\mathcal{L} B_{1,0})'= 3A_{1,0}'' + f'(Q)A_{1,0} - 3a_{1,0}Q''+ f'(Q)
  \end{align}

This first system is easily solvable, as shows the following
%%%%%%%%%%%%%%%%%%%%%
%%%%%%%%%%%%%%%%%%%%%

\begin{lem}[Resolution of $(\Omega_{1,0})$]\label{lem:omega10}

There exists a solution $(A_{1,0}, B_{1,0}, a_{1,0})$ of (\ref{A10})-(\ref{B10}) of the form (\ref{eq:st}) and such that $A_{1,0}\in \mathcal{Y}$ is even (and $\ga_{1,0}=0$), $B_{1,0}$ is odd and $a_{1,0}, b_{1,0}$ are given by the formulae 

\be\label{eq:solOmega10}
a_{1,0}   = \frac{ \int_\R \Lambda Q  }{\int_\R \Lambda Q Q} 
,\quad b_{1,0}  = \frac 12 a_{1,0}\int_\R \Lambda Q + \frac 12 \int_\R \bar P.
\ee

Moreover, $A_{1,0}$ is given by
\be\label{A}
A_{1,0} = \bar P -a_{1,0} \hat P.
\ee
(cf. (\ref{A10}), (\ref{hatp}) and (\ref{barp})). Finally, we choose $B_{1,0}$ such that $\displaystyle{\int_\R Q'B_{1,0} =0}$. 
\end{lem}

\begin{rem}
Note that from the value of $\hat P=-(xQ' + \Lambda Q)$ and (\ref{A}) we get 
\be\label{balt}b_{1,0}  = \frac 12 \Big[ a_{1,0}\int_\R Q +  \int_\R A_{1,0}\Big]. 
\ee
\end{rem}

\begin{proof}

Note that both $(f'(Q))'$ and $f'(Q)$ are odd and even $\mathcal{Y}$-functions respectively, so thanks to 
Proposition \ref{prop:ex}, a solution with the desired properties does exist. We will chose $\gamma_{1,0}: =0$.
The value of $a_{1,0}$ and $b_{1,0}$ comes from (\ref{a})-(\ref{b}), after some simple computations. These computations have been carried out in \cite{MMcol2}, but by completeness we rewrite them. Indeed, note that we only need to verify that 
$$
\int_\R f'(Q)(Q +  P) =-\int_\R \Lambda Q.
$$
In fact, from (\ref{P}), the explicit value of $P$ and Claim \ref{surL} (ii), we have
$$
\int_\R f'(Q)(Q +  P) =- \int_\R f'(Q) (xQ'+\Lambda Q) = \int_\R f(Q) -\int_\R ( 1-\mathcal L1) \Lambda Q = \int_\R (f(Q)-Q) -\int_\R \Lambda Q,
$$
but $f(Q)-Q = -Q''$, so we are done.

On the other hand, note that $\mathcal{L} (1+ \bar P) = 1$, thus
$$
\int_\R f'(Q)(\bar P +1) = \int_\R \mathcal L \bar P (1+\bar P) = \int_\R \bar P 1.
$$
This give finally the expected value of $b_{1,0}$. 
\smallskip

Finally, the constant $\kappa_{1,0}$ in the expression of $B_{1,0}$ is a free parameter
that we will fix such that $\int_\R B_{1,0} Q'=0$ for convenience in some future
computations (see Proposition \ref{prop:ex} and (\ref{unique}) in Appendix \ref{D}). We have
\begin{align*}
0 = \int_\R Q'B_{1,0}&= \int_\R \widetilde B_{1,0} Q'  + b_{1,0}\int_\R \varphi Q'+\kappa_{1,0}\int_\R Q'^2.
\end{align*}
where we can obtain $\kappa_{1,0}$.
\end{proof}

\subsubsection{Resolution of the system $(\Omega_{2,0})$}

From Proposition \ref{prop:decompbis} (iii) in Appendix \ref{sec:B}, the system $(\Omega_{2,0})$ is given by
  \begin{align}
\label{eq:A20}
& (\mathcal{L}A_{2,0})'=  a_{2,0}(3Q - 2f(Q))' +F_{2,0}, \\
\label{eq:B20}
& (\mathcal{L} B_{2,0})'= 3A_{2,0}'' + f'(Q)A_{2,0} - 3a_{2,0}Q''+ G_{2,0}
  \end{align}
  
where the source terms are given by
\begin{enumerate}
\item Case $m=2$,
\begin{eqnarray}
F_{2,0} & = &  -(3A_{1,0}'+3B_{1,0}'' +f'(Q)B_{1,0}) +\frac 12 (f''(Q)(2A_{1,0} +A_{1,0}^2))' \nonumber \\
& & \quad  - a_{1,0} ( 3A_{1,0}''  -Q +f'(Q)(1+A_{1,0})  )' + 3a_{1,0}^2 Q^{(3)} + \frac 12 (f''(Q)-2)' \label{F202},
\end{eqnarray}
and
\begin{eqnarray}
G_{2,0} & = &  \frac 12 (f''(Q)-2) -(A_{1,0}+3B_{1,0}')  + \frac 12 f''(Q)(2A_{1,0}+A_{1,0}^2) + \frac 32 a_{1,0}^2 Q ''\nonumber \\
& &    -\frac 12 a_{1,0}(9A_{1,0}' +3B_{1,0}'' +f'(Q)B_{1,0})'  +\frac 12 (f''(Q)(B_{1,0} +A_{1,0}B_{1,0}))' .\label{G202}
\end{eqnarray}
\item Case $m=3$,
\be\label{F203}
F_{2,0}  =  (\frac 12 f''(Q)(1+A_{1,0})^2)' +3a_{1,0}^2 Q^{(3)} -a_{1,0} (f'(Q) +3A_{1,0}'' +f'(Q)A_{1,0})'  
\ee
and
\bea\label{G203}
G_{2,0} & = & \frac 12 f''(Q)(1+A_{1,0})^2 +\frac 32 a_{1,0}^2 Q'' -\frac 12 a_{1,0}(9A_{1,0}' +3B_{1,0}'' +f'(Q)B_{1,0})'  \nonumber\\
& & + \frac 12 (f''(Q)(1+A_{1,0})B_{1,0})'.
\eea
\end{enumerate}

%%%%%%%%%%%%%%%%%%%%%%%%%%%%%%%%
%%%%%%%%%%%%%%%%%%%%%%%%%%%%%%%%

\begin{prop}[Resolution of $(\Omega_{2,0})$]\label{lem:b20}
Let $f$ be as in (\ref{gKdV0b}). There exists a constant $\ve_0>0$ not depending on $c$ such that the following holds.

\begin{enumerate}
\item (Case $m=2$)
There exists a solution $(A_{2,0},B_{2,0},a_{2,0})$ of 
$(\Omega_{2,0})$ satisfying (\ref{eq:st}) and such that
$$ \lim_{+\infty} A_{2,0}= - \frac 1{2} b_{1,0}^2=\gamma_{2,0},\quad
A_{2,0} -\gamma_{2,0} \in \mathcal{Y},$$
$$ \lim_{+\infty} B_{2,0}= b_{2,0},\quad
B_{2,0}- b_{2,0} \varphi \in \mathcal{Y},$$
but for all $|\ve|\in (0,\ve_0)$
\begin{equation}\label{diff0}
d(\ve) := b_{2,0}(f) + \frac 1{6} b_{1,0}^3(f) = c_{2,p}\ve + o(\ve), \quad \hbox{ with } \quad c_{2,p}\neq 0\;  \hbox{ for all } \; p\geq 4.
\end{equation}
\item (Case $m=3$)
There exists a solution $(A_{2,0},B_{2,0},a_{2,0})$ of 
$(\Omega_{2,0})$ such that $A_{2,0}\in \mathcal{Y}$ is even, $B_{2,0}$ is bounded, odd and
$$ \lim_{+\infty} B_{2,0}= b_{2,0},\quad
B_{2,0}- b_{2,0} \varphi \in \mathcal{Y},$$
but for all $\abs{\ve}\in (0,\ve_0)$, 
\begin{equation}\label{diff2}
 d(\ve):= b_{2,0}(f)= c_{3,p}\ve + o(\ve), \quad \hbox{ and } \quad  c_{3,p}\neq 0 \; \hbox{ for all }\; p\geq 4.
\end{equation} 
\end{enumerate}
Moreover, in both cases the solution found satisfies {\bf (IP)}.
\end{prop}

\begin{rem}
Note that in the case $m=2$, one has $c_{2,p} = 0$ for $p=3$ (see (\ref{c2p})). This cancelation is consequence of the complete integrability of the Gardner equation.   
\end{rem}
%%%%%%%%%%%%%%%%%%%%%%%%%%%%%%%%
%%%%%%%%%%%%%%%%%%%%%%%%%%%%%%%%

\begin{proof}
Note that in both cases, $m=2$ and $m=3$ the source terms  $F_{2,0}, G_{2,0}$ belongs to $\mathcal{Y}$, with the former being an odd function and the last one being even. Thus the existence of solutions to (\ref{eq:A20})-(\ref{eq:B20}) with the desired properties follows directly from Proposition
\ref{prop:ex} above. 

In particular we will choose $\ga_{2,0}: = -\frac 12 b_{1,0}^2$ for the quadratic case and $\ga_{2,0}:=0$ in the cubic one.

Let us now check that, being fixed $\ga_{1,0}, a_{1,0}, b_{1,0}$ and  $\ga_{2,0} $, the value of $a_{2,0}$ and $b_{2,0}$ is uniquely determined. Indeed, from (\ref{a})-(\ref{b})

\be\label{a20}
a_{2,0} = -\frac{1}{\int_\R \Lambda Q Q }\Big[ \ga_{2,0} \int_\R P +\int_\R G_{2,0}Q  - \int_\R  F_{2,0}\int_0^x P\Big] ,
\ee
and
\be\label{b20}
b_{2,0} = \frac 12\Big[ \ga_{2,0}\int_\R \bar P + a_{2,0}\int_\R \Lambda Q -\int_\R F_{2,0}\int_0^x \bar P+ \int_\R G_{2,0}\Big] .
\ee

We claim (\ref{diff0}) and (\ref{diff2}) with
\be\label{c2p}
c_{2,p} := -\Big[ \frac{(p-3)(2p-1) (24 -23p +3 p^2  +2p^3)}{36(p^2 -1)(p-2)} \Big] \int_\R \Big[ \frac{3}{2\cosh^2(\frac 12 x)}\Big]^p,
\ee
and
\be\label{c3p}
c_{3,p} := -\Big[ \frac{(p-1)(p-3)(p^2-3p+8) }{8(p-2)(p+1)}\Big]\int_\R \Big[\frac{\sqrt{2}}{\cosh x }\Big]^p.
\ee
The end of the proof of (\ref{diff0})-(\ref{diff2}), and (\ref{c2p})-(\ref{c3p}) is a lengthy but straightforward computation. For the sake of continuity we postpose the proof to Appendix \ref{D}. %The proof of Proposition \ref{lem:b20} is now complete.
\end{proof}

\begin{rem}\label{smallep}
An explicit expression for $d(\ve)$ for any nonlinearity has escaped to us (see Claim \ref{CL}), and we only have in our hands an asymptotic expression for small values of $\ve$. We believe, however, that it may exist a --necessarily-- large $\ve_0$ for which $d(\ve_0)=0$, and even more, a pure 2-soliton solution may exist at any order. 
\end{rem}

\begin{rem}\label{Sch0}
The expressions (\ref{diff0})-(\ref{diff2}) above say roughly speaking that the second order linear system $(\Omega_{2,0})$ has a solution that does not obey (at third order derivatives) the Taylor expansion of a small soliton shifted. Indeed,
$$
 Q_c(y_c + b_{1,0}\varphi )\sim Q_c(y_c) + b_{1,0}\varphi Q_c'(y_c) + \frac 12 b_{1,0}^2 Q_c''(y_c) + \frac 16 b_{1,0}^3 Q_c^{(3)}(y_c) .
$$
Note that (cf. (\ref{taylor0}) and (\ref{taylor1}))
$$
Q_c''(y_c) \sim cQ_c(y_c) - Q_c^m(y_c), \quad Q_c^{(3)}(y_c) \sim cQ_c'(y_c) - (Q_c^m)'(y_c),
$$
and thus for a perfect collision we should have $b_{2,0} = -\frac 16 b_{1,0}^3$ for $m=2$ and $b_{2,0}=0$
for $m=3$, as in the integrable cases. This formal discussion will be  justified in the proof of Proposition \ref{lem:z}.

\end{rem}

%%%%%%%%%%%%%%%%%%%%%%%%%%%%%%%%%
%%%%%%%%%%%%%%%%%%%%%%%%%%%%%%%%%

\medskip

\subsubsection{Resolution of system $(\Omega_{1,1})$, cases $m=2,3$}
Now we consider the first mixed system, $(\Omega_{1,1})$. Note that this system has a different order depending on the power of leading nonlinearity: for $m=2$, $ cQ_c$ is of quadratic order in $Q_c$, meanwhile, in the cubic one, $cQ_c$ is a term of cubic order.

From Proposition \ref{prop:decompbis} the system $(\Omega_{1,1})$ is given by
\begin{align}
\label{A11}
& (\mathcal{L}A_{1,1})'=  a_{1,1}(3Q - 2f(Q))' + (3A_{1,0}' + 3B_{1,0}'' + f'(Q)B_{1,0}) , \\
\label{B11}
& (\mathcal{L} B_{1,1})'= 3A_{1,1}'' + f'(Q)A_{1,1} - 3a_{1,1}Q'' + 3B_{1,0}'.
  \end{align}

For this system, we recall its source terms
\be\label{F11G11}
F_{1,1}:=3A_{1,0}' + 3B_{1,0}'' + f'(Q)B_{1,0}, \quad G_{1,1} := 3B_{1,0}'.
\ee
Note that as $(k,l)=(1,0)$ satisfies the {\bf (IP)} property, we have both $F_{1,1}, G_{1,1}\in \mathcal{Y}$.  

\begin{lem}[Resolution of $(\Omega_{1,1})$, $m=2,3$]\label{lem:b11}
There exists a solution $(A_{1,1}, B_{1,1}, a_{1,1})$ of 
$(\Omega_{1,1})$ such that $A_{1,1}$ is even, $B_{1,1}$ is odd and
$$ \lim_{+\infty} A_{1,1}  =\gamma_{1,1}:=\frac 12 b_{1,0}^2 ,\quad
A_{1,1} -\gamma_{1,1} \in \mathcal{Y},$$
$$ \lim_{+\infty} B_{1,1}= b_{1,1},\quad
B_{1,1}- b_{1,1} \varphi \in \mathcal{Y}.$$
Besides, this solution implies that {\bf (IP)} holds for $(k,l)=(1,1)$. 
\end{lem}

\begin{proof} From Proposition \ref{prop:decomp}, it is clear that
$F_{1,1}$ and $G_{1,1}$ given in (\ref{F11G11}) satisfy the assumptions of Proposition \ref{prop:ex}.
The choice of $\gamma_{1,1}$ will be justified in Proposition \ref{lem:z}.
In the rest of this paper, we will not need the expression of $b_{1,1}$
(note that it would be possible to compute it as in the proof of Proposition \ref{lem:b20}).
\end{proof}

\subsubsection{Resolution of high order systems, quadratic case}

From now on, we consider the triplet $$(A_{k,l},B_{k,l}, a_{k,l})$$ defined for all $(k,l)\in \Sigma_m$, $1\leq k+l\leq 2$ in Lemma \ref{lem:omega10}, Proposition \ref{lem:b20} and Lemma \ref{lem:b11}.
We now solve the systems $(\Omega_{k,l})$ for $k+l=3$ . Denote $\delta_{33}:= 1$ and $\delta_{p3} :=0$ for $p\geq 4$.

\begin{lem}[Resolution of $(\Omega_{k,l})$ for $k+l=3$ and $m=2$]\label{lem:2.5}

For all $(k,l)\in \Sigma_2$  such that $k+l=3$, $F_{k,l}$ is odd and $G_{k,l}$ even; both are in the class $\mathcal{Y}$, and there exists a solution $(A_{k,l}, B_{k,l}, a_{k,l}) $ of 
$(\Omega_{k,l})$ such that $A_{k,l}$ is even, $B_{k,l}$ is odd and
$$ \lim_{+\infty} A_{k,l}  =\gamma_{k,l},\quad
A_{k,l} -\gamma_{k,l} \in \mathcal{Y},$$
$$ \lim_{+\infty} B_{k,l}= b_{k,l},\quad
B_{k,l}- b_{k,l} \varphi \in \mathcal{Y}.$$
Moreover, we will choose the particular values
\begin{align*}
&\gamma_{3,0}:=\frac 5 {36} b_{1,0}^4+ \frac {10}{3} d(\ve)b_{1,0} + \frac 12 \mu(\ve) b_{1,0}^2,\quad \gamma_{2,1}:= \frac 1{24} b_{1,0}^4 -b_{1,0}b_{1,1} - 4 d(\ve) b_{1,0}
,\\
&  \gamma_{1,2}:=-\frac 3{24} b_{1,0}^4 + b_{1,0}b_{1,1},
\end{align*}

where $d(\ve)$ satisfies (\ref{diff0})-(\ref{diff2}).
\end{lem}

\begin{proof}
The proof of this result is easy after the validity of the following claim:
\begin{equation}\label{eq:cl}
\hbox{For all $(k,l) \in \Sigma_2$   such that $k+l=3$, we have
$F_{k,l}\in \mathcal{Y}$ is odd, $G_{k,l} \in \mathcal{Y}$ is even. }
\end{equation}
Assuming \eqref{eq:cl}, 
Lemma \ref{lem:2.5} is a direct consequence of Proposition \ref{prop:ex}.

Let us prove \eqref{eq:cl}. 
From the Appendix \ref{sec:B} and Proposition \ref{prop:decompbis} several (bounded but) nonlocalized terms appear in the expression of 
$F_{k,l}$ and  $G_{k,l}$ for $k+l=3$, but {\bf all these terms eventually cancel}. 

Indeed, thanks to the {\bf (IP)} property, terms containing {\bf derivatives} of $B_{1,0}$, $A_{1,1}$ and $A_{2,0}$ are in $\mathcal{Y}$ as well as terms of the kind $f'(Q)B_{1,0}$ and so on.
Thus, we focus on the terms containing only $B_{1,0}$, $A_{1,1}$ and
$A_{2,0}$ without derivatives nor multiplication by functions of $Q$. Note also that $A_{1,0}\in \mathcal{Y}$, so we also discard it.  For simplicity of notation, we will skip the variables $\ys$ and $y$.

\smallskip

Now, we recollect all the non-localized terms (due to $B_{1,0}$, $A_{1,1}$ and $A_{2,0}$) in $S[\tilde u]$ of order $c^l \wqs^k$ or $c^l (\wqs^k)' $ with $k+l=3$. We have only three cases: the pairs $(3,0), (2,1)$ and $(1,2)$.  From Proposition \ref{prop:decompbis} we obtain

\begin{enumerate}
\item (Case $(3,0)$). Here 
$$
F_{3,0} = \widetilde{F}_{3,0}, \quad G_{3,0} =  \widetilde{G}_{3,0} -\frac 23 (B_{1,0}^2 + 2A_{2,0}), \, \hbox{ with } \widetilde{F}_{3,0}, \widetilde{G}_{3,0} \in\mathcal{Y};
$$
\item (Case $(2,1)$). Here 
$$
F_{2,1} = \widetilde{F}_{2,1}, \quad G_{2,1} =  \widetilde{G}_{2,1} + ( B_{1,0}^2 + A_{1,1} + 3A_{2,0}), \, \hbox{ with } \widetilde{F}_{2,1}, \widetilde{G}_{2,1} \in\mathcal{Y};
$$
\item (Case $(1,2)$). Here $F_{1,2} , G_{1,2}\in\mathcal{Y}$.
\end{enumerate}

Using the following relations among the limits of $A_{2,0}$, $A_{1,1}$
and $B_{1,0}^2$ at $\pm \infty$ (see Proposition \ref{lem:b20} and Lemma \ref{lem:b11}):
$$
\lim_{\pm\infty} A_{2,0}=-\frac 1{2} \lim_{\pm\infty} B_{1,0}^2,\quad
\lim_{\pm\infty} A_{1,1}=- \lim_{\pm\infty} A_{2,0},
$$
we observe that the source functions in $(\Omega_{k,l})$ are in fact all localized. This proves (\ref{eq:cl}).
\end{proof}

\medskip

\subsubsection{Resolution of high order systems, cubic case}

Finally we claim the existence of bounded solutions for the \emph{third} and \emph{fourth} order systems in the cubic case. The proof of these results is identical to the previous Lemma.

\begin{lem}[Resolution of $(\Omega_{3,0})$, $(\Omega_{4,0})$ and $(\Omega_{2,1})$ for $m=3$]\label{lem:3.5}
For all $(k,l)\in \Sigma_3$ with $k\geq 2$
there exists a solution $(A_{k,l}, B_{k,l}, a_{k,l})$ of 
$(\Omega_{k,l})$ such that $A_{k,l}$ is even, $B_{k,l}$ is odd and
$$ \lim_{+\infty} A_{k,l}  =\gamma_{k,l},\quad
A_{k,l} -\gamma_{k,l} \in \mathcal{Y},$$
$$ \lim_{+\infty} B_{k,l}= b_{k,l},\quad
B_{k,l}- b_{k,l} \varphi \in \mathcal{Y}.$$
In particular, we choose
\be\label{g3}
\ga_{3,0}:=-\frac 12 b_{1,0}^2, \; \ga_{2,1}:= -4b_{1,0}d(\ve),\; \ga_{4,0} := 3d(\ve) b_{1,0}+\frac 12 \ve b_{1,0}^2\delta_{p4}.
\ee
In this case $d(\ve):= b_{2,0}(\ve)$ (cf. (\ref{diff2})).
\end{lem}

\begin{proof}
We note that, thanks to the {\bf (IP)} property and Proposition \ref{prop:decompbis}, the only a priori non localized source term is
$$
G_{2,1} = 3A_{2,0} + \tilde G_{2,0}, \hbox{ with } \tilde G_{2,0}\in \mathcal{Y}.
$$
Then the conclusion of the Lemma follows from the fact that, from Proposition \ref{lem:b20}, in the cubic case, we have a priori chosen $A_{2,0}\in \mathcal{Y}$. 
\end{proof}

\medskip

For further purposes, we recall the important quantities (see (\ref{Tc}) and (\ref{diff0})-(\ref{diff2}))
\be\label{eq:dtau}
        T_c =c^{-\frac 12 -\frac 1{100}}, \quad
        d(\ve) = b_{2,0}(\ve)+\frac 1{6} b_{1,0}^3(\ve)\delta_{m2}, 
\end{equation}
with $\delta_{m2}=0$ for $m=3$, and $\delta_{22} =1$.

%%%%%%%%%%%%%%%%

\subsection{Recomposition of the approximate solution. Proof of Proposition \ref{Constr}}\label{sec:3.7}

Having solved several linear systems we now are able to prove Proposition \ref{Constr}. Indeed, we have now 
the enough knowlegde about the notation, so we can go further and claim the following improved result on $\tilde u$.

\begin{prop}[Construction of a symmetric approximate solution of gKdV, improved version]\label{lem:z}\
The solution $\tilde u $ above constructed satisfies, for any $0<c<c_{0}$, the following properties:
\begin{enumerate}
\item For all $(t,x)$ $\tilde u(t,x)= \tilde u(-t,-x)$.
\item For every time $t\in [-T_c,T_c]$,
\begin{equation}\label{eq:z2}
\left\|  S[\tilde u ](t) \right\|_{H^2(\R)}\leq K c^{\frac 3{m-1} + \frac 34}.
\end{equation}
\item Closeness to the sum of two soliton solution: For all time $t\in [-T_c, T_c]$, the function $\tilde u $ is in $H^1(\R)$ and satisfies the estimate
\be\label{8}
\norm{\tilde u(t) -Q(y) -Q_c(y_c) }_{H^1(\R)}\leq K_0 c^{\frac 1{m-1}}.
\ee
\item Closeness to a shifted two soliton solution plus a strange term: Denote
\begin{equation}\label{eq:z3}
\Delta_1 : =\sum_{(k,l)\in \Sigma_m} a_{k,l}\, c^l \int_\R \wqs^k,\quad
\tilde b_{1,1}:=b_{1,1} -\frac 16 b_{1,0}^3,
\quad \Delta_2 := 2 (b_{1,0}+ c \tilde b_{1,1}\delta_{m2} ).
\end{equation}
Then $\tilde u$ satisfies at time $\pm T_c$
\begin{equation}\label{eq:z1}
\begin{split}
& \| \tilde u (\pm T_c)- Q(\cdot \mp \frac 12\Delta_1) - \wqs(\cdot  \pm (1-c)T_c \mp \frac 12 \Delta_2) \\
&\quad  \quad \quad   \pm d(\ve) (\wqs^2)'(\cdot  \pm (1-c) T_c  \mp \frac 12 \Delta_2)\}\|_{H^1(\R)}  \leq K c^{\frac{3}{m-1} + \frac 14},
\end{split}
\end{equation}
provided for each $(k,l)\in \Sigma_m$, the constants $\ga_{k,l}$ must be chosen as in Lemma \ref{lem:omega10}, Proposition \ref{lem:b20}, Lemmas \ref{lem:b11}, \ref{lem:2.5} and \ref{lem:3.5}. Recall that $d(\ve)$ satisfies (\ref{diff0})-(\ref{diff2}).
\end{enumerate}
\end{prop}

\begin{rem}
The quantity $\bar b_{1,1}$ in (\ref{eq:z3}) represents the difference between the expected value of $b_{1,1}$ given by the integrable case and the actual one; namely, for $\ve=0$ we have $\bar b_{1,1} =0$.
\end{rem}

\begin{proof}

Let us start by proving (\ref{eq:z2}).  This follows from from Proposition \ref{prop:decompbis}, and the choice of $a_{k,l}$, $A_{k,l}$, $B_{k,l}$ for $(k,l)\in \Sigma_m$, solving each linear system 
$(\Omega_{k,l})$, so that  
$$S[\tilde u] =\mathcal{E}(t,x).$$

\smallskip

Now we deal with (\ref{8}). This is an easy consequence of the fact that $y=x-\al$, $y_c =x+(1-c)t$, and
$$
\tilde u(t) -Q(y) -Q_c(y_c) = W(t,x), \quad \norm{W(t)}_{H^1(\R)}\leq K_0 c^{\frac 1{m-1}}.
$$

\smallskip

Proof of \eqref{eq:z1}. We begin with some preliminary estimates.

\begin{Cl}\label{cl:30}
\begin{equation}\label{eq:al}
\|\alpha\|_{L^\infty}\leq K c^{\frac1{m-1} -\frac 12},\quad \|\alpha'\|_{L^\infty}\leq K c^{\frac 1{m-1}}.
\end{equation}
Suppose $f=f(y)\in \mathcal Y$. Then for all $t\in [-T_c,T_c]$,
\be\label{fQ}
\norm{f(y)Q_c^k(y_v)}_{L^2(\R)} + \frac 1{\sqrt{c}} \norm{f(y)(Q_c^k)'(y_v)}_{L^2(\R)} \leq K c^{\frac k{m-1}} e^{-(1-c)\sqrt{c}\abs{t}},
\ee
and for $g=g(y)\in L^{\infty} (\R)$,
\be\label{gQ}
\norm{g(y)Q_c^k(y_v)}_{L^2(\R)} + \frac{1}{\sqrt{c}} \norm{g(y)(Q_c^k)'(y_v)}_{L^2(\R)} \leq K c^{\frac{k}{m-1} -\frac 14}.
\ee
In particular, if $t=T_c$ and $f\in \mathcal{Y}$, we have, for $ c>0$ small,
\begin{equation}\label{eq:de}
\|f(y)\wqs(\ys)\|_{H^1(\R)}\leq K c^{10},
\end{equation}
\begin{equation}\label{eq:Qy}
\|Q(y)-Q(x-\tfrac 12 \Delta_1)\|_{H^1(\R)}\leq K c^{10}.
\end{equation}
\end{Cl}

\begin{proof}
The proof of these estimates are similar to Claim C.1 in the Appendix C of \cite{MMMcol}. See also Claim 2.6 in \cite{MMcol1}. In particular for the proof we use Lemma \ref{surQc2} from Appendix \ref{sec:B}. We skip the details.
\end{proof}

We continue the proof of (\ref{eq:z1}). 

For the sake of brevity, we will prove only the case $m=3$. The case $m=2$ is identical to Lemma 2.6 in \cite{MMMcol}.

Note that from Claim \ref{cl:30},
\begin{equation}\label{eq:c2}
\| \wqs(\ys-b_{1,0})  - \wqs  + b_{1,0} \wqs' - \frac 12 b_{1,0}^2 \wqs'' \|_{H^1(\R)}  \leq K c^{\frac {7}4},
\end{equation}
and
\begin{equation}\label{eq:c2bis}
\|(\wqs^2)'(\ys-b_{1,0}) -(\wqs^2)' + b_{1,0}(\wqs^2)'' \|_{H^1(\R)}\leq K c^{\frac {9}4},
\end{equation}
(here we have used the fact $\| Q_c^{(3)}\|_{H^1(\R)} \leq K c^{\frac {7}4}$ and $\|(Q_c^2)^{(3)}\|_{H^1(\R)} \leq K c^{\frac {9}4}$). From the identities
$$
Q_c'' =cQ_c-Q_c^3 - \ve Q_c^p + O(Q_c^{p+1}), \; (Q_c^2)'' = 4cQ_c^2 -3Q_c^4 + O(Q_c^5),
$$
we obtain 
\begin{equation}\label{eq:c3}
\begin{split}
&  \big\| \wqs(\ys-b_{1,0}) - d(\ve)(\wqs^2)'(\ys - b_{1,0})\\
& \quad \quad  -[ Q_c -b_{1,0}Q_c' + \frac 12 b_{1,0}^2 c Q_c -\frac 12 b_{1,0}^2 Q_c^3 - \frac 12 \ve b_{1,0}^2\delta_{p4} Q_c^4   ] \\
& \quad \quad  + d(\ve)[ (Q_c^2)'   -4b_{1,0}cQ_c^2 +3b_{1,0}Q_c^4 ] \big\|_{H^1(\R)}\leq K c^{\frac {7}4}.
\end{split}
\end{equation}
On the other hand, using  the fact that $\lim_{+\infty} A_{k,l}=\ga_{k,l}$, $\lim_{+\infty} B_{k,l}=b_{k,l}$, and Claim \ref{cl:30} we get
\bee
& & \big\|\tilde u(T_c) -Q -Q_c -b_{1,0}Q_c' -\ga_{2,0}Q_c^2 -b_{2,0}(Q_c^2)' -\ga_{1,1}cQ_c \\
& & \quad \quad -\ga_{2,1}cQ_c^2 -\ga_{3,0}Q_c^3 -\ga_{4,0}Q_c^4 \big\|_{H^1(\R)} \leq Kc^{7/4}.
\eee
Combining this estimate and (\ref{eq:c3}), we find
\begin{equation*} \begin{split}
&  \big\| \tilde u(T_c) - \{Q(y)+\wqs(\ys-b_{1,0}) - d(\ve)(\wqs^2)'(\ys - b_{1,0})\} \\
& +(\gamma_{1,1}-\frac 12 b_{1,0}^2) c \wqs +\gamma_{2,0}\wqs^2   + (b_{2,0} -d(\ve))(Q_c^2)'+ (\gamma_{2,1}+4d(\ve)b_{1,0}) c \wqs^2 \\
& + (\gamma_{3,0} +\frac 12 b_{1,0}^2 ) \wqs^3 + (\ga_{4,0}  -3d(\ve)b_{1,0} - \frac 12 \ve b_{1,0}^2 \delta_{p4}) Q_c^4  \big\|_{H^1(\R)}\leq K c^{\frac {7}4}.
\end{split}
\end{equation*}
It follows that with the choice
$$
\ga_{1,1}=\frac 12 b_{1,0}^2,\quad \ga_{2,0}=0, \quad \ga_{3,0}=-\frac 12 b_{1,0}^2, \quad b_{2,0} =d(\ve), 
$$
$$
\ga_{2,1} = -4d(\ve)b_{1,0} \quad \hbox{ and }  \quad \ga_{4,0}= 3d(\ve)b_{1,0} +\frac 12 \ve b_{1,0}^2\delta_{p4}.
$$
we obtain
\begin{equation}\label{eq:c9}
\| \tilde u (T_c)- Q(y)- \wqs(\ys - b_{1,0}) + d(\ve) (\wqs^2)'(\ys - b_{1,0}) \|_{H^1(\R)}
\leq K c^{\frac {7} 4}.
\end{equation}
The case $t=-T_c$ is similar and we left the proof to the reader.

Together with \eqref{eq:Qy}, we complete the proof of \eqref{eq:z1}. 
This justifies in particular the choices of $\gamma_{k,l}$, $(k,l)\in \Sigma_m$ done in preceding Lemmas.
\end{proof}

\medskip

\subsection{Existence of the approximate pure $2$-soliton collision solution}\label{se:26}

The fact that $d(\ve)\neq 0$ (see Proposition \ref{lem:b20}) in Proposition \ref{lem:z} means formally that the collision is not elastic and that the residue due to the collision is of order 
$(\wqs^2)'$.
However, the approximate solution $\tilde u(t,x)$ given in Lemma \ref{lem:z} being symmetric, it contains the residue at both $-T_c$ and $T_c$ (see \eqref{eq:z1}). To match the solution $u(t)$ considered in Theorem \ref{PUREgen}, which is pure at $-\infty$, we need to introduce a modified approximate solution, which, at  main order, will contain a residue only at time $t=T_c$. This will be clear after the following

\begin{prop}\label{cor:1}
There exists a function  $\hat u=\hat u (t,x)$, of the form given by (\ref{defW}) such that for some constants $K,c_{0}>0$ and $0<c<c_{0}$, the following estimates hold:
\begin{enumerate}
\item $\hat u(t,x)\not\equiv \hat u(-t,-x)$ for every $t,x$. 
\item Almost solution. For any $t\in [-T_c,T_c]$, 
\begin{equation}\label{eq:approx1}
\left\| S[\hat u](t) \right\|_{H^1(\R)}\leq K c^{\frac 2{m-1} +\frac 34} [c^{\frac 1{m-1}} + \abs{d(\ve)}c^{\frac 14}],
\end{equation}
(recall that $d(\ve)$ measures the residue after the collision, introduced in (\ref{diff0})-(\ref{diff2})).

\item Closeness to a two-soliton solution at time $t=-T_c$. With the defintions of shifts given in (\ref{eq:z3}), the modified function $\hat u$ is close to a two solitons solution at time $-T_c$:
\be\label{eq:zmT}
 \| \hat u(-T_c) - \{Q (\cdot + \tfrac 12 \Delta_2)
+Q_{c}(\cdot + (1-c)T_c +\tfrac 12\Delta_2)\}
\|_{H^1(\R)}\leq K c^{\frac 2{m-1} +\frac 14} [ c^{\frac 1{m-1}} + \abs{d(\ve)}c^{\frac 12}]. 
\ee
 
\item Non-matching with  a two-soliton solution at time $t=T_c$: 
\be\label{eq:zT}
\begin{split}
& \| \hat u(T_c) - Q (\cdot  -\frac 12 \Delta_1) -Q_{c}(\cdot +(1-c) T_c -\frac 12 \Delta_2) \\
& \quad + 2 d(\ve) (Q_{c}^2)'(\cdot  +(1-c)T_c-\frac 12 \Delta_2) \|_{H^1(\R)} \leq  K c^{\frac 2{m-1} +\frac 14} [ c^{\frac 1{m-1}} + \abs{d(\ve)}c^{\frac 12}]. 
\end{split}
\ee
where, from (\ref{diff0})-(\ref{diff2}), 
$$\forall 0< \abs{\ve} \leq  \ve_0, \quad d(\ve)= c_{m,p}\ve + o(\ve),$$
and
\begin{equation}\label{eq:Dt}
\abs{\Delta_1 -   a_{1,0}\int_\R Q_c } \leq K c^{\frac 2{m-1} - \frac 12}, \quad \abs{\Delta_2 - 2 b_{1,0}} \leq K c.
\end{equation}

\item Comparison residue versus error terms: The residue in (\ref{eq:zT}) satisfies
\bee\label{def}
\big\| 2 d(\ve) (Q_{c}^2)'(\cdot  +(1-c)T_c-\frac 12 \Delta_2)\big\|_{H^1(\R)} & \sim & \abs{d(\ve)}c^{\frac 2{m-1} +\frac 14}\\
& \gg & c^{\frac 2{m-1} +\frac 14} [ c^{\frac 1{m-1}} + \abs{d(\ve)}c^{\frac 12}],
\eee
provided $c^{\frac 1{m-1}} \ll \abs{d(\ve)}$.
\end{enumerate}
\end{prop}

\begin{rem}
The approximate solution $\hat u$ above mentioned describes the collision of two pure solitons that at time $t\sim T_c$ (after colliding) differ by a term of order $ \abs{d(\ve)}c^{\frac 2{m-1} +\frac 14}$ of the ingoing solitons before the collision, at time $t\sim -T_c$, provided (\ref{Cond0}) holds. 

For even small values of $\ve$ such that condition (\ref{Cond0}) does not hold, we need to go further in our approximate solution and solve even more linear systems. We believe that in this case, more involved, the conclusions of this paper are the same. 
\end{rem}

Let us return to the proof of Proposition \ref{cor:1}.

\begin{proof}
Let $\hat u:= \tilde u + w_{\#}$, where 
\be\label{gatow}
w_{\#}(t,x):= - d(\ve)(Q_c^2)'(y_c)(1+\bar P(y)), 
\ee
and $\bar P$ was defined in (\ref{barp}). Now $w_\#$ can be expressed in the form
$$
w_\# (t,x) = Q(y) + Q_c(y_c) + \sum_{(k,l)\in \Sigma_m }c^l \{ \hat A_{k,l}(y) Q_c^k(y_c) + \hat B_{k,l}(y) (Q_c^k)'(y_c) \},
$$
where $\hat A_{k,l} = A_{k,l}$, $\hat B_{k,l} = B_{k,l} + w_\# \delta_{(k,l),(2,0)} $. Here $ \delta_{(2,0),(2,0)} =1$ and $ \delta_{(k,l),(2,0)} =0$ otherwise.  

Let us prove (\ref{eq:zT}).  
Replacing $\tilde u =\hat u - w_{\#}$ in (\ref{eq:z1}), we have
\begin{equation*}\begin{split}
&\| \hat u(T_c)- Q(\cdot -\frac 12 \Delta_1) 
        - \wqs( \cdot +(1-c) T_c - \frac 12 \Delta_2)\\
        & \quad + d(\ve) (\wqs^2)'(\cdot +(1-c) T_c -\frac 12 \Delta_2) 
        - w_{\#}(T_c) \|_{H^1(\R)}  \leq K c^{\frac 3{m-1}  + \frac {1}4} .
\end{split}\end{equation*}
Thus, using \eqref{eq:de} (note that $\bar P\in \mathcal{Y}$)
\bee
 & & \| \hat u(T_c)- Q( \cdot -\frac 12 \Delta_1) - \wqs(\cdot +(1-c)T_c - \frac 12 \Delta_2) + 2 d(\ve) (\wqs^2)'(\cdot +(1-c) T_c - \frac 12 \Delta_2)  \|_{H^1(\R)}\\
 & & \leq K c^{\frac 3{m-1} + \frac {1}4} +\| d(\ve) (\wqs^2)'(\cdot +(1-c) T_c -\frac 12 \Delta_2) 
 + w_{\#}(T_c) \|_{H^1(\R)} \\
 & &    \leq K c^{\frac 3{m-1} + \frac {1}4} + K  \abs{ d(\ve)} \|  (\wqs^2)'(\cdot -\frac 12 \Delta_2)  -(\wqs^2)' \|_{H^1(\R)}\\
  & & \leq K c^{\frac 3{m-1} +\frac 14 } +K \abs{ d(\ve)} c^{\frac 2{m-1} +\frac 34}.
\eee
Similarly, at time $t=-T_c$
\bee               
& & \|\hat u(-T_c)- Q(\cdot +\frac 12 \Delta_1) - \wqs(\cdot -(1-c) T_c +\frac 12\Delta_2) \\
& & \quad - d(\ve) (\wqs^2)'(\cdot -(1-c)T_c +\frac 12 \Delta_2)-w_{\#}(-T_c) \|_{H^1(\R)}  \leq  K c^{\frac 3{m-1}+ \frac 14},
\eee
so that
\bee
& & \| \hat u(-T_c)- Q(\cdot +\frac 12 \Delta_1) - \wqs(\cdot -(1-c)T_c +\frac 12 \Delta_2)\|_{H^1(\R)}\\
& & \quad \quad \quad  \leq  K c^{\frac 3{m-1}+ \frac 14} + K\abs{d(\ve)} \|  (\wqs^2)'(\cdot +\frac 12 \Delta_2) -(\wqs^2)' \|_{H^1(\R)} \\
& &\quad\quad\quad  \leq Kc^{\frac 3{m-1}+ \frac 14} + K\abs{d(\ve)}c^{\frac 2{m-1} +\frac 34}.
\eee

Note that \eqref{eq:Dt} is clearly a consequence of \eqref{eq:z3}.

\medskip

Finally, we prove  (\ref{eq:approx1}). Note that (cf. Appendix B for the definitions)
\bee
S[\hat u] & = &S[\tilde u + w_{\#}]\\
& = & S[\tilde u] + {\bf III}(w_\#) + [ f(\tilde u + w_\#) -f(\tilde u) -f'(Q)w_\#]_x 
\eee

 The following estimates allow to conclude (\ref{eq:approx1}). We claim
\begin{Cl}\label{6}
With the choice of $w_\#$ given in (\ref{gatow}), 
\be\label{eq:dS1}
\| {\bf III}(w_{\#}) \|_{H^1(\R)}\leq K \abs{d(\ve)}c^{1+\frac 2{m-1}}.
\ee
and
\be\label{eq:dS2}
\norm{[ f(\tilde u + w_\#) -f(\tilde u) -f'(Q)w_\#]_x}_{H^1(\R)}  \leq K\abs{d(\ve)} c^{\frac 12 + \frac{3}{m-1}}
\ee
\end{Cl}

\begin{proof} 
The proof is similar to the proof of Proposition \ref{prop:decompbis} above. We only sketch the main ideas. 

Let us prove (\ref{eq:dS1}). First, note that for $\bar P$ defined in (\ref{barp})
\be\label{anul}
(\mathcal{L} (1+ \bar P))' = (1-f'(Q)+f'(Q))' =0.
\ee
This property will be useful in what follows. From the calculations performed in (\ref{B3}), (\ref{anul}) and the fact that $(1+\bar P)' \in \mathcal{Y}$, we note that (cf. (\ref{eq:sion2}) and (\ref{Sigm}) for the definition of ${\bf III}(\cdot)$ and $\Sigma_m'$ respectively)
\bee
{\bf III}((1+\bar P)(Q_c^2)' ) & = & -(\mathcal L (1+\bar P))' (Q_{c}^2)' -c (1+\bar P)(Q_{c}^2)'' +(1+\bar P)(Q_{c}^2)^{(3)} \\
& & + \sum_{(k,l)\in \Sigma'_{m}} c^l [\tilde F_{k,l}Q_{c}^k + \tilde G_{k,l}(Q_{c}^k)' ] + O(cQ_{c}^3 +Q_{c}^5 + c^2Q_{c})\\
& = &  \sum_{(k,l)\in \Sigma'_{m}} c^l [ \tilde F_{k,l}Q_{c}^k + \tilde G_{k,l}(Q_{c}^k)' ] +  O(c(Q_{c}^2)' + cQ_{c}^3 +Q_{c}^5 + c^2Q_{c}),
\eee
 where both $\tilde F_{k,l}$ and $\tilde G_{k,l}$ are in $\mathcal{Y}$. Moreover, $\tilde F_{3,0}=0$. From here, the definition of $w_\#$ in (\ref{gatow}) and Claim \ref{cl:30}, we obtain
 $$
 \norm{{\bf III}(w_{\#})}_{H^2(\R)} \leq K \abs{d(\ve)} c^{\frac 2{m-1}+1}.
 $$

Now, we deal with (\ref{eq:dS2}). 
We note that 
\bee
(\ref{eq:dS2}) & = & [ f(\tilde u + w_{\#}) -f(\tilde u) -f'(\tilde u)w_{\#} ]_{x}  + [ (f'(\tilde u) - f'(Q))w_{\#} ]_{x} \\
& & [ \frac 12 f''(\tilde u)w_{\#}^2 + O(w_{\#}^3) ]_{x} + [ f''(Q)(\tilde u -Q)w_{\#} + \frac 12 f^{(3)}(Q)(\tilde u -Q)^2 w_{\#} + O( (\tilde u -Q)^3w_{\#}) ]_{x}. 
\eee
From here, using the expresion for $w_{\#}$  and Claim \ref{cl:30}, we obtain
$$
\norm{(\ref{eq:dS2})}_{H^2} \leq K \abs{d(\ve)} [ \abs{d(\ve)}c^{\frac 32 + \frac{4}{m-1}} + c^{\frac 12 + \frac{4}{m-1}} + c^{\frac 12 + \frac{3}{m-1}} ] \leq K \abs{d(\ve)} c^{\frac 12 + \frac{3}{m-1}}.
$$
This finishes the proof.

\end{proof}

This Claim allows us to finish the proof of the Proposition.

\end{proof}

\bigskip

%\begin{rem}\label{Sch}
%We give here an intuitive explanation of the non-appareance of small solitons on the left of the two solitons after the collision. We believe that this phenomenon could be related to the appearance of linearly growing functions as solutions of the linear system $(\Omega_{1,0})$, or more generally to the existence of solutions not satisfying the {\bf (IP)} property for some linear system $(\Omega_{k,l})$ with $(k,l)\in \Sigma_m$; and this is indeed not the case. 

%As an example, suppose that for the linear system $(\Omega_{1,0})$ every polynomially growing solution satisfies $A_{1,0}\in \mathcal Y$, $B_{1,0}(x) - b_{1,0}\varphi + 2x \in \mathcal{Y} $. Then, as in Remark \ref{Sch0}, at time $t= T_c$ and discarding localized terms in $\mathcal{Y}$,
%\bee
%\tilde u(T_c,x) & \sim & Q(x) + Q_c(y_c) + b_{1,0} Q_c'(y_c) - 2xQ'(y_c)\\
%& \sim & Q(x) + Q_c(y_c + b_{1,0})  + Q_c(x-cT_c-2x) \\
%& \sim & Q(x) + Q_c(y_c + b_{1,0})  + Q_c(x+cT_c).
%\eee
%Although extremely informal this computation suggests the existence of a small soliton at time $t=T_c$ going to the left on the line $x=-ct$. However, the resolution of the linear systems $(\Omega_{k,l})$ satisfying the {\bf (IP)} property finally rules out this scenario at this order.     
%\end{rem}

%\end{document}
%
\section{Preliminary results for stability of the 2-soliton structure}\label{sec:3}

In this section several stability results will allow to study the long time behavior of the 2-soliton soliton solution. First of all, we recall a general result proved in \cite{MMcol2} concerning the existence and properties of an actual function $u=u(t,x)$, solution of (\ref{gKdV}) in the interval $[-T_c,T_c]$ and close enough to our approximate solution $\hat u$.
This will be done in the next subsection.

Next, we study the stability of a solution $u(t)$ of (\ref{gKdV0}) for long time, namely $t\geq T_c$. These results have been proved in great generality by Martel and Merle in \cite{MMas1}, \cite{MMcol2}, and \cite{MMas2}. In particular, we will use the \emph{stability} and \emph{asymptotic stability} of the two solitons (Proposition \ref{ASYMPTOTIC}) to show the persistence of the 2-soliton structure for long time.

Finally, a key result is the \emph{decomposition result} from Lemma \ref{LEMMEB1}, which will be essential to show the persistence of the residual term (cf. (\ref{eq:zT})) at infinity.

\subsection{Dynamic stability in the interaction region}\label{sec:3.1}

For any $c>0$ sufficiently small, we will consider the function $\hat u(t)$ of the form
\begin{equation*}
\hat u(t,x)=Q(y)+\wqs(y_c)+\sum_{(k,l)\in \Sigma_m} c^l \left\{\wqs^k(y_c) \hat A_{k,l}(y)+(\wqs^k)'(y_c) \hat B_{k,l}(y)\right\}
\end{equation*}
defined in Proposition \ref{cor:1} (the notation was introduced in (\ref{defv}) and (\ref{defW})).  
Recall the error term
$$S[\hat u](t)=\hat u_t +(\hat u_{xx} - u+ f(\hat u) )_x.$$

%%%%%%%%%%%%%%%%%%%%%%%%%%%%%
%%%%%%%%%%%%%%%%%%%%%%%%%%%%%
%

\begin{prop}[Exact solution close to the approximate solution $\hat u$, \cite{MMcol2}]\label{prop:I}
Let $\theta>\frac 1{m-1}$ and $\ve$ small enough such that (\ref{stab}) holds for $Q$. There exists $c_0>0$ such that the following holds for any $0<c <c_0$.
    Suppose that for all $ t\in [-T_c,T_c]$
 \begin{equation}\label{INTkl}
        \left\| S[\hat u](t)\right\|_{H^2(\R)}\leq K \frac {c^\theta}{T_c},    
 \end{equation}
   and   for some $T_0\in [-T_c,T_c]$,
    \begin{equation}\label{hypINT}
        \| u(T_0) - \hat u(T_0) \|_{H^1(\R)}\leq K c^{\theta},
    \end{equation}
where  $u(t)$ is an $H^1$ solution of (\ref{gKdV}). 
    Then, there exist $K_0=K_0(\theta,K,f)$ and a $C^1$ function $\rho:[-T_c,T_c]\rightarrow \R$ such that, for all $t\in [-T_c,T_c]$,
    \begin{equation}\label{INT41}
        \|u(t)-\hat u(t, \cdot-\rho(t)) \|_{H^1(\R)} \leq K_0 c^{\theta},\quad |\rho'(t) -1|\leq K_0 c^{\theta}.
    \end{equation}
\end{prop}

\smallskip

\begin{rem}
The proof of the above Proposition is nontrivial and requires some refined techniques such as  modulation theory, coercivity properties and the introduction of a modified energy functional adapted to a two soliton collision. It is necessary to emphasize that one of the key elements in the proof is the smallness of the error term $S[\hat u]$ along the collision. For the sake of completeness, we will draw the main lines of the argument, see \cite{MMcol2} for the actual complete proof.
\end{rem}

\begin{proof}
It suffices to show the result on the interval  $[T_0,T_c]$. By using the transformation $x\to -x$, $t\to -t$, the proof is the same on $[-T_c,T_0]$.

Let $K^*>1$ be a constant to be fixed later. Since $\|u(T_0)-\hat u(T_0)\|_{H^1(\R)}\leq c^\theta$, by continuity in time in $H^1(\mathbb{R})$, there exists $T_0<T^*\leq T_c$ such that
\bee
    T^*& = & \sup\big\{T\in [T_0,T_c] \hbox{ such that for all }  t\in [T_0,T], \hbox{ there exists } r(t)\in \mathbb{R} \hbox{ with } \\
    & & \quad \quad \quad     \|u(t) -  \hat u(t, \cdot  - r(t))\|_{H^1(\R)}\leq K^* c^{\theta} \big\}.
\eee
The objective is to prove that $T^*=T_c$ for $K^*$ large. For this, we argue by contradiction, assuming that $T^*<T_c$ and reaching a    contradiction with the definition of $T^*$ by proving some independent estimates on $\|u(t)- \hat u(t,\cdot -r)\|_{H^1(\R)}$ on $[T_0,T^*]$.

An argument using the Implicit function theorem allows to construct a modulation parameter and to estimate its variation in time:
 
\begin{Cl}\label{DEFZ} Assume that $0<c<c(K^*)$ small enough.
There exists a unique $C^1$ function $\rho(t)$ such that, for all $t\in [T_0,T^*]$,
\begin{equation*}\label{defz}
z(t,x)=u(t,x+\rho(t))-\hat u(t,x) \quad \text{satisfies}\quad 
 \int_\R z(t,x) Q'(y) dx=0.
\end{equation*}
Moreover, we have,  for all $t\in [T_0,T^*]$,
\begin{equation*}\label{TRANS3}
              |\rho(T_0)|+\|z(T_0)\|_{H^1(\R)}\leq K c^{\theta}, \ \|z(t)\|_{H^1(\R)}\leq  2 K^* c^{\theta},
\end{equation*}
\begin{equation*}\label{eqz}
 z_t +( z_{xx} -z + f(z+\hat u)-f(\hat u))_x= -S[\hat u](t) + (\rho'(t)-1) (\hat u+z)_x.
\end{equation*}
\begin{equation*}\label{TRANS3bis}
 |\rho'(t)-1|\leq K \|z(t)\|_{H^1(\R)}+K \|S[ \hat u] (t)\|_{H^1(\R)},
\end{equation*}
\end{Cl}

\medskip

The purpose of the modulation theory is to establish a lower bound in the following energy functional for $z(t)$:
$$\mathcal{F}(t): =\frac 12 \int_\R \left((\partial_x z)^2 + (1+\alpha'(y_c))z^2\right)
- \int_\R (F(\hat u+z) - F(\hat u) - f(\hat u) z).$$

Indeed, this functional enjoys two useful properties: it has a very small time variation and it is coercive up to the direction $Q$:

\begin{lem}[Coercivity of $\mathcal{F}$]\label{varF}
Assume that $0<c<c(K^*)$ small enough.
    There exists $K>0$ (independent of $K^*$ and $c$) such that
    \begin{enumerate}
    \item Coercivity of $\mathcal{F}$ under orthogonality conditions:
    \begin{equation*}\label{posf1}
       \forall t\in [T_0,T^*],\quad\|z(t)\|_{H^1(\R)}^2\leq K \mathcal{F}(t) + K \left|\int_\R z(t) Q(y)\right|^2.
    \end{equation*}
    \item Control of the direction $Q$:
    \begin{equation*}\label{Qdir1} \forall t\in [T_0,T^*],\quad
        \left|\int_\R  z(t)Q(y)\right|\leq  K c^\theta + K c^{\frac 1{p-1} -\frac 14} 
        \|z(t)\|_{L^2}+ K\|z(t)\|_{L^2}^2.
\end{equation*}
\item Control of the variation of the energy functional:
\begin{equation*}\label{varF1}
        \mathcal{F}(T^*)-\mathcal{F}(T_0)\leq  K c^{2 \theta}\left((K^*)^2(1+K^*)c^{\frac 1{2(m-1)} -\frac 18} + K^*\right).
    \end{equation*}
\end{enumerate}
\end{lem}

\medskip

These estimates allow us, after fixing $K^*$ large enough and possibly taking $c$ even smaller, to show that actually 
\begin{equation*}
    \|z(T^*)\|_{H^1(\R)}^2 \leq \frac 12 (K^*)^2 c^{2 \theta}.
\end{equation*}
contradicting the definition of $T^*$, thus proving that $T^*=T_c$.
\end{proof}

%%%%%%%%%%%%%%%%%%%%%%%%%%%%%
%%%%%%%%%%%%%%%%%%%%%%%%%%%%%

Once the existence of an actual solution (close to our approximate solution $\hat u$ in the interval $[-T_c, T_c]$) is established, one would like to investigate the behavior in long time of this solution. We treat this problem in the next subsection. 

\medskip

\subsection{Stability and asymptotic stability for large time}\label{sec:3.2}

Here we consider the stability of the $2$-soliton structure after the collision, and for a long time.  Let $T _c$ be defined in (\ref{Tc}). We start with an important

\begin{rem}\label{stabeps}
Since (\ref{stab}) holds for $f(s)=s^m$, $m=2,3$, it is clear by a perturbation argument that (\ref{stab}) holds also for $f$ as in (\ref{gKdV0b}) for all $0<c<1$, provided $0<\abs{\ve}<\ve_{0}$ is small enough. 
\end{rem}

\begin{prop}[Stability of two decoupled solitons, \cite{MMas1}, \cite{MMas2}]\label{ASYMPTOTIC}

Let $\ve$ small enough such that (\ref{stab}) holds for $Q$. Then there exist constants $c_0, K>0$, such that for any $0<c<c_0$ and for any $\omega>0$, the following holds.
Let $u(t)$ be an $H^1$ solution of (\ref{gKdV}) such that for some time $t_1\in \R$ and  $\frac 12 T_c \leq X_0 \leq \frac 32 T_c$,

\begin{equation}
\label{D25}
\|u(t_1)-Q -Q_{c}(\cdot + X_0)\|_{H^1(\R)}\leq c^{\frac 14 + \frac{1}{m-1} +\omega}.
\end{equation}

Then there exist $C^1$-functions $\rho_1(t)$, $\rho_2(t)$ defined on $[t_1,+\infty)$  such that
\begin{enumerate}

\item Stability:
\begin{equation}\label{huit}
\sup_{t\geq t_1} \| u(t)- Q(\cdot -\rho_1(t))- Q_{c}(\cdot -\rho_2(t))  \|_{H^1(\R)}
\le K c^{-\frac 14 + \frac{1}{m-1} +\omega},
\end{equation}
and for all $t\geq t_1$,
\begin{equation}\label{suppl}%\begin{split}
 \frac 12 \leq \rho_1'(t)-\rho_2'(t) \leq \frac 32, \quad  | \rho_1(t_1)|\leq K c^{\frac 14 + \frac 1{m-1}+ \omega },\quad |\rho_2(t_1) + X_0|\leq K c^{\omega}.
%\end{split}
\end{equation}

\item Asymptotic stability:
There exist $c_1^+, c_2^+ > 0$  such that on the right hand side limit
\begin{equation}\label{neuf}
\lim_{t\to +\infty}\|  u(t)- Q_{c_1^+}(x-\rho_1(t)) -Q_{c_2^+}(x-\rho_2(t)))\|
_{H^1(x> \frac 1{10} c t)} = 0,
\end{equation}
with
\begin{equation}\label{sept}
\big| c_1^+  - 1\big| \leq   K  c^{\frac 14+ \frac 1{m-1} + \omega } ,\quad
 \Big| \frac{c_2^+}{c}-1 \Big|  \leq K  c^\omega .
\end{equation}
\end{enumerate}
\end{prop}

\subsection{A decomposition result}
Recall a more precise decomposition of $u(t)$ used in the proof of Proposition \ref{ASYMPTOTIC} in
\cite{MMas1}, \cite{MMas2}.

\begin{lem}[Decomposition of the solution, \cite{MMas2}]\label{LEMMEB1}
Suppose (\ref{stab}) holds for $Q$. Let $u=u(t)$ be a solution of the gKdV equation (\ref{gKdV0}) such that the estimate (\ref{D25}) holds. Then there exist $C^1$-functions $\rho_1(t)$, $\rho_2(t)$, $c_1(t),$ $c_2(t)$, defined on $[t_1,+\infty)$,  such that the function
\begin{equation*}
\eta(t,x):= u(t,x) - R_1(t,x)-R_2(t,x),
\end{equation*}
where, for $ j=1,2,$ $R_j(t,x):=Q_{c_j(t)}(x-\rho_j(t)),$ satisfies for all $t\geq t_1,$
\begin{eqnarray}&&
\int_\R R_j(t)\eta(t) =\int_\R (x-\rho_j(t))  R_j(t)\eta(t)=0,\quad j=1,2,\label{dix}\label{eq:o1}\\ &&
\|\eta(t)\|_{H^1(\R)}+\left| c_1(t)-1\right|
+c^{\frac 1{m-1}-\frac 14}\left|\frac {c_2(t)}{c} -1\right|
\le K   c^{\omega+\frac 1{m-1}-\frac 14}, \label{onze}\\
& & \hbox{ and for all } t\geq t_1 \quad  |\rho'_2(t)| + |\rho_1'(t)-1|\leq \frac 1{10}, \quad \rho_1(t)-\rho_2(t) \geq \frac 12 t + \frac 14T_c.
\end{eqnarray}
Moreover, we have the convergence $\lim_{t\to +\infty} \bar c_j(t) =c_j^+$ for $j=1,2$.
\end{lem}

At this moment we have all the necessary information about the 2-soliton solution of (\ref{gKdV}). Indeed, recall from the sketch of proof (Subsection \ref{sop}) that the asymptotic in long time will be treated using the tools from this section, more precisely using Proposition \ref{ASYMPTOTIC} and Lemma \ref{LEMMEB1}. On the other hand the collision region will be described by Proposition \ref{prop:I}. This is the purpose of the next section.   

\medskip

\section{Proof of the Theorem \ref{PUREgen}}

Now we are in a position to prove the main Theorem of this work.

\medskip

\begin{proof}[Proof of Theorem \ref{PUREgen}]
Let $1=c_1<c_*(f)$ such that \eqref{stab} holds and $0<c<c_0(\ve)$ small enough (depending on $\ve$).
Let $u(t)$ be the \emph{unique} solution of (\ref{gKdV0}) such that (see Theorem 1 and Remark 2 in \cite{Martel})
$$
\lim_{t\to -\infty} \|u(t)- Q (x - t)-Q_{c}(x- ct)\|_{H^1(\R)}=0.
$$

\medskip

\noindent{\bf 1. Behavior at $-T_c$.}
We claim that for all   $t< - \frac 1{32} T_c$, 
\begin{equation}\label{ETI1}
\|  u(t)-Q (\cdot - t) - Q_{c}(\cdot -c t)\|_{H^1(\R)} \leq K e^{\frac 14 \sqrt{c}(1-c)t}.
\end{equation}
This is a consequence of the proof of existence of $u(t)$ in \cite{Martel}. 
See Proposition 5.1 in \cite{MMcol1} for a proof in the power case.

\smallskip

Now, using (\ref{ETI1}), we will match the function $u$ with the collision solution $\hat u$ constructed in Proposition \ref{cor:1}. For this, we will translate $u$ in time and space, as follows.

Let $\Delta_1$, $\Delta_2$ be defined in Proposition \ref{lem:z} and
$$
T_c^- := T_c + \frac 12 \frac {\Delta_1-\Delta_2}{1 - c},\quad
a :=\frac 12 \Delta_1 - T_c^- .
$$
Since from (\ref{eq:Dt}) 
$$
|\Delta_1|\leq K c^{\frac 1{m-1} -\frac 12},\quad  \hbox{ and } \quad  \abs{\Delta_2} \leq K,
$$
we have $-T_c^-  \leq -\frac 1{32} T_c$, and thus, from (\ref{ETI1}) for $c$ small enough, and after a translation by $a$, we get 
\begin{equation*}
\|u(-T_c^-, \cdot + a )-Q(\cdot +\frac {\Delta_1}2)-Q_{c}(\cdot -(1-c)T_c +\frac {\Delta_2}2)\|_{H^1(\R)} 
\leq K e^{- \frac 14 \sqrt{c}(1- c)T_c^-}\leq K c^{10}.
\end{equation*}
%Let $\bar u(t,x) := u( t +T_c-T_c^-, x + a)$. Then $\bar u(t,x)$ is also solution of \eqref{gKdV0} and satisfies
By translation invariance, we may assume $T_c^- =T_c$ and $a=0$, such that
\begin{equation}\label{atmoinstc}
\| u(-T_c)-Q (\cdot +\frac {\Delta_1}2)
-Q_{c}(\cdot -(1-c)T_c +\frac {\Delta_2}2)\|_{H^1(\R)} \leq  K c^{10}.
\end{equation}

\medskip

\noindent{\bf 2. Behavior at $+T_c$.}  Now, possibly taking a smaller $c$, consider $\hat u=\hat u_{1,c}$ constructed in Proposition \ref{cor:1}. 
By \eqref{eq:zmT} and \eqref{atmoinstc}, we have
$$
\|u(-T_c) - \hat u(-T_c)\|_{H^1(\R)} \leq K c^{\frac 2{m-1} +\frac 14} [c^{\frac 1{m-1}} + \abs{d(\ve)}c^{\frac 12}].
$$
Applying Proposition \ref{prop:I} with 
$$T_0=-T_c,\quad c^\theta := c^{\frac 2{m-1} +\frac 14-\frac 1{100}} [c^{\frac 1{m-1}} +\abs{d(\ve)}c^{\frac 14}],$$
it follows that there exists a function $\rho(t)$ such that for all $t\in [-T_c,T_c]$,
$$
\|u(t)-\hat u(t,\cdot -\rho(t))\|_{H^1(\R)} \leq   K c^\theta.
$$
In particular, for $r:=\rho(T_c)$, we have 
$$
\|u(T_c)-\hat u(T_c,\cdot -r)\|_{H^1(\R)} \leq K c^\theta.
$$
Using (\ref{eq:zT}) and triangular inequality , we obtain
\begin{equation}\label{atplustc}
\| u(T_c)- Q (\cdot -r_1) -Q_{c}(\cdot -r_2)  -2d(\ve)(Q_{c}^2)' (\cdot -r_2) \|_{H^1(\R)} \leq K c^{\theta}.
\end{equation}
Here, 
$$
r_1 :=  \frac 12 \Delta_1 + r, \quad r_2:= (c-1) T_c +\frac 12 \Delta_2 + r,
$$
so that $r_1-r_2 = (1-c)T_c +\frac 12 (\Delta_1-\Delta_2)$ satisfies
\be\label{r1mr2}
\frac 12 (1 -c) T_c \leq r_1-r_2 \leq 32 (1-c)T_c.
\ee
Moreover, note that
$$
\norm{ (Q_{c}^2)'  }_{H^1(\R)}\leq Kc^{\frac{2}{m-1} +\frac 14 },
$$
so that 
\be\label{diff1}
\| u(T_c)- Q(\cdot - r_1) -Q_{c}(\cdot - r_2)  \|_{H^1(\R)} \leq K c^{\frac 2{m-1} + \frac 14}[\abs{d(\ve)} + c^{\frac 1{m-1}}] \leq K \abs{d(\ve)}c^{\frac 2{m-1} +\frac 14},
\ee
provided
\be\label{cond1}
\abs{d(\ve)}\geq \kappa_0 c^{\frac 1{m-1}-\frac 1{100}},
\ee
for some $\kappa_0>0$ large enough but fixed. We have thus arrived to time $t=+T_c$ with a stability property of the 2-soliton structure, namely (\ref{diff1}).

\medskip

\noindent{\bf 3. Behavior as $t\to +\infty$.}
From \eqref{diff1}, it follows that we can apply Proposition \ref{ASYMPTOTIC} to
$u(t, \cdot + r_1)$ for $t\geq T_c$ (that is, $t_1 :=T_c$), with $X_0:=r_1-r_2$, and
$$c^\omega:= \abs{d(\ve)} c^{\frac 1{m-1}}.$$

It follows that there exist $\rho_1(t)$, $\rho_2(t)$, $c_1^+>0$, $c_2^+>0 $ so that
\be\label{wp}
w^+(t,x) := u(t,x ) -Q_{c_1^+}(x-r_1-\rho_1(t))-Q_{c_2^+}(x-r_1-\rho_2(t))
\ee
satisfies
\be\label{ETI5}
 \sup_{t\geq T_c} \|w^+(t)\|_{H^1(\R)} \leq  K  \abs{d(\ve)} c^{\frac 2{m-1}-\frac 14} , \quad \lim_{t\to +\infty} \|w^+(t)\|_{H^1(x>\frac {c}{10} t +r_1)} =0,
\ee
and
\be\label{ETI3}
|c_1^+ - 1|\leq K \abs{d(\ve)} c^{\frac 2{m-1}+\frac 14},\quad |c_2^+ - c|\leq K  \abs{d(\ve)} c^{1+\frac 1{m-1}}.
\ee

\medskip

In particular, the behavior of the 2-soliton structure remains stable at infinity, modulo the emergency of a possible $H^1$-nonzero residual term. This proves the upper bound in (\ref{restgen}). At this stage, we do not know if this residual term (that is, $w^+$) can be bounded by below uniformly in time. This is the purpose of the following key step.

\medskip

%%%%%%%%%%%%%%%%%%%%%%%%%%%%%%%%%%%
%%%%%%%%%%%%%%%%%%%%%%%%%%%%%%%%%%%

\noindent
{\bf 4. Lower bound on $w^+(t)$ for $t>T_c$ large.}
Consider the decomposition of $u(\cdot,\cdot + r_1)$ defined in Lemma ~\ref{LEMMEB1}, i.e.
the center of mass $\bar\rho_1(t)$, $\bar\rho_2(t)$, the scaling parameters
$\bar c_1(t)$, $\bar c_2(t)$ such that, for $t>T_{c}$,
\be\label{etaa}
\eta(t,x) := u(t,x) - Q_{\bar c_1(t)}(x-r_1 -\bar \rho_1(t) ) -Q_{\bar c_{2}(t)}(x-r_1 -\bar \rho_2(t) )
\ee
satisfies
\begin{equation}\label{eq:tt}\begin{split}
        &\sup_{t\geq T_{c}}\|\eta(t)\|_{H^1(\R)}\leq K \abs{d(\ve)} c^{\frac 2{m-1} - \frac 14},\quad
        |\bar c_1(T_{c})-1| \leq K \abs{d(\ve)} c^{\frac 2{m-1} -\frac 14},\\
        & \bar\rho_1(t)-\bar \rho_2(t) \geq \frac 12 t +\frac 14 T_c ,\quad 
        |\bar c_2(t)-c |\leq K \abs{d(\ve)} c^{\frac 1{m-1}+\frac 32}, 
\end{split}\end{equation}
and
\be\label{last}
\abs{\bar\rho_2(T_c) + r_1 - r_2} \leq  K \abs{d(\ve)} c^{\frac 1{m-1}}.
\ee
Moreover, we have for $j=1,2$
\be\label{conv}
\lim_{t\to +\infty} \bar c_j(t) = c_j^+.
\ee

\smallskip

First, as a consequence of \eqref{atplustc}, we claim the following lower bound at $t=T_{c}$: for $K_0>0$, independent of $c>0$,
\begin{equation}\label{eq:aT}
        \int_{x<\bar\rho_2(T_{c})+r_1+\frac 14 T_{c}} \eta^2(T_{c},x) dx \geq K_0  \abs{d(\ve)}^2 c^{\frac 4{m-1} + \frac 12}.
\end{equation}

\emph{Proof of \eqref{eq:aT}}. The proof will proceed by a contradiction argument. Indeed, suppose that for any $\alpha >0$ there exists $c>0$ small enough such that (\ref{eq:aT}) does not hold properly, namely
\be\label{eqcontr}
\norm{\eta(T_{c})}_{L^2(x<\bar\rho_2(T_{c}) +r_1+\frac 14T_{c})}\leq \alpha \abs{d(\ve)} c^{\frac 2{m-1} + \frac 14}.
\ee
 
Replacing $$u(T_{c},x) =Q_{\bar c_1(T)}(x-r_1-\bar\rho_1(T_{c})) + Q_{\bar c_2(T_{c})}(x-r_1-\bar\rho_2(T_{c}))+\eta(T_{c},x)$$ in \eqref{atplustc}, we find
\begin{equation*}\begin{split}
& \big\|[ Q_{\bar c_1(T_{c})}(\cdot-r_1 -\bar\rho_1(T_{c}))-Q (\cdot -r_1)] +[Q_{\bar c_2(T_{c})}(\cdot -r_1-\bar\rho_2(T_{c}))-Q_{c}(\cdot -r_2)]\\ 
& \quad + \eta(T_{c})+2d(\ve) (Q_{c}^2)'(\cdot-r_2)\big\|_{H^1(\R)} 
 \leq K c^{\frac 2{m-1} + \frac 14}[ c^{\frac 1{m-1}} + \abs{d(\ve)}c^{\frac 12} ].
\end{split}\end{equation*}

By the decay properties of $Q$, (\ref{eq:tt}) at time $t=T_c$ and $r_1-r_2\geq \frac 12 (1-c) T_{c}$ (see (\ref{r1mr2})), we obtain
\bee
   \|[Q_{\bar c_2(T_{c})}(\cdot-r_1- \bar\rho_2(T_{c}))-Q_{c}(\cdot -r_2)]
 + \eta(T_{c})+2d(\ve)  (Q_{c}^2)'(\cdot - r_2)\|_{L^2(x<\bar \rho_2(T_{c})+ r_1+\frac 14 T_{c})} 
    \quad  \quad   & & \\
 \leq   K c^{\frac 2{m-1} + \frac 14}[ c^{\frac 1{m-1}} + \abs{d(\ve)}c^{\frac 12} ]. & & 
\eee

Then, using (\ref{eqcontr}) and (\ref{cond1}),
\bee
  \|[Q_{\bar c_2(T_{c})}(\cdot -r_1- \bar\rho_2(T_{c}))-Q_{c}(.-r_2)]
  +2d(\ve)  (Q_{c}^2)'(\cdot -r_2)\|_{L^2(x<\bar \rho_2(T_{c})+r_1+\frac 14 T_{c})} \quad \quad \quad && \\
 \leq   c^{\frac 2{m-1} + \frac 14}\big[ Kc^{\frac 12} + 2\alpha + \frac K{\kappa_0}c^{\frac 1{100}} \big] \abs{d(\ve)}. & &
\eee
By scaling and translation, and decay of $Q$, we obtain
\bee
  \|\bar Q-Q + 2 d(\ve) c^{\frac 12 + \frac1{m-1}}  (Q^2)'\|_{L^2(\R)} & \leq & c^{\frac 1{m-1} + \frac 12}\big[ Kc^{\frac 12} + 2\alpha + \frac K{\kappa_0}c^{\frac 1{100}} \big]  \abs{d(\ve)}  \\
   & & \quad + \|\bar Q-Q + 2 d(\ve) c^{\frac 12 + \frac1{m-1}}  (Q^2)'\|_{L^2(x>\beta)},
\eee
where $\bar Q(x)= \lambda Q\left( \mu x - \xi\right),$ and
$$
\lambda:=\Big[\frac {\bar c_2(T_{c})}{c}\Big]^{\frac 1{m-1}}, \quad \mu:=\sqrt{\frac{\bar c_2(T_{c})}{c}},
$$
(do not be confused with $\mu$ of Theorem \ref{PUREgen}), and
$$
 \xi= \sqrt{\bar c_{2}(T_{c})}(\bar\rho_2(T_{c}) +r_1-r_2), \quad \beta := \sqrt{c}(\frac 14 T_{c} +\bar \rho_{2}(T_{c}) +r_1-r_{2}).
$$
Note that from (\ref{eq:tt}) and (\ref{last}), 
$$
\beta\geq \frac 18 \sqrt{c} T_c\geq \frac 18 c^{-\frac 1{100}}  , \quad  \|\bar Q-Q + 2 d(\ve) c^{\frac 12 + \frac1{m-1}}  (Q^2)'\|_{L^2(x>\beta)}\leq Kc^{10}.
$$
Moreover, note that $\bar Q (x) = Q_{\mu}(x-\frac{\xi}{\mu})$, and that by \eqref{eq:tt}, we have
$$
|\mu-1|\leq K \abs{d(\ve)}c^{\frac 12 +\frac 1{m-1}},\quad  |\xi|\leq K \abs{d(\ve)}c^{\frac 12 + \frac 1{m-1}}.
$$
Expanding $\bar Q$ in $\mu-1$, and $\xi /\mu$, and using parity properties, we find
$$
\| \xi Q'  + 2 d(\ve)c^{\frac 12 + \frac 1{m-1}}  (Q^2)' \|_{L^2(\R)} \leq \big[ Kc^{\frac 12} + 3\alpha + \frac K{\kappa_0}c^{\frac 1{100}} \big]  \abs{d(\ve)} c^{\frac 12 +\frac 1{m-1}},
$$
so that for some constant $\bar \xi\in \R$,
$$
\| \bar \xi Q'  + 2 d(\ve)(Q^2)' \|_{L^2(\R)} \leq \big[ Kc^{\frac 12} + 4\alpha + \frac K{\kappa_0}c^{\frac 1{100}}\big]  \abs{d(\ve)}.
$$
Note that exists $\kappa_{1}>0$, independent of $\ve$ and $c$, such that 
$$
\inf_{\bar \xi \in \R} \| \bar \xi Q'  + 2 d(\ve)(Q^2)' \|_{L^2(\R)} \geq \kappa_{1} \abs{d(\ve)},
$$  
since $Q'\neq \gamma (Q^2)' $ for all $\gamma \in \R$. By choosing $\kappa_0$ large enough in (\ref{cond1}), depending only on $\kappa_{1}$, and $c$ small enough, we find a contradiction for $\alpha$ small. This contradiction proves \eqref{eq:aT}.

\medskip

Now, we finish the proof of the lower bound by proving the following
\begin{lem}\label{5.1}
There exists $K_0>0$ such that 
\begin{equation}\label{eq:lw}
        \liminf_{t\to +\infty} \|w^+(t)\|_{H^1_{c}(\R)} \geq K_0 \abs{d(\ve)}c^{\frac {3}4+ \frac 2{m-1}}.
\end{equation}
\end{lem}

Note that \eqref{eq:lw} combined with (\ref{diff0})-(\ref{diff2}) prove the lower bound in (\ref{restgen}). Thus, we are now reduced to prove \eqref{eq:lw}.

\begin{proof}

We argue by contradiction. Assume that for any $\alpha>0$, there exist
arbitrarily large $T_0$ and $c$ arbitrarily close to $0$ such that
\begin{equation}\label{eq:ct}
\|w^+(T_0)\|_{H^1_{c}(\R)} \leq \alpha  \abs{d(\ve)} c^{\frac {3}4+ \frac 2{m-1}}.
\end{equation}
By (\ref{conv}), we can choose $T_0>T_{c}$ large enough so that
\begin{equation}\label{eq:yy}
        \|\eta(T_0)\|_{H_{c}^1(x < m(T_0)+\frac{T_0}4)} \leq 2 \alpha  \abs{d(\ve)} c^{\frac {3}4+ \frac 2{m-1}}.
\end{equation}
Here $m(t) :=r_1 + \frac 12 (\bar\rho_1(t)+\bar\rho_2(t))$  is the \emph{middle point} between the two solitons at time $t$.

We need to estimate some local in space conservation laws. For this reason we introduce a sort of cutoff function \emph{supported} on the small soliton. Let
\begin{equation}\label{eq:ph}\begin{split}
        & \psi(x)=\frac 2 \pi \arctan(\exp(x/\kappa)), \quad \text{so that} \quad 
        \lim_{-\infty} \psi=0,\ \lim_{\infty} \psi=1, \hbox{ and for all } x\in \R,\\
        & \psi(-x)=1-\psi(x),\quad
 \psi'(x)=\frac 1{\pi\kappa \cosh(x/\kappa)},
 \quad  |\psi'''(x)|\leq \frac 1{\kappa^2} |\psi'(x)|.
%        &  m(t) :=\frac 12 (\rho_1(t)+\rho_2(t)).
\end{split}\end{equation}

Let 
$$
a:= \frac {E(Q_{\bar c(T_0)}) - E(Q_{\bar c(T_{c})})}
{M(Q_{\bar c(T_{c})}) - M(Q_{\bar c(T_0)})}.
$$
We set 
\bee%\label{eq:vq}
\mathcal{G}(t) &: = &\frac 12 a \int_\R u^2(t,x) (1-\psi(x-m(t))) dx + \frac 12 \int_\R (u_x^2- 2F(u))(t,x) (1-\psi(x-m(t))) dx \\
& = & a M(u(t)) + E(u(t)) - (a \mathcal{M}_1(t) + \mathcal{E}_1(t)),
\eee
where 
$$
\mathcal{M}_1(t) := \frac 12 \int_\R u^2(t,x)  \psi(x-m(t)) dx,\quad  \mathcal{E}_1(t) := \frac 12 \int_\R (u_x^2- 2F(u))(t,x) \psi(x-m(t)) dx.
$$
We claim the following results on $m(t)$, $a$ and $\mathcal{G}(t)$.

\begin{Cl}\label{Clf}
The following estimates hold
\be\label{m1}
\frac 12 \leq m'(t) \leq  \frac 32.
\ee
and for a positive constant $k_m$, 
\be\label{a1}
a = k_mc + o(c).
\ee
(Here $o(c)$ means $\abs{c^{-1}o(c)}\to 0$ as $c\to 0$.)
\end{Cl}

\begin{proof}
To prove (\ref{m1}), it is enough to consider Lemma \ref{LEMMEB1} on the interval $[T_c, T_0]$ to have
$$
m'(t)\geq 1-\frac 1{10} \geq \frac 12; \quad m'(t)\leq 1 + \frac 1{10}\leq \frac 32.
$$
Let us now treat (\ref{a1}). It is easy to show that 
$$
M(Q_{c}) = c^{\frac 2{m-1} -\frac 12} \int_{\R} Q^2 + o(c^{\frac 2{m-1} -\frac 12}).
$$  
On the other hand,
$$
E(Q_{c}) = \frac 12 c^{\frac 2{m-1} + \frac 12}\Big[ \int_{\R} Q'^2 - \frac 2{m+1} \int_{\R}Q^{m+1}\Big] + o(c^{\frac{m+1}{m-1} +\frac 12}).
$$
Thus, from (\ref{eq:tt}) and the fact that $E(Q)<0$, a Taylor expansion and L'Hopital rule gives
\be\label{asympa}
a = -\frac{\partial_c E(Q_c) }{\partial_c M(Q_c)}\Big|_{c=\bar c(T_c)} + O(\abs{\bar c(T_0)-\bar c(T_c)})  = -\frac{E(Q)}{M(Q)} \bar c(T_c) + o(c) = k_m c + o(c).
\ee
where $k_m:= -\frac{E(Q)}{M(Q)}>0$  is a constant depending on $m$. 
\end{proof}

\begin{lem}\label{le:mo}
For $0<c<c_0$ small enough,
$$
\mathcal{G}(T_{c})-\mathcal{G}(T_0) \leq K c^{10}.
$$
\end{lem}
\begin{proof}
We will need the following

\begin{Cl}
Define $h:=u_{xx}+f(u)$, such that $u_t = -h_x$. Then
$$
\mathcal{M}_1'(t)  =  -\frac 32 a\int_\R u_x^2 \psi' +\frac a2 \int_\R u^2 (\psi''' - m' \psi') + a\int_\R (uf(u)-F(u))\psi',
$$
and
\bee
\mathcal{E}_1'(t) & = & -\frac 32 \int_\R h^2 \psi' - \frac 12\int_\R u_x^2( m' \psi' -\psi''') + \int_\R F(u)(m' \psi' +\psi''')\\
& &  \quad +\int_\R f^2(u)\psi' -\int_\R u_x^2 f'(u)\psi' .
\eee
\end{Cl}

\begin{proof} A direct computation, see for example \cite{MMMcol}.
\end{proof}

Now, we follow the proof contained in Appendix D, \cite{MMMcol}. From above Claim, we have
\bee
 \mathcal{G}'(t) & = &  \frac 32 \int_\R h^2 \psi'  + \frac 12\int_\R u_x^2( m' \psi' -\psi''' +3a\psi' ) + \frac a2 \int_\R u^2 ( m' \psi' -\psi''' )\\
& & - a\int_\R (uf(u)-F(u))\psi' - \int_\R F(u)(m' \psi' +\psi''') -\int_\R f^2(u)\psi' + \int_\R u_x^2 f'(u)\psi' .
\eee
From  Claim \ref{Clf} we choose $\kappa>0$ large enough such that $m' \psi' -\psi'''  \geq \frac 1{4}\psi' $.  From here, 
$$
 \frac 32 \int_\R h^2 \psi'  + \frac 12\int_\R u_x^2( m' \psi' -\psi''' +3a\psi' ) + \frac a2 \int_\R u^2 ( m' \psi' -\psi''' )\geq \frac c4\int_\R (u_x^2 + u^2) \psi' 
$$
Let us consider now the nonlinear terms in the second row of $\mathcal{G}'(t)$. For this, let $$I:=[r_1+ \bar\rho_2(t) +\frac 18 T_c , r_1+ \bar\rho_1(t) -\frac 18 T_c]$$ an interval between the two solitons. We have two cases: $x\in I$ and $x\not\in I$.

In the first case, from (\ref{etaa}) we have for all $t\geq T_c$
$$
\abs{u(t)}\leq \abs{Q_{\bar c_1}}(t) + \abs{Q_{\bar c_2}}(t) +\abs{\eta}(t) \leq K \abs{d(\ve)} c^{\frac 2{m-1} - \frac 14}.
$$
Thus,
\bee
& & \abs{- a\int_I (uf(u)-F(u))\psi' - \int_I F(u)(m' \psi' +\psi''') -\int_I f^2(u)\psi' + \int_I u_x^2 f'(u)\psi' } \\
&  & \leq K \Big[\norm{u(t)}_{L^\infty(I)}^{m-1} + \norm{u(t)}_{L^\infty(I)}^{2(m-1)} \Big]\int_\R (u^2+ u_x^2)\psi' \\
& & \leq K \abs{d(\ve)}^{m-1} c^{2 - \frac 14(m-1)}  \Big[\abs{d(\ve)}^{m-1} c^{2 - \frac 14(m-1)} + 1\Big]\int_\R (u^2+ u_x^2)\psi'\\
& &  \leq K c^{\frac 32} \int_\R (u^2+ u_x^2)\psi' .
\eee
In the second case, we have $\abs{x-m(t)}\geq \frac t4$ and thus $\psi'(x-m(t))\leq K e^{- \ga t}$, with $\ga>0$ a fixed constant. From here,
\bee
\abs{- a\int_{x\not\in I} (uf(u)-F(u))\psi' - \int_{x\not\in I} F(u)(m' \psi' +\psi''') -\int_{x\not\in I} f^2(u)\psi' + \int_{x\not\in I} u_x^2 f'(u)\psi' } \leq  Ke^{-\ga t}.
\eee
In conclusion, putting together above estimates, we get for all $t\in [T_c, T_0]$,
$$
 \mathcal{G}'(t)\geq -Ke^{-\ga t},
$$ 
and after integration we obtain the desired result. The proof is now complete.
\end{proof}

Now, define
$$
\mathcal{H}(t) := \int_\R \big[ a\eta^2 + \eta_x^2 - f'(R_2) \eta^2\big] (1-\psi).
$$ 
We have the

\begin{lem}\label{le:qd}
For $0<c<c_0$ small enough,
\begin{enumerate}
\item Small variation:
\bea\label{sv}
 \mathcal{G}(T_c)-\mathcal{G}(T_0) & = &\frac 12 (\mathcal{H}(T_c)-\mathcal{H}(T_0))
+ O(\alpha^2 \abs{d(\ve)}^3 c^{\frac 6{m-1}+\frac 14 } c^{\frac{m-2}{m-1}})\nonumber  \\
& & \quad  + O( \abs{d(\ve)} c^{\frac 2{m-1}-\frac 14} c^{\frac{m-2}{m-1}}\int_\R \eta^2(T_c) (1-\psi) ) +O(c^{10}), 
\eea
\item Coercivity:
\be\label{coe}
\mathcal{H}(t) \geq  \sigma_0 \int_\R \left[c \eta^2 + \eta_x^2 \right](t,x) (1-\psi) dx.
\ee
for some $\sigma_0>0$ independent of $c$.
\end{enumerate}
\end{lem}

\begin{proof}
Let us first prove (\ref{sv}). We replace $u= R_1+R_2+\eta$ in the definition of $\mathcal G$. We obtain
\bee
\mathcal{M}(t) & = &\frac 12 \int_\R ( R_1 +R_2 +\eta)^2(1-\psi)\\  
& = & \frac 12 \int_\R R_2^2 (1-\psi) + \int_\R \eta R_2 (1-\psi) +  \frac 12 \int_\R \eta^2 (1-\psi) +O(c^{10}).
\eee
Here we have used the estimate for $t\geq T_c$
$$
\abs{\int_\R R_1(t) (1-\psi)} \leq K e^{-\frac 12 t} \leq Kc^{10},
$$
among other similar estimates.

In the same way,
$$
\frac 12 \int_\R u_x^2(1-\psi)  = \frac 12 \int_\R (R_2)_x^2 (1-\psi) + \int_\R \eta_x (R_2)_x (1-\psi)  +  \frac 12 \int_\R \eta_x^2 (1-\psi) +  O(c^{10}).
$$
Finally, using the character exponentially decreasing of $R_1$ where $1-\psi$ is away from zero,
\bee
\int_\R F(u)(1-\psi) & =  & \int_\R F(R_1+R_2+\eta)(1-\psi )\\
 &=&  \int_\R [F(R_1+R_2+\eta) - F(R_2 +\eta)  ] (1-\psi )\\
% & & +  \int [ F(R_1+R_2)  -F(R_1)-f(R_1)R_2 ](1-\psi )\\ 
 & & +  \int_\R [ F (R_2 + \eta )  - F(R_2)  -f(R_2)\eta - \frac 12 f'(R_2)\eta^2 ](1-\psi )\\ 
& & +  \int_\R [ F(R_2) +  f(R_2)\eta +  \frac 12f'(R_2)\eta^2  ](1-\psi )\\ 
%& & +  \int [ F(R_1) +  f(R_2+\eta)R_1 +f(R_1)R_2   ](1-\psi )\\ 
&=&  \int_\R [ f(R_2+\eta)R_1 + O(R_1^2)](1-\psi )
% & & +  \int [ F(R_1+R_2)  -F(R_1)-f(R_1)R_2 ](1-\psi )\\ 
 +  O(\norm{R_2}_{L^{\infty}(\R)}^{m-2}\int_\R \abs{\eta(t)}^3 (1-\psi))  \\ 
& & +  \int_\R [ F(R_2) +  f(R_2)\eta +  \frac 12f'(R_2)\eta^2  ](1-\psi )\\ 
& = &  \int_\R [ F(R_2) +  f(R_2)\eta +  \frac 12f'(R_2)\eta^2  ](1-\psi ) \\
& & + O\big[\norm{R_2}_{L^{\infty}(\R)}^{m-2}\norm{\eta(t)}_{H^1}\int_\R \eta^2(t) (1-\psi)\big] + O(c^{10}).
\eee

From this, 
$$
\mathcal{G}(t) = \mathcal{G}[R_2](t) + \mathcal{H}(t) + O\big[\norm{R_2}_{L^{\infty}(\R)}^{m-2}\norm{\eta(t)}_{H^1}\int_\R \eta^2(t) (1-\psi)\big]+ O(c^{10}).
$$
Putting together these estimates, using the value of $a$, evaluating at times $t=T_c$ and $t=T_0$ and using (\ref{eq:yy}) and (\ref{eq:tt}), we obtain the desired result.

\smallskip

The proof of (\ref{coe}) is standard, see e.g. \cite{MMas2} Appendix B.3.

\end{proof}

\medskip

Combining Lemmas \ref{le:mo} and \ref{le:qd}, we find
\bea
  \int_\R \left[c \eta^2 + \eta_x^2 \right](T_{c}) (1-\psi) & \leq &  K \mathcal{H}(T_c) \nonumber \\
& \leq &  K \mathcal{H}(T_0) + K \abs{d(\ve)} c^{\frac 1{m-1}+\frac 34} \int_\R \eta^2(T_c) (1-\psi) \nonumber \\
 & &  \quad +\ K\alpha^2 \abs{d(\ve)}^3 c^{\frac 5{m-1} + \frac 54}  +O(c^{10}) +K(\mathcal{G}(T_{c})-\mathcal{G}(T_0)) \nonumber  \\
 & \leq &  K\alpha^2 \abs{d(\ve)}^2 c^{\frac 4{m-1} + \frac 32} + K\alpha^2  \abs{d(\ve)}^3 c^{\frac 5{m-1}+\frac 54} + Kc^{10}. \label{merc}
% & & + K\abs{d(\ve)}c^{\frac 2{m-1}-\frac 14} \int_\R \eta^2(T_c) (1-\psi)
\eea %and \eqref{eq:xx}.
The last inequality is consequence of 
$$
\abs{d(\ve)}c^{\frac 1{m-1} + \frac 34} \ll \frac cK,  \quad m=2,3 \, \hbox{ and }\,  4;
$$
therefore the term  
$$
 K \abs{d(\ve)} c^{\frac 1{m-1}+\frac 34} \int_\R \eta^2(T_c) (1-\psi)
$$
can be sent to the left hand side of (\ref{merc}).
Using \eqref{eq:aT} we finally get
$$
\abs{d(\ve)}^2 c^{\frac 4{m-1} + \frac 32} \leq K \alpha^2 \abs{d(\ve)}^2 c^{\frac 4{m-1} + \frac 32} +  K\alpha^2  \abs{d(\ve)}^3 c^{\frac 5{m-1}+\frac 54}.
$$
But this estimate is a contradiction for $\alpha>0$ small enough and $0<c<c_0$ small enough (it is enough to put $\ve$ even smaller). The proof of Claim \ref{5.1} is now complete.
\end{proof}

\medskip

{\bf 5. Lower bounds on the parameters}. 
We finally prove (\ref{scagen}). This result is a consequence of Theorem \ref{PURE}, (\ref{th02bis}), (\ref{th03}) and (\ref{th04}), see also \cite{MMcol2} for the proof. Indeed, from (\ref{th02bis}) and (\ref{restgen}), we have
$$
\frac 1K \abs{d(\ve)}^2 c^{\frac 32 +\frac 4{m-1}} \leq 2E^+ + cM^+ \leq \frac 1K \abs{d(\ve)}^2 c^{\frac 12 +\frac 4{m-1}}.
$$
The final conclusion follows from (\ref{th03}), (\ref{th04}) and (\ref{diff0})-(\ref{diff2}).

This finishes the proof of the Theorem \ref{PUREgen}.
\end{proof}

%%%%%%%%%%%%%%%%%%%%%%%%%%%
%%%%%%%%%%%%%%%%%%%%%%%%%%%
%%%%%%%%%%%%%%%%%%%%%%%%%%%

\appendix

\section{Proof of Proposition \ref{prop:decomp}}\label{sec:B}

The proof is similar to Proposition 2.2 in \cite{MMcol2} and Appendix in \cite{MMcol1}. The main difference consists in the fact that we need to know explicitly all linear systems up to order $m+1$ to show the nonexistence of growing solutions. We will discard several trivial terms by using the property {\bf (IP)}.
For this purpose it is better to state an improved version of Proposition \ref{prop:decomp}. Before that we introduce a useful notation.

\begin{defn}
Consider $f,g:\R\to \R$ given functions. We say that $f=g \mod \mathcal{Y}$ if there exists 
$h\in \mathcal{Y}$ such that $f=g+h$. 
\end{defn}

%Of course, this definition mimics the classical equivalence relation mod. 

In our case, this definition will be useful to discard localized functions in the source terms. Indeed, 

\begin{prop}[Decomposition of $S(\tilde u)$, improved version]\label{prop:decompbis}

Assume that $f$ is of class $C^{m+2}$.  Let
\begin{equation}\label{defLy1}
    \mathcal{L} w = - w_{yy} + w -f'(Q) w.
\end{equation}
 Then,
    \begin{align}\label{S}
            S[\tilde u](t,x) & = 
        \sum_{(k,l)\in \Sigma_m}
        c^l Q_c^k(y_c)    \Big[a_{k,l} (-3 Q+2 f(Q))'(y)    -(\mathcal{L} A_{k,l})'(y) +F_{k,l}(y)  \Big] \nonumber
        \\& \quad
        + \sum_{(k,l)\in \Sigma_m}
        c^l (Q_c^k)'(y_c)    
        \Big[a_{k,l} (-3 Q'')(y) + \left(3A_{k,l}'' +f'(Q) A_{k,l}\right)(y)    - (\mathcal{L} B_{k,l})'(y) + G_{k,l}(y) \Big] \nonumber
        \\& \quad +  \mathcal{E}(t,x)
    \end{align}
    where 
    $F_{k,l} $, $G_{k,l} $ and $\mathcal{E}$ satisfy, for any $(k,l)\in \Sigma_m$,
    \begin{itemize}
        \item[{\rm (i)}] Dependence property of $F_{k,l}$ and $G_{k,l}$:
        The expressions of $F_{k,l}$ and $G_{k,l}$ depend only on $(a_{k',l'})$, $(A_{k',l'})$, $(B_{k',l'})$ for $(k',l') < (k,l)$.
        \item[{\rm (ii)}] Parity property of $F_{k,l}$ and $G_{k,l}$:  Assume that for any $(k',l')$ such that $(k',l')< (k,l)$
         $A_{k',l'}$ is even and    $B_{k',l'}$ is odd, then 
        $F_{k,l} $ is odd and  $G_{k,l}$ is  even.

       \item[{\rm (iii)}] Explicit source terms: We have $F_{1,0}=(f'(Q))'$ and $G_{1,0} = f'(Q)$, 
       \bee
F_{2,0} & = &  -(3A_{1,0}'+3B_{1,0}'' +f'(Q)B_{1,0}) +\frac 12 (f''(Q)(2A_{1,0} +A_{1,0}^2))' \\
& & \quad  - a_{1,0} ( 3A_{1,0}''  -Q +f'(Q)(1+A_{1,0})  )' + 3a_{1,0}^2 Q^{(3)} + \frac 12 (f''(Q)-2)',
\eee
and
\bee
G_{2,0} & = &  \frac 12 (f''(Q)-2) -(A_{1,0}+3B_{1,0}')  + \frac 12 f''(Q)(2A_{1,0}+A_{1,0}^2) + \frac 32 a_{1,0}^2 Q '' \\
& &    -\frac 12 a_{1,0}(9A_{1,0}' +3B_{1,0}'' +f'(Q)B_{1,0})'  +\frac 12 (f''(Q)(B_{1,0} +A_{1,0}B_{1,0}))' .
\eee
for the case $m=2$, and
\bee
F_{2,0}  =  (\frac 12 f''(Q)(1+A_{1,0})^2)' +3a_{1,0}^2 Q^{(3)} -a_{1,0} (f'(Q) +3A_{1,0}'' +f'(Q)A_{1,0})'  
\eee
and
\bee
G_{2,0} & = & \frac 12 f''(Q)(1+A_{1,0})^2 +\frac 32 a_{1,0}^2 Q'' -\frac 12 a_{1,0}(9A_{1,0}' +3B_{1,0}'' +f'(Q)B_{1,0})' \\
& & + \frac 12 (f''(Q)(1+A_{1,0})B_{1,0})'.
\eee
in the case $m=3$. If property {\bf (IP)} holds for $(k,l)=(1,0)$, then each term above is in $\mathcal{Y}$.

       \item[{\rm (iv)}] Explicit high order source terms modulo $\mathcal{Y}$: Suppose property {\bf (IP)} holds for $(k,l)\in \Sigma_m$ with $k+l\leq 2$. Then, for the quadratic case,
       $$
       F_{1,2}, G_{1,2}, F_{2,1} \hbox{ and } F_{3,0}  \in \mathcal{Y}; \quad G_{3,0}=-\frac 23(B_{1,0}^2 +2A_{2,0}) \mod \mathcal{Y},
       $$
       and
       $$
       G_{2,1} = B_{1,0}^2 +A_{1,1}+ 3A_{2,0} \mod \mathcal{Y}.
       $$
        For the cubic case,
        $$
        F_{3,0}, G_{3,0}, F_{2,1}, F_{4,0} \hbox{ and } G_{4,0}  \in \mathcal{Y}, \quad G_{2,1} = 3A_{2,0} \mod \mathcal{Y}.
        $$
       
       \item[{\rm (v)}] Improved estimate on $\mathcal{E}$: Suppose in addition that property {\bf (IP)} holds for any $(k,l)\in \Sigma_m$, then for all $j=0,1,2$ 
       $$
       \norm{\partial_x^j \mathcal{E}(t,x)}_{H^1(\R)} \leq K c^{\frac 34 +\frac 3{m-1}}.
       $$
    \end{itemize}
\end{prop}

\begin{proof}
Expansion  (\ref{S}), and items (i) and (ii) were proven in \cite{MMcol2}, so in what follows we deal with (iii)-(v). For this it is necesary to improve the computation done in \cite{MMcol2}.
  
We start with an important lemma concerning the algebra of $Q_c$. 

\begin{lem}[Properties of $Q_c$, see Lemma 2.1 in \cite{MMcol2}]\label{surQc2}
Suppose $0<c\leq 1$, $0<\ve\leq \ve_0$ small, $ k\in \{1,\ldots,k_0\}$, and $m=2,3$. Then

\begin{enumerate}
\item There exists a positive constant $K=K(\ve)>0$ such that 
\be\label{decay}
\frac 1K c^{\frac 1{m-1}} e^{-\sqrt{c}|x|} \leq Q_c(x)\leq K c^{\frac 1{m-1}} e^{-\sqrt{c}|x|},\quad
|Q_c'(x)|\leq K c^{\frac 1{m-1}+\frac 12} e^{-\sqrt{c}|x|}.
\ee
\item For $F$ defined in (\ref{F}) and any $k\geq 1$,
\be\label{taylor0}
Q_c'' = cQ_c -f(Q_c), \quad Q_c'^2 = cQ_c^2 - 2F(Q_c).
\ee
\be\label{taylor1}
(Q_c^k)'' = c k^2  Q_c^k -2k(k-1)Q_c^{k-2}F(Q_c) -kf(Q_c)Q_c^{k-1}.
\ee
\end{enumerate}
\end{lem}

\medskip

We recall the notation introduced in Subsection \eqref{sec:2-1}: 
\begin{equation*}
  S[\tilde u] =\tilde u_t +(\tilde u_{xx} -\tilde u +f(\tilde u))_x.
\end{equation*}
We easily verify that 
\begin{equation}\label{eq:sion}
S[\tilde u]= \bf I +II +III + IV,
\end{equation}
where (we omit the dependence on $t,x$) 
$$
{\bf I} := S[R], \quad {\bf II} := (f(R+R_c) -f(R)-f(R_c))_x,
$$
and
\begin{equation}\label{eq:sion2}\begin{split}
& \mathcal{L}=-\partial_x^2 + 1- f'(Q),\\
& {\bf III} = {\bf III}(W) := W_t - (\mathcal{L}W)_x , \\
&{\bf IV} : = \left\{ f(R+R_c +W) - f(R +R_c) -f(R)W\right\}_x.
\end{split}\end{equation}
Since $\wqs(\ys)$ is a solution to (\ref{ellipticf}),
we have $S(\wqs)=0.$

\begin{Cl}
  \label{cl:SKdV1}
Let $A=A(y)$ and $q=q(y_c)$ be $C^3$-functions with $y,y_c$ defined in Section \ref{sec:2-1}. Then 
\bee
{\bf III}(Aq)  & = & -q(\mathcal L A)' + q'(3A''+f'(Q)A)  \\
& & + q(-3\beta A^{(3)} -\beta A' -3\beta_x A'' -A'\beta_{xx} +\beta A' -\beta (f'(Q)A)' ) \\
& & + q( 3\beta^2 A^{(3)} +3\beta\beta_x A'' -\beta^3 A^{(3)} + c\beta A') \\
& &  + q'( -cA-6A''\beta -3A'\beta_x + 3A''\beta^2) \\
& & + q''( 3A' - 3A'\beta) + Aq^{(3)}.
\eee
\end{Cl}
\begin{proof}
Direct differentiation, see \cite{MMcol2}, Proposition 2.2.
\end{proof}

\begin{Cl}
\label{cl:beta}
Recall from (\ref{defALPHA}),
\be\label{beta}
\beta  =  \sum_{(k,l)\in \Sigma_m} a_{k,l}  c^l \wqs ^k  (y_c) 
\ee
Then, for some fixed numbers $\hat a_{k,l}^1, \hat a_{k,l}^2, \bar a_{k,l}, \tilde a_{k,l}$ with $(k,l)\in \Sigma_m$, depending only on $a_{k',l'}$ with $(k',l')\leq (k,l)$, we have 
$$\begin{cases}
\beta_x  =   \sum_{(k,l)\in \Sigma_m} a_{k,l}  c^l (\wqs ^k )' (y_c), \\
\beta_{xx}  =  \sum_{\substack{(k,l)\in \Sigma_m\\ l \geq 1} } \hat a^1_{k,l}  c^l \wqs ^k (y_c) + \sum_{\substack{(k,l)\in \Sigma_m\\ k\geq m}} \hat a^2_{k,l}  c^l \wqs ^k  (y_c) + O(Q_c^5 + cQ_c^3 ),\\
\beta^2  =  \sum_{\substack{(k,l)\in \Sigma_m\\ k\geq 2} } \bar a_{k,l}  c^l \wqs ^k  (y_c) + O(Q_c^5 + cQ_c^3 )\\
(\beta^2)_x  =  \sum_{\substack{(k,l)\in \Sigma_m\\k\geq 2}} \bar a_{k,l}  c^l (\wqs ^k )' (y_c) + O(Q_c^5 + cQ_c^3 ), \hbox{ and }\\
\beta^3   =  \sum_{\substack{(k,l)\in \Sigma_m\\ k\geq 3} } \tilde a_{k,l}  c^l \wqs ^k  (y_c) + O(Q_c^5 + cQ_c^3 ) .
\end{cases}$$
\end{Cl}

\begin{proof}
The proof follows by elementary calculations from \eqref{defALPHA}.
\end{proof}

In the next lemmas, we expand the terms in \eqref{eq:sion}.

\begin{lem}
\label{lem:SQ}
\begin{equation}
  \label{eq:SQ}\begin{split}
{\bf I}= &\sum_{(k,l)\in \Sigma_m} c^l\left[ \wqs^k(\ys)a_{k,l}( 2f(Q) -3Q )'(y)
+  (\wqs^k)'(\ys) (-3a_{k,l}Q''(y)) \right]
\\ &+\sum_{(k,l)\in \Sigma_m}c^l\left(\wqs^k(\ys)F_{k,l}^I(y)
+ (\wqs^k)'(\ys)G_{k,l}^I(y)\right)+c^3 O(\wqs(\ys)),
\end{split}
\end{equation}
where
\begin{align*}
& F_{1,0}^I = G_{1,0}^I = F_{1,1}^I= G_{1,1}^I = 0,\\
& F_{2,0}^I= 3 a_{1,0}^2Q^{(3)} + a_{1,0}Q' \delta_{m2} , \quad G_{2,0}^I= \frac{3}{2}a_{1,0}^2 Q'',
\end{align*}
and for all $(k,l)\in \Sigma_m$, $F_{k,l}^I\in \mathcal{Y}$ is odd, $G_{k,l}^I\in \mathcal{Y}$ is even and depend only on $a_{k',l'}$ for $(k',l')<(k,l).$
\end{lem}

\begin{proof}[Proof of Lemma \ref{lem:SQ}.]
We have (here $'$ denotes derivative with respect to $y$)
\bee
{\bf I} & = & R_t + (R_{xx} -R +f(R))_x\\
&= & -(1-c)\beta Q'  + (f(Q))'(1-\beta)  -Q'(1-\beta) +( Q''(1-\beta)^2 -Q'\beta_x )_x \\
& = & (Q''-Q +f(Q))' + Q^{(3)}(-3\beta +3\beta ^2 -\beta^3) -3Q''(\beta_x-\beta \beta_x ) - \beta_{xx}Q' -\beta (f(Q))' +c\beta Q' \\
&  = & -[ 3\beta Q^{(3)} + 3Q'' \beta_x + \beta (f(Q))' ] + 3\beta^2 Q^{(3)} + 3\beta\beta_x Q'' -\beta_{xx}Q' + c\beta Q' -\beta^3 Q^{(3)}. 
\eee
Hence using Claim \ref{cl:beta}, we obtain
\begin{align*}
{\bf I}&= a_{1,0}(2f(Q)-3Q)'  \wqs(\ys) + a_{1,0}(-3Q'')\wqs'(\ys) \\
 & + \big(a_{2,0}(2f(Q)-3Q)' + 3a_{1,0}^2Q^{(3)} +a_{1,0}Q'\delta_{m2} \big) \wqs^2(\ys)\\ 
&+ \big( a_{2,0}(-3Q'')+\tfrac{3}{2}a_{1,0}^2Q''\big) (\wqs^2)'(\ys) \\
& + \sum_{k+l=3,4} c^l\big( a_{k,l} (2f(Q) -3Q)'(y) \wqs^k(\ys) +  a_{k,l}(-3Q'')(y) (\wqs^k)'(\ys) \big) \\
& +\sum_{k+l=3,4} c^l\big( F^I_{k,l} \wqs^k(\ys) + G^I_{k,l} (\wqs^k)'(\ys)\big) +c^3 O(\wqs),
\end{align*}
where for all $k+l=3$, $F_{k,l}^I\in \mathcal{Y}$ and $G_{k,l}^I\in \mathcal{Y}$, as claimed in the statement of the Lemma.
\end{proof}

\begin{lem}
  \label{lem:SintII}
  \begin{align*}
{\bf II} =\sum_{(k,l)\in \Sigma_m}c^l
\left(\wqs^k(\ys)F_{k,l}^{II}(y)+(\wqs^k)'(\ys)G_{k,l}^{II}(y)\right)+  O(\wqs^{m+2}),
  \end{align*}
where for all $(k,l)\in \Sigma_m$ and for all $p\geq m+1$, $F_{k,l}^{II}, G_{k,l}^{II} \in \mathcal{Y}$ and are odd and even respectively. Moreover, for $m=2$,
\begin{align*}
&  F_{1,0}^{II}= (f'(Q))',\quad G_{1,0}^{II}=f'(Q),\quad
F_{1,1}^{II}=G_{1,1}^{II}=0,\\
& F_{2,0}^{II} =( \frac 12 f''(Q) -a_{1,0}f'(Q) )' ,\quad  G_{2,0}^{II}= \frac 12 f''(Q)-1.
\end{align*}
Finally, if $m=3$,
\begin{align*}
&  F_{1,0}^{II}= (f'(Q))',\quad G_{1,0}^{II}=f'(Q),\quad
F_{1,1}^{II}=G_{1,1}^{II}=0,\\
& F_{2,0}^{II} = (\frac 12 f''(Q) - a_{1,0}f'(Q))', \quad G_{2,0}^{II} = \frac 12 f''(Q).
\end{align*}

\end{lem}

\begin{proof}
First define ${\bf \tilde{II} } := f(R+R_c) -f(R) - f(R_c)$. Note that
$$
{\bf \tilde{II} }  =   f'(R)R_c  +  \frac 12 f''(R)R_c^2  + \frac 16 f^{(3)}(R)R_c^3 - f(R_c)  +  \frac 1{24} f^{(4)}(R)R_c^4+ O(R_c^5),
$$
Thus taking derivative
\bea
{\bf II } & = &  (f'(Q))'(1-\beta)Q_c + f'(Q)Q_c' +  \frac 12 (f''(Q))'(1-\beta)Q_c^2 + \frac 12 f''(Q)(Q_c^2)'   + \frac 16 (f^{(3)}(Q))' (1-\beta )Q_c^3\nonumber \\ 
& &  + \frac 16 f^{(3)}(Q) (Q_c^3)'  +\frac 1{24}(f^{(4)}(Q))' Q_c^4 +\frac 1{24}f^{(4)}(Q)(Q_c^4)'  -(f(Q_c))' + O(Q_c^5).\label{II3}
\eea
%Note that from $Q_c^5 = c^{\frac{4}{m-1}}O(Q_c)$. 
Here we have to identify two different results, depending on the value of $m$. For $m=2$, namely, the quadratic case, we will need only up to third order terms. After replacing the value of $\beta$ given by (\ref{beta}), we will obtain (recall that $p\geq 3$)
\bee
{\bf II } & = &  (f'(Q))' Q_c + f'(Q)Q_c' + ( \frac 12 f''(Q) -a_{1,0}f'(Q) )' Q_c^2 + (\frac 12 f''(Q)-1)(Q_c^2)'   -a_{1,1} (f'(Q))' cQ_c^2 \\
& & + ( \frac 16 f^{(3)}(Q) -\frac 12 a_{1,0}f''(Q) -a_{2,0} f'(Q) )' Q_c^3  + \frac 16 f^{(3)}(Q) (Q_c^3)'-\ve (Q_c^p)' - (f_1(Q_c))'  +  O(Q_c^4).
\eee
It is easy to check that every term depending on $y$ up to order $Q_c^3, (Q_c^3)'$  is indeed in the class $\mathcal{Y}$. Even in the worst case, $p=3$, we will have the cancelation
$$
\frac 16 f^{(3)}(Q) (Q_c^3)'-\ve (Q_c^p)'  = \frac 16 f_1^{(3)}(Q) (Q_c^3)',
$$
 with $ \frac 16 f_1^{(3)}(Q)\in \mathcal{Y}$.
 
 Let us consider now the cubic case, $m=3$. The procedure is completely similar, although we must keep the fourth order terms. We start by replacing $\beta $ in (\ref{II3}) and collecting similar terms
 \bee
{\bf II } & = &  (f'(Q))'(1-\beta)Q_c + f'(Q)Q_c' +  \frac 12 (f''(Q))'(1-\beta)Q_c^2 + \frac 12 f''(Q)(Q_c^2)'   + \frac 16 (f^{(3)}(Q))' (1-\beta )Q_c^3\\ 
& &  + \frac 16 f^{(3)}(Q) (Q_c^3)'  +\frac 1{24}(f^{(4)}(Q))' Q_c^4 +\frac 1{24}f^{(4)}(Q)(Q_c^4)'  -(f(Q_c))' + O(Q_c^5).\\
& = & (f'(Q))'Q_c + f'(Q)Q_c' + (\frac 12 f''(Q) - a_{1,0}f'(Q))' Q_c^2+ \frac 12 f''(Q)(Q_c^2)'   \\
& &  ( \frac 16 f^{(3)}(Q) - \frac 12 a_{1,0}f''(Q) - a_{2,0}f'(Q) )' Q_c^3   + (\frac 16 f^{(3)}(Q) -1) (Q_c^3)'  -a_{1,1}(f'(Q))' cQ_c^2\\
& &  +(\frac 1{24} f^{(4)}(Q) -\frac 16 a_{1,0}f^{(3)}(Q) -\frac 12 a_{2,0}f''(Q) -a_{3,0}f'(Q) )' Q_c^4 +\frac 1{24}f^{(4)}(Q)(Q_c^4)' \\
& &  -( \ve Q_c^p +f_1(Q_c))' + O(cQ_c^3 + cQ_c^4 + Q_c^5).
\eee
It is straightforward to check that every function depending on $y$ is indeed in $\mathcal{Y}$. The only complicated terms are (note that $p\geq 4$)
$$
\frac 16 f^{(3)}(Q) -1 = \frac 16 p(p-1)(p-2)\ve Q^{p-3}  + \frac 16 f_1^{(3)}(Q) \in \mathcal{Y},
$$
which is in front of $(Q_c^3)'$; and for $p=4$, facing $(Q_c^4)'$ we have
$$
\frac 1{24} f^{(4)}(Q) -\ve = \ve + \frac 1{24}f_1^{(4)}(Q) -\ve =  \frac 1{24}f_1^{(4)}(Q) \in \mathcal{Y}. 
$$
\end{proof}

\begin{lem}
  \label{lem:dSKdVw}
  \begin{align*}
{\bf III} &= \sum_{(k,l)\in\Sigma_m } c^l
\left( \wqs^k(\ys)(-\mathcal{L} A_{k,l})'(y)+ (\wqs^k)'(\ys) ((-\mathcal{L} B_{k,l})'+3A_{k,l}''+ f'(Q) A_{k,l})(y)\right)
\\&
+\sum_{(k,l)\in\Sigma_m}c^l
\left(\wqs^k(y_c) F_{k,l}^{III}(y)+(\wqs^k)' (y_c)G_{k,l}^{III}(y)\right)+ O(\wqs^{m+2}),
  \end{align*}
where
\begin{align*}
&  F_{1,0}^{II}=0,\quad
G_{1,0}^{II}=0,\quad 
 F_{1,1}^{II}= 3A_{1,0}' + 3B_{1,0}'' +  f'(Q)B_{1,0},\quad
 G_{1,1}^{II}= 3B_{1,0}',\\
& F_{2,0}^{II}= -a_{1,0}(3A_{1,0}''+ f'(Q) A_{1,0})' -(3A_{1,0}' + 3B_{1,0}'' + f'(Q) B_{1,0})\delta_{m2}\\
& G_{2,0}^{II}=-\frac 12 a_{1,0} \left(9A_{1,0}'+3 B_{1,0}''+ f'(Q) B_{1,0}\right)' - (A_{1,0}+3B_{1,0}')\delta_{m2},
\end{align*}
and for $(k,l)\in \Sigma_m$, $F_{k,l}^{II}$, $G_{k,l}^{II}$ depend on 
$A_{k',l'}$, $B_{k',l'}$ such that $(k',l')< (k,l)$.
Moreover, if $A_{k',l'}$ are even and $B_{k',l'}$ are odd then 
$F_{k,l}^{II}$ are odd and $G_{k,l}^{II}$ are even.

Finally, the following important property holds. Suppose {\bf (IP)} holds for any $(k,l)\in \Sigma_m$ with $k+l\leq 2$.  Then we have a sharp decomposition for each high order source term: 

\begin{enumerate}
\item For $m=2$,
\be\label{decayIII2}
F_{3,0}^{III} = 0 , \quad G_{3,0}^{III} =- \frac{10}3A_{2,0} , \quad F_{2,1}^{III} = 0,\quad G_{2,1}^{III}=-A_{1,1} + 3A_{2,0}, \quad F_{1,2}^{III} = G_{1,2}^{III}= 0 \mod \mathcal{Y}.
\ee
\item For $m=3$,
\be\label{decayIII3}
F_{3,0}^{III}, G_{3,0}^{III}, F_{2,1}^{III} \in \mathcal Y ,\quad G_{2,1}^{III}= 3A_{2,0}, \quad F_{4,0}^{III} =0 , \quad G_{4,0}^{III}= -3A_{2,0} \mod \mathcal{Y}.
\ee
\end{enumerate}

\end{lem}

\begin{proof}
We have, thanks to the linearity of the operator ${\bf III}(\cdot)$,
\begin{equation*}
{\bf III}(W) =\sum_{(k,l)\in \Sigma_m}c^l
\left(  {\bf III}(A_{k,l}(y)\wqs^k(\ys)) + {\bf III}(B_{k,l}(y)(\wqs^k)'(\ys))\right).
\end{equation*}
In what follows, for commodity of notation we omit the variables $y,y_c$, if there is no related confusion.  First, we compute ${\bf III}(A_{1,0}(y)\wqs(\ys))$.
By Claim \ref{cl:SKdV1} and the definition of $\beta$, we have
\bee
{\bf III} (A_{1,0}\wqs ) & = & -Q_c(\mathcal L A_{1,0})' + Q_c'(3A_{1,0}''+f'(Q)A_{1,0})  \\
& & + Q_c(-3\beta A_{1,0}^{(3)} -\beta A_{1,0}' -3\beta_x A_{1,0}'' -A_{1,0}'\beta_{xx} +\beta A_{1,0}' -\beta (f'(Q)A_{1,0})' ) \\
& & + Q_c( 3\beta^2 A_{1,0}^{(3)} +3\beta\beta_x A_{1,0}'' -\beta^3 A_{1,0}^{(3)} + c\beta A_{1,0}') \\
& &  + Q_c'( -cA_{1,0}-6A_{1,0}''\beta -3A_{1,0}'\beta_x + 3A_{1,0}''\beta^2) \\
& & + Q_c''( 3A_{1,0}' - 3A_{1,0}'\beta) + A_{1,0}Q_c^{(3)}\\
& = & -(\mathcal L A_{1,0})' Q_c + (3A_{1,0}''+f'(Q)A_{1,0})Q_c' + 3A_{1,0}'' cQ_c \\
& & - ( 3a_{1,0} A_{1,0}'' + a_{1,0}f'(Q)A_{1,0} + 3A_{1,0}\delta_{m2})' Q_c^2 -  (\frac 92 a_{1,0} A_{1,0}'' + A_{1,0}\delta_{m2} )(Q_c^2)' \\
& & + \sum_{3\leq k+l \leq 4} c^l ( F_{k,l} Q_c^k + G_{k,l} (Q_c^k)') + O(cQ_c^3 + Q_c^5 + c^2Q_c).
\eee
Moreover, by hypothesis $A_{1,0}\in \mathcal{Y}$ so we have all the source terms $F_{k,l}, G_{k,l}\in \mathcal{Y}$, as can be verified directly.

Now, we compute ${\bf III}(B_{1,0}(y)\wqs'(\ys))$ in a similar way:
\bee
{\bf III} (B_{1,0}\wqs' ) & = & -Q_c'(\mathcal L B_{1,0})' + Q_c''(3B_{1,0}''+f'(Q)B_{1,0})  \\
& & + Q_c'(-3\beta B_{1,0}^{(3)} -\beta B_{1,0}' -3\beta_x B_{1,0}'' -B_{1,0}'\beta_{xx} +\beta B_{1,0}' -\beta (f'(Q)B_{1,0})' ) \\
& & + Q_c'( 3\beta^2 B_{1,0}^{(3)} +3\beta\beta_x B_{1,0}'' -\beta^3 B_{1,0}^{(3)} + c\beta B_{1,0}') \\
& &  + Q_c''( -cB_{1,0}-6B_{1,0}''\beta -3B_{1,0}'\beta_x + 3B_{1,0}''\beta^2) \\
& & + Q_c^{(3)}( 3B_{1,0}' - 3B_{1,0}'\beta) + B_{1,0}Q_c^{(4)}\\
& = & -(\mathcal L B_{1,0})' Q_c' - Q_c^2 ( 3B_{1,0}'' + f'(Q)B_{1,0}) \delta_{m2} +  (3B_{1,0}''+f'(Q)B_{1,0})cQ_c  \\
& & - ( 3a_{1,0} A_{1,0}'' + a_{1,0}f'(Q)A_{1,0} + 3A_{1,0}\delta_{m2})' Q_c^2 -  (\frac 92 a_{1,0} A_{1,0}'' + A_{1,0}\delta_{m2} )(Q_c^2)' \\
& & + \sum_{3\leq k+l \leq 4} c^l ( F_{k,l} Q_c^k + G_{k,l} (Q_c^k)') + O(cQ_c^3 + Q_c^5 + c^2Q_c).
\eee

Suppose now that $2\leq k+l \leq 4$. Here we will use {\bf (IP)} for $k+l\leq 2$ to discard several terms of a tedious but direct computation. Indeed, from Claim \ref{cl:SKdV1} we have
\bee
{\bf III}(A_{k,l}  \wqs^k )  & = & -Q_c^k(\mathcal L A_{k,l})' + (Q_c^k)'(3A_{k,l}''+f'(Q)A_{k,l})  \\
& & + Q_c^k(-3\beta A_{k,l}^{(3)} -\beta A_{k,l}' -3\beta_x A_{k,l}'' -A_{k,l}'\beta_{xx} +\beta A_{k,l}' -\beta (f'(Q)A_{k,l})' ) \\
& & + Q_c^k( 3\beta^2 A_{k,l}^{(3)} +3\beta\beta_x A_{k,l}'' -\beta^3 A_{k,l}^{(3)} + c\beta A_{k,l}') \\
& &  + (Q_c^k)'( -cA_{k,l}-6A_{k,l}''\beta -3A_{k,l}'\beta_x + 3A_{k,l}''\beta^2) \\
& & + (Q_c^k)''( 3A_{k,l}' - 3A_{k,l}'\beta) + A_{k,l}(Q_c^k)^{(3)}\\
& = &  -Q_c^k(\mathcal L A_{k,l})' + (Q_c^k)'(3A_{k,l}''+f'(Q)A_{k,l})  + A_{k,l}(Q_c^k)^{(3)} -A_{k,l}c(Q_c^k)' \\
& & + \sum_{\substack{(k',l')\in \Sigma'_{m}\\ (k,l)\leq (k',l')}} c^{l'} ( F_{k',l'} Q_c^{k'} + G_{k',l'} (Q_c^{k'})') + O(cQ_c^3 + Q_c^5 + c^2Q_c).
\eee
Here $\Sigma'_{m} \subseteq \Sigma_{m}$ is the set of indices of \emph{third order} in $\Sigma_{m}$. More specificaly, 
\be\label{Sigm}
\Sigma'_{2} :=\{ (1,2), (2,1), (3,0) \}, \quad \Sigma'_{3} :=\{(2,1), (3,0), (4,0)\}.
\ee
The terms describing $ F_{k',l'}$ and $ G_{k',l'}$ with $(k',l')\in \Sigma'_{m}$ are in $\mathcal{Y}$ provided {\bf (IP)} is satisfied for every $(k,l)\in \Sigma_{m}\backslash \Sigma'_{m}$.

Now note that from (\ref{taylor1})
$$
(Q_c^k)^{(3)} = k^2 (cQ_c^k)' -\frac{k(2k+m-1)}{m+1} (Q_c^{k+m-1})' -\ve k \frac{k(2k+p-1)}{p+1} (Q_c^{k+p-1})' + O(Q_c^{k+p}).
$$

We can finally conclude that 
\bee
{\bf III}(A_{k,l}  \wqs^k ) & = &  -Q_c^k(\mathcal L A_{k,l})' + (Q_c^k)'(3A_{k,l}''+f'(Q)A_{k,l})   +  (k^2-1)A_{k,l}c(Q_c^k)'  \\ 
& & -\frac{k(2k+m-1)}{m+1} A_{k,l}(Q_c^{k+m-1})'  + \sum_{(k',l')\in \Sigma'_{m}} c^{l'} ( F_{k',l'} Q_c^{k'} + G_{k',l'} (Q_c^{k'})') \\
& & + O(cQ_c^3 + Q_c^5 + c^2Q_c),
\eee
where, as described above, the terms $ F_{k',l'}$ and $ G_{k',l'}$ with $(k',l')\in \Sigma'_{m}$ are in $\mathcal{Y}$ provided {\bf (IP)} is satisfied for every $(k,l)\in \Sigma_{m}\backslash \Sigma'_{m}$.

\smallskip

On the other hand, the terms of the form \be\label{B3}{\bf III}(B_{k,l}  (\wqs^k)' ), \quad(k,l)\in \Sigma_{m}, \quad  2 \leq k+l \leq 4, \ee  can be treated in the same way as above, and we only write the final result (see the computation of ${\bf III}(B_{1,0}Q_c' )$ for example):
\bee
{\bf III}(B_{k,l}  (\wqs^k)' ) & = &  -(Q_c^k)'(\mathcal L B_{k,l})'   + B_{k,l} (Q_c^k)^{(4)}- B_{k,l} (cQ_c^k)'' \\
& & + \sum_{1\leq k' \leq 4} c^{l'} ( F_{k',l'} Q_c^{k'} + G_{k',l'} (Q_c^{k'})') + O(cQ_c^3 + Q_c^5 + c^2Q_c)\\
&= & -(Q_c^k)'(\mathcal L B_{k,l})'  + \sum_{\substack{(k',l')\in \Sigma'_{m}\\ (k,l)\leq (k',l')}} c^{l'} ( F_{k',l'} Q_c^{k'} + G_{k',l'} (Q_c^{k'})') + O(cQ_c^3 + Q_c^5 + c^2Q_c).
\eee
%({\bf I think it is necessary to explain a little bit more})
 
 To obtain (\ref{decayIII2}) and (\ref{decayIII3}) we only evaluate the expressions for ${\bf III}(A_{k,l}  \wqs^k )$ and ${\bf III}(B_{k,l}  (\wqs^k)' )$ for each $(k,l)\in \Sigma_m$ with $2\leq k+l$. 
The final result follows from the sum of each term ${\bf III}(A_{k,l}  \wqs^k )$, ${\bf III}(B_{k,l}  (\wqs^k)' )$ for $(k,l)\in \Sigma_m$, discarding localized terms.
This concludes the proof.
\end{proof}

The final term reads

%%%%%%%%%%%%%%%%%%%%%%%%%%%%%%%%%%
%%%%%%%%%%%%%%%%%%%%%%%%%%%%%%%%%%

\begin{lem}
  \label{lem:Sint}
\be\label{IV}
{\bf IV} =\sum_{(k,l)\in \Sigma_m}c^l \left(\wqs^k(\ys)F_{k,l}^{IV}(y)+(\wqs^k)'(\ys)G_{k,l}^{IV}(y)\right)+c^3 O(\wqs),
\ee
where
\begin{align*}
&  F_{1,0}^{IV}=G_{1,0}^{IV}=0,\quad F_{1,1}^{IV}=G_{1,1}^{IV}=0,\\
& F_{2,0}^{IV}=\frac12  (f''(Q)  ( 2A_{1,0} + A_{1,0}^2) )',
\quad
G_{2,0}^{IV}=\frac12 \big[ f''(Q) ( 2A_{1,0} + A_{1,0}^2) + (f''(Q)(B_{1,0}+ A_{1,0}B_{1,0}) )'  \big] ,
\end{align*}
and for $(k,l)\in \Sigma'_m$ (see (\ref{Sigm})),   $F_{k,l}^{IV}$, $G_{k,l}^{IV}$ depend on 
$A_{k',l'}$, $B_{k',l'}$ for $(k',l')\in  \Sigma_m $ with $(k',l')< (k,l)$.
Moreover, if $A_{k',l'}$ are even and $B_{k',l'}$ are odd then 
$F_{k,l}^{IV}$ are odd and $G_{k,l}^{IV}$ are even.

Finally, suppose {\bf (IP)} holds for $(k,l)\in \Sigma_m$ with $k+l\leq 2$. Then the only non localized terms for $(k,l)\in \Sigma'_m$ are given by 
\be\label{145}
F_{2,1}^{IV}=0,\quad G_{2,1}^{IV} = \frac 12 f''(Q) ( 2A_{1,1} + B_{1,0}^2)\;  \mod \mathcal{Y},   
\ee
and
\be\label{146}
F_{3,0}^{IV} =  0,%f''(Q)(A_{2,0} -\frac 13 B_{1,0}^2),  
\quad G_{3,0}^{IV} =  f''(Q) ( A_{2,0} - \frac 13 B_{1,0}^2) \; \mod \mathcal{Y},
\ee
for the quadratic case, and
\be\label{147}
G_{4,0}^{IV}  = \frac 12 f^{(3)}(Q) A_{2,0} \; \mod \mathcal{Y}.
\ee
in the cubic case.
\end{lem}

\begin{proof}
As above, first define ${\bf \tilde{IV} } := f(R+R_c+W) -f(R+R_c) - f'(R)W$. Note that, using that $R :=Q(y)$ and $R_c :=Q_c(y_c)$,
\bee
{\bf \tilde{IV} } & = & (f'(Q+Q_c) -f'(Q))W + \frac 12 f''(Q+Q_c)W^2 + \frac 16 f^{(3)}(Q+Q_c)W^3 \\
&  &  + \frac 1{24}f^{(4)}(Q+Q_c)W^4  + O(W^5)\\ 
& = & [f''(Q)Q_c +\frac 12 f^{(3)}(Q)Q_c^2 +\frac 16f^{(4)}(Q)Q_c^3 + O(Q_c^4)  ]W \\
&  & + \frac 12[f''(Q) +f^{(3)}(Q) Q_c + \frac 12 f^{(4)}(Q) Q_c^2 + O(Q_c)^3 ] W^2   \\
&  &  + \frac 16[ f^{(3)}(Q) + f^{(4)}(Q)Q_c + O(Q_c^2) ] W^3 + \frac 1{24}f^{(4)} (Q)W^4 + O(Q_c^5)\\
& = &  f''(Q)(Q_cW +\frac 12 W^2) +\frac 12 f^{(3)}(Q)(Q_c^2W+W^2Q_c+\frac 13W^3 ) \\
& & + \frac 12 f^{(4)}(Q) (\frac 13 Q_c^3W +\frac 12 Q_c^2 W^2 +\frac 13 Q_c W^3 + \frac 1{12} W^4) + O(Q_c^5).
\eee

Now, the final value of ${\bf IV}$ depends on the different values of $m$. We will proceed carefully in both cases.

\smallskip

\noindent
\emph{Case $m=2$}. Here we consider only up to third order, namely
\bee
{\bf \tilde{IV} }&  =& \frac 12 f''(Q)( 2Q_cW +  W^2) +\frac 12 f^{(3)}(Q)(Q_c^2W+W^2Q_c+\frac 13W^3 ) + O(Q_c^4)\\
& = : &  {\bf \tilde{IV} }_2 + {\bf \tilde{IV} }_3 + O(Q_c^4).
\eee
First of all let us consider the third order term $ {\bf \tilde{IV} }_3$.  A quickly computation using (\ref{defW}) gives us
$$
Q_c^2 W = A_{1,0} Q_c^3 + \frac 13 B_{1,0}(Q_c^3)' + O(c^3Q_c),  
$$
and
$$
W^2Q_c = A_{1,0}^2Q_c^3 +\frac 23 A_{1,0} B_{1,0}(Q_c^3)' + O(c^3Q_c) , \quad W^3 = A_{1,0}^3 Q_c^3 + A_{1,0} B_{1,0} (Q_c^3)' + O(c^3 Q_c).
$$
If we suppose $A_{1,0}\in \mathcal{Y}$ and $B_{1,0}$ bounded  (this is actually the case), we will obtain
$$
 {\bf \tilde{IV} }_3 = F_{3,0}^{\tilde{IV}_3} Q_c^3 + (G_{3,0}^{\tilde{IV}_3} + \frac 16 f^{(3)}(Q) B_{1,0} ) (Q_c^3)' + O(c^3Q_c),
$$
where $F_{3,0}^{\tilde{IV}_3}, G_{3,0}^{\tilde{IV}_3}\in \mathcal{Y}$. Moreover, for $p\geq 4$, actually $f^{(3)}(Q) -6\mu(\ve ) \in \mathcal{Y}$, because of
$$
f^{(3)}(Q) = 6\mu(\ve)   + p(p-1)(p-2)\ve Q^{p-3} + f_1^{(3)}(Q).
$$
Now, let us compute in detail the term ${\bf \tilde{IV} }_2$. These terms above are important because they will give us source terms of second order. Now from the definition of $W$ in (\ref{defW}) it is easy to check that, up to third order,  
$$
Q_c W = A_{1,0}Q_c^2 + \frac 12 B_{1,0}(Q_c^2)' + A_{1,1}cQ_c^2 + \frac 12 B_{1,1}c(Q_c^2)' + A_{2,0}Q_c^3 + \frac 23 B_{2,0} (Q_c^3)' + O(c^3Q_c), 
$$
and
\bee
W^2 &  = & A_{1,0}^2 Q_c^2 +  A_{1,0}B_{1,0}(Q_c^2)'  +(  2A_{1,0}A_{1,1} +  B_{1,0}^2 )cQ_c^2 + (A_{1,0}B_{1,1} +B_{1,0}A_{1,1}) c(Q_c^2)' \\ 
& &    +  ( 2A_{1,0}A_{2,0} -\frac 23 B_{1,0}^2)Q_c^3 + \frac 23 (A_{1,0}B_{2,0} + B_{1,0}A_{2,0} )(Q_c^3)' .
\eee
From here, 
\bee
{\bf \tilde{IV}_2 } & = &  \frac 12 f''(Q) ( 2A_{1,0} +A_{1,0}^2)Q_c^2 +  \frac 12 f''(Q)( B_{1,0} + A_{1,0}B_{1,0}) (Q_c^2)'  \\
& & +  \frac 12 f''(Q) (F_{2,1}^{\tilde{IV}_2} + 2A_{1,1} + B_{1,0}^2) cQ_c^2 + \frac 12 f''(Q) (G_{2,1}^{\tilde{IV}_2} + B_{1,1} + A_{1,1}B_{1,0}) c(Q_c^2)' \\
& & +  \frac 12 f''(Q) (F_{3,0}^{\tilde{IV}_2} + 2A_{2,0} -\frac 23 B_{1,0}^2) Q_c^3 + \frac 12 f''(Q) (G_{3,0}^{\tilde{IV}_2} + \frac 43 B_{2,0} + \frac 23 A_{2,0}B_{1,0}) (Q_c^3)' \\
& & + O(c^3Q_c),
\eee 
where for $k+l=3$ it is satisfied $F_{k,l}^{\tilde{IV}_2}, G_{k,l}^{\tilde{IV}_2}\in \mathcal{Y}$, provided $A_{1,0}\in \mathcal{Y}$ and $B_{1,0}$ is bounded (namely $(k,l)=(1,0)$ satisfies {\bf (IP)}).

Putting all this information together, and after derivation, we obtain (note that $p\geq 3$ and $(f''(Q))'\in \mathcal{Y}$)
\bee
{\bf {IV}} &= & ({\bf \tilde{IV} }_2 + {\bf \tilde{IV} }_3 + O(Q_c^4))_x\\
& =&  \frac 12( f''(Q))' (1-\beta)( 2A_{1,0} +A_{1,0}^2)Q_c^2 +  \frac 12 f''(Q)( 2A_{1,0} +A_{1,0}^2)'(1-\beta)Q_c^2   \\
& &+ \frac 12 f''(Q)( 2A_{1,0} +A_{1,0}^2)(Q_c^2)' +  \frac 12 (f''(Q))'(1-\beta)( B_{1,0} + A_{1,0}B_{1,0}) (Q_c^2)'  \\
& & +  \frac 12 f''(Q)(1-\beta)( B_{1,0} + A_{1,0}B_{1,0})' (Q_c^2)'    \\
& & + F_{2,1}^{IV} cQ_c^2 + [ G_{2,1}^{IV} + \frac 12 f''(Q) ( 2A_{1,1} + B_{1,0}^2) ](cQ_c^2)'   \\
& & +  \big[F_{3,0}^{IV} +  f''(Q)( A_{2,0} -\frac 13 B_{1,0}^2)\big] (Q_c^3)' + \big[G_{3,0}^{IV}  +   f''(Q) ( \frac 23 B_{2,0} + \frac 13 A_{2,0}B_{1,0})\big] (Q_c^3)'' \\
& &  + O(c^3Q_c),\\
& = &  \frac 12\big[ f''(Q)( 2A_{1,0} +A_{1,0}^2)\big]'Q_c^2 \\
& & + \frac 12 \big[ f''(Q)( 2A_{1,0} +A_{1,0}^2) +  (f''(Q)( B_{1,0} + A_{1,0}B_{1,0}))' \big] (Q_c^2)'\\
& & + F_{2,1}^{IV} cQ_c^2 + [ G_{2,1}^{IV} + \frac 12 f''(Q) ( 2A_{1,1} + B_{1,0}^2) ](cQ_c^2)'   + F_{3,0}^{IV} Q_c^3 \\
& & +  \big[ G_{3,0}^{IV} +  f''(Q)(A_{2,0} -\frac 13 B_{1,0}^2)\big] (Q_c^3)'   %+   f''(Q) ( \frac 23 B_{2,0} + \frac 13 A_{2,0}B_{1,0})\big] (Q_c^3)' 
  + O(c^3Q_c),
\eee
where $ F_{2,1}^{IV}, G_{2,1}^{IV}, F_{3,0}^{IV} $ and $G_{3,0}^{IV}$ are  $\mathcal Y$-functions provided property {\bf (IP)} holds for $k+l\leq 2$. We finally get the Lemma in the quadratic case, the decomposition (\ref{IV}), (\ref{145}) and (\ref{146}), with the desired properties.

\smallskip

\noindent
\emph{Case $m=3$}. Here we consider up to fourth order in our computations. First of all, we write
\bee
{\bf \tilde{IV}}& = &  \frac 12f''(Q)(2Q_cW + W^2) +\frac 12 f^{(3)}(Q)(Q_c^2W+W^2Q_c+\frac 13W^3 ) \\
& & + \frac 12 f^{(4)}(Q) (\frac 13 Q_c^3W +\frac 12 Q_c^2 W^2 +\frac 13 Q_c W^3 + \frac 1{12} W^4) + O(Q_c^5)\\
& =: & {\bf \tilde{IV}_2} + {\bf \tilde{IV}_3} + {\bf \tilde{IV}_4} + O(Q_c^5).
\eee

From now on, and for the sake of simplicity in our computations, we will consider that property {\bf (IP)} holds for any $(k,l)\in \Sigma_3$ and that $A_{1,0}\in \mathcal{Y}$. We recall that in the cubic case, our correction term is given by
\bee
W & = & A_{1,0}Q_c + B_{1,0}Q_c' + A_{1,1}cQ_c +B_{1,1}cQ_c' \\  
& & + A_{2,0}Q_c^2 + B_{2,0}(Q_c^2)' + A_{3,0} Q_c^3 + B_{3,0}(Q_c^3)'\\
& & + A_{2,1}cQ_c^2 + B_{2,1}c(Q_c^2)' + A_{4,0}Q_c^4 + B_{4,0}(Q_c^4)'.
\eee

Let us first consider the term ${\bf \tilde{IV}_2}$.
Note that 
\bee
Q_c W & = & A_{1,0}Q_c^2 + \frac 12 B_{1,0}(Q_c^2)' +  A_{1,1}cQ_c^2 + \frac 12 B_{1,1}c(Q_c^2)' + A_{2,0}Q_c^3 + \frac 23 B_{2,0}(Q_c^3)'\\
& & + A_{3,0}Q_c^4 + \frac 34B_{3,0} (Q_c^4)' +O(cQ_c^3 + Q_c^5).
\eee
Using (\ref{taylor0}) we get
\bee
W^2 & = & A_{1,0}^2 Q_c^2 + B_{1,0}^2 Q_c'^2  + A_{2,0}^4Q_c^4 \\
& & + A_{1,0}B_{1,0}(Q_c^2)' + 2A_{1,0}A_{1,1} cQ_c^2 + A_{1,0}B_{1,1}c(Q_c^2)' \\
& & + 2A_{1,0}A_{2,0}Q_c^3 + \frac 43 A_{1,0}B_{2,0}(Q_c^3)' + 2A_{1,0}A_{3,0}Q_c^4 + \frac 32 A_{1,0}B_{3,0} (Q_c^4)' \\
& & +   A_{1,1}B_{1,0}c(Q_c^2)'   +  \frac 23 A_{2,0}B_{1,0}(Q_c^3)' + \frac 12 A_{3,0}B_{1,0}(Q_c^4)' + O(cQ_c^3 + Q_c^5) \\
& = & A_{1,0}^2 Q_c^2 + A_{1,0}B_{1,0} (Q_c^2)'  +( 2A_{1,0}A_{1,1} +B_{1,0}^2) cQ_c^2\\
& & + (A_{1,1}B_{1,0} + A_{1,0}B_{1,1})c(Q_c^2)' +  2A_{1,0}A_{2,0}Q_c^3 + \frac 23( 2A_{1,0}B_{2,0}+   A_{2,0}B_{1,0})(Q_c^3)' \\
& & + ( 2A_{1,0}A_{3,0} + A_{2,0}^4 - \frac 12 B_{1,0}^2)Q_c^4 + \frac 12 ( 3A_{1,0}B_{3,0}+  A_{3,0}B_{1,0})(Q_c^4)' + O(cQ_c^3 + Q_c^5).
\eee

From here, and using the {\bf (IP)} property, we get (note that $f''(Q)\in \mathcal{Y}$)
\bee
{\bf \tilde{IV}_2} & = & \frac12 f''(Q) ( 2A_{1,0} + A_{1,0}^2)Q_c^2 + \frac12 f''(Q)(B_{1,0}+ A_{1,0}B_{1,0}) (Q_c^2)' \\
 & & +  \sum_{(k,l)\in \Sigma_3^{'}}c^l ( F_{k,l}^{\tilde{IV}_2} Q_c^k+G_{k,l}^{\tilde{IV}_2}(Q_c^k)' ) + O(cQ_c^3 + Q_c^5),
\eee
where $\Sigma_3^{'}$ was introduced in (\ref{Sigm}), and  $F_{k,l}^{\tilde{IV}_2}, G_{k,l}^{\tilde{IV}_2}\in \mathcal{Y}$.

Now we deal with ${\bf \tilde{IV}_3}$. Here we have
$$
Q_c^2 W = A_{1,0}Q_c^3 + \frac 13 B_{1,0}(Q_c^3)' +  A_{2,0}Q_c^4 + \frac 14 B_{2,0}(Q_c^4)'  + O(cQ_c^3 + Q_c^5).
$$
and
\bee
Q_cW^2 &=&A_{1,0}^2 Q_c^3 + \frac 23 A_{1,0}B_{1,0} (Q_c^3)' +  2A_{1,0}A_{2,0}Q_c^4 \\
& &   + \frac 12( 2A_{1,0}B_{2,0}+   A_{2,0}B_{1,0})(Q_c^4)'  + O(cQ_c^3 + Q_c^5).
\eee
Finally
\bee
W^3 & = & W^2 W\\
& = & A_{1,0}^3 Q_c^3 + \frac 23 A_{1,0}^2B_{1,0} (Q_c^3)' +  2A_{1,0}^2A_{2,0}Q_c^4 + \frac 12( 2A_{1,0}^2B_{2,0}+   A_{1,0}A_{2,0}B_{1,0})(Q_c^4)'  \\
& &   + \frac 13 B_{1,0} A_{1,0}^2 (Q_c^3)'  +  \frac 12A_{1,0}B_{1,0} A_{2,0}(Q_c^4)'  \\
& & + A_{2,0}A_{1,0}^2 Q_c^4 +  \frac 12  A_{2,0}A_{1,0}B_{1,0}(Q_c^4)'  +  \frac 12 A_{1,0}^2 B_{2,0} (Q_c^4)' + O(cQ_c^3 + Q_c^5) \\
& = & A_{1,0}^3 Q_c^3  +  A_{1,0}^2B_{1,0}  (Q_c^3)' + 3A_{1,0}^2A_{2,0} Q_c^4   + \frac 32( A_{1,0}^2B_{2,0}+   A_{1,0}A_{2,0}B_{1,0}  )(Q_c^4)' \\
& &+ O(cQ_c^3 + Q_c^5).
\eee
From here, and using the {\bf (IP)} property, we get 
% Q_c^2 W + W^2Q_c + \frac 13 W^3
\bee
{\bf \tilde{IV}_3} & = & \frac12 f^{(3)}(Q) [ \frac 13 B_{1,0}(Q_c^3)' +  A_{2,0}Q_c^4 + \frac 14(  2A_{2,0}B_{1,0} + B_{2,0})(Q_c^4)' ] \\
 & & +  \sum_{(k,l)\in \Sigma_3^{'}}c^l ( F_{k,l}^{\tilde{IV}_2} Q_c^k+G_{k,l}^{\tilde{IV}_2}(Q_c^k)' ) + O(cQ_c^3 + Q_c^5),
\eee
where  $\Sigma_3^{'}$ was introduced in (\ref{Sigm}), and  $F_{k,l}^{\tilde{IV}_3}, G_{k,l}^{\tilde{IV}_3}\in \mathcal{Y}$.

\smallskip

Finally, fourth order terms are easy to compute:
$$
Q_c^3 W = A_{1,0}Q_c^4 + \frac 14 B_{1,0}(Q_c^4)' + O(Q_c^5 + cQ_c^3) ,
$$
$$
Q_c^2 W^2 = A_{1,0}^2 Q_c^4 + \frac 12A_{1,0}B_{1,0} (Q_c^4)' + O(Q_c^5 + cQ_c^3) ,
$$
$$
Q_c W^3 = A_{1,0}^3 Q_c^4 + \frac 34 A_{1,0}^2 B_{1,0} (Q_c^4)' + O(Q_c^5 + cQ_c^3),
$$
and
$$
W^4 = A_{1,0}^4 Q_c^4 + A_{1,0}^3 B_{1,0} (Q_c^4)' + O(Q_c^5 + cQ_c^3 ).
$$
As we have supposed $A_{1,0}\in \mathcal{Y}$ ($(1,0)$ satisfies {\bf (IP)}), we will obtain
$$
{\bf \tilde{IV}_4} =  F_{4,0}^{ \tilde{IV}_4} Q_c^4 +   [ G_{4,0}^{ \tilde{IV}_4} + \frac 1{24} f^{(4)}(Q) B_{1,0} ](Q_c^4)' + O(Q_c^5 + cQ_c^3), 
$$
where $F_{4,0}^{ \tilde{IV}_4}, G_{4,0}^{ \tilde{IV}_4}\in \mathcal{Y}$.

We finally collect  the expansions of ${\bf \tilde{IV}_2}, {\bf \tilde{IV}_3}$ and ${\bf \tilde{IV}_4}$. We derivate to obtain
\bee
{\bf IV} & = & ({\bf \tilde{IV}_2}+{\bf \tilde{IV}_3}+{\bf \tilde{IV}_4} +O(Q_c^5) )_x\\
&  &  \frac12  (f''(Q)  ( 2A_{1,0} + A_{1,0}^2) )' Q_c^2   + \frac12 f^{(3)}(Q) A_{2,0}(Q_c^4)'   \\
& &  + \frac12 \big[ f''(Q) ( 2A_{1,0} + A_{1,0}^2) + (f''(Q)(B_{1,0}+ A_{1,0}B_{1,0}) )'  \big] (Q_c^2)'  \\
 & & +  \sum_{(k,l)\in \Sigma_3^{'}}c^l ( F_{k,l}^{\tilde{IV}} Q_c^k+G_{k,l}^{\tilde{IV}}(Q_c^k)' ) + O(cQ_c^3 + Q_c^5),
\eee
where, as we have emphasized, $F_{k,l}^{\tilde{IV}}, G_{k,l}^{\tilde{IV}}\in \mathcal{Y}$  provided $A_{k',l'}, B_{k',l'}$ satisfy the {\bf (IP)} property for $(k',l')< (k,l)$, as is the case. Here we have also used that $(f''(Q))', (f^{(4)}(Q))' \in \mathcal{Y}$ for all $p\geq 4$. The set $ \Sigma_3^{'}$ was defined in (\ref{Sigm}).

\medskip

Let us finally prove (v). From (i), the rest term $\mathcal{E}(t,x)$
is a finite sum of terms of the type $c^l\wqs^k(\ys)f(y)$ or $c^l(\wqs^k)'(\ys)f(y)$, where $(k,l)\not\in \Sigma_m$. More specifically, this means $k+l\geq 4$ for $m=2$ and $(k,l)=(1,2), (3,1)$ or higher order terms (excluding $(k,l)=(4,0)$) in the case $m=3$ (see the definition of $\Sigma_m$ in Section \ref{sec:2-1}). Here  $f$ is a bounded function such that $f'\in \mathcal{Y}$. Thus, we easily conclude, using Claim \ref{cl:30},
$$ \norm{\mathcal{E}(t)}_{H^1(\R)} \leq K c^{3(\frac 1{m-1} +\frac 14)},$$
as desired. This finishes the proof.

\end{proof}

Putting together Lemmas \ref{lem:SQ}--\ref{lem:Sint}, we obtain Proposition \ref{prop:decompbis}, in
particular, the explicit expressions of $F_{k,l}$ and $G_{k,l}$ for $1\leq k+l\leq 2$.

\end{proof}

\section{End of proof of Proposition \ref{lem:b20}}\label{D}

Continuing with the proof of Proposition \ref{lem:b20}, we show now the existence of a nonzero residual term appearing after the collision.

\subsection{General computations}

We proceed to compute the constants $b_{2,0}$ more explicitly. In the course of the proof we will made use several times of the equations satisfied by the functions $A_{k,l}, B_{k,l}, a_{k,l}$ for $(k,l)=(1,0)$ and $(2,0)$, cf. (\ref{A10})-(\ref{B10}) for the system $(\Omega_{1,0})$ and (\ref{eq:A20})-(\ref{eq:B20}) for the second one.

\begin{Cl}[Explicit value of $b_{2,0}$]\label{CL}
Suppose $f$ as in (\ref{surf}). Then the following expressions for the $b_{2,0}$ coefficient hold.
\begin{enumerate}
\item Case $m=2$.
\bea\label{b202new}
b_{2,0}& = & -\frac 12 b_{1,0}^3 + \frac 14 \int_\R (f''(Q) -2)(1+A_{1,0})^3   -2b_{1,0}  +  \frac 12\int_\R A_{1,0}(1+A_{1,0}^2) \nonumber \\
& & \quad  - \frac 12 a_{1,0}\int_\R Q A_{1,0} -\frac 34 a_{1,0}^3 \int_\R Q'^2 + \frac 12 a_{1,0}^2\int_\R Q[Q -f'(Q)(1+A_{1,0}) ] \nonumber \\
& & \quad -\frac 34a_{1,0} \int_\R [ f'(Q)(1+A_{1,0}) + 3A_{1,0}''  ]A_{1,0}   + \frac 12 \int_\R B_{1,0}[ 3A_{1,0}' + f'(Q)\int_0^x (A_{1,0}  + a_{1,0}Q) ] \nonumber \\
& & \quad  + 3 a_{1,0}^2 \int_\R Q''A_{1,0} .
\eea

\item Case $m=3$.
\bea\label{b20new}
b_{2,0} & = &  \frac 14 \int_\R f''(Q)(1+A_{1,0})^3 -\frac 34 a_{1,0}\int_\R f'(Q)(1+A_{1,0})A_{1,0} + \frac 94 a_{1,0} \int_\R A_{1,0}'^2  \nonumber \\
& & -\frac 12 a_{1,0}^2 \int_\R f'(Q)Q (1+A_{1,0}) + 3a_{1,0}^2 \int_\R A_{1,0}Q'' -\frac 34 a_{1,0}^3 \int_\R Q'^2 .
\eea 
\end{enumerate}
\end{Cl}

\begin{proof}

We treat first the cubic case, being easier. Let us start with (\ref{a20}) and (\ref{b20}). In this case, we have a priori chosen $A_{2,0}\in \mathcal{Y}$, so that  $\ga_{2,0}=0$, and then from (\ref{eq:solOmega10})
\be\label{b3}
b_{2,0} = \frac 12\Big[  -a_{1,0} \int_\R G_{2,0}Q - a_{1,0}\int_\R \tilde F_{2,0}  P  + \int_\R  \tilde F_{2,0} \bar P  + \int_\R G_{2,0} \Big],
\ee 
where $\tilde F_{2,0}' = F_{2,0}$, $\tilde F_{2,0}\in \mathcal{Y}$. More precisely,
$$
\tilde F_{2,0} :=  \frac 12 f''(Q)(1+A_{1,0})^2 +3a_{1,0}^2 Q'' -a_{1,0} (3A_{1,0}'' +f'(Q)(1+A_{1,0})).
$$
First, it is easy to see from (\ref{G203}) by using the {\bf (IP)} property for $(k,l)=(1,0)$, that 
\be\label{G1}
\int_\R G_{2,0} =  \frac 12\int_\R f''(Q)(1+A_{1,0})^2 .
\ee
Secondly, from (\ref{A}), (\ref{P}) and (\ref{hatp}),
$
\bar P -a_{1,0}P = A_{1,0} +a_{1,0}Q,
$
and thus
$$
b_{2,0} = \frac 12\Big[  a_{1,0}\int_\R (\tilde F_{2,0} -G_{2,0}) Q  + \int_\R  \tilde F_{2,0}A_{1,0}  + \frac 12 \int_\R  f''(Q)(1+A_{1,0})^2 \Big].
$$
It is clear that
\bee
\tilde F_{2,0} -G_{2,0} &  = &  \frac 32 a_{1,0}^2 Q''  + a_{1,0} [ \frac 32  A_{1,0}'' - f'(Q)(1+A_{1,0}) + \frac 32 B_{1,0}^{(3)} + \frac 12 (f'(Q)B_{1,0})' ] \\
& & - \frac 12 (f''(Q)(1+A_{1,0})B_{1,0})' .
\eee
From here, after several integration by parts,
\begin{eqnarray}
\int_\R (\tilde F_{2,0} -G_{2,0}) Q &  = &  -\frac 32 a_{1,0}^2 \int_\R Q'^2  + a_{1,0}\int_\R [ \frac 32  A_{1,0}'' - f'(Q)(1+A_{1,0})]Q \nonumber \\
& & - \frac 12 a_{1,0} \int_\R B_{1,0} [ 3 Q'' + f(Q) ]' + \frac 12 \int_\R B_{1,0} (f'(Q))'(1+A_{1,0}) \nonumber \\
& = &  -\frac 32 a_{1,0}^2 \int_\R Q'^2  + a_{1,0}\int_\R [ \frac 32  A_{1,0}'' - f'(Q)(1+A_{1,0})]Q  \nonumber \\
& & +\frac 12 \int_\R B_{1,0}[(f'(Q))'(1+A_{1,0})  - a_{1,0} ( 3Q -2f(Q))' ]. \label{BF}
\end{eqnarray}
But from (\ref{A10}), 
$ (\mathcal L(1+A_{1,0}))' = (1-f'(Q) + \mathcal L A_{1,0})' %= -(f'(Q))' -a_{1,0}(3Q-2f(Q))' + (f'(Q))' 
= -a_{1,0}(3Q-2f(Q))' . $ On the other hand, expanding $(\mathcal L(1+A_{1,0}))'$, we get
$$
(\mathcal L(1+A_{1,0}))' = -A_{1,0}^{(3)} +A_{1,0}'  - (f'(Q))'(1+A_{1,0}) - f'(Q)A_{1,0}' .
$$
From here, the quantity in front of $B_{1,0}$ in (\ref{BF}) is nothing but $\mathcal L A_{1,0}'$. Coming back to (\ref{BF}), and using the equation for $B_{1,0}$ (\ref{B10}), we obtain
\bea
\int_\R (\tilde F_{2,0} -G_{2,0}) Q & = &  -\frac 32 a_{1,0}^2 \int_\R Q'^2  + a_{1,0}\int_\R [ \frac 32  A_{1,0}'' - f'(Q)(1+A_{1,0})]Q -\frac 12 \int_\R A_{1,0} (\mathcal L B_{1,0})'\nonumber \\
& = &  -\frac 32 a_{1,0}^2 \int_\R Q'^2 + \frac 32 \int_\R A_{1,0}'^2 + 3a_{1,0}\int_\R Q''A_{1,0} \nonumber \\
& & -a_{1,0}\int_\R f'(Q)Q(1+A_{1,0})  -\frac 12 \int_\R f'(Q)(1+A_{1,0})A_{1,0}.\label{FQ} 
\eea

Finally, an easy computation shows that
\bea
\int_\R \tilde F_{2,0} A_{1,0} & = & \frac 12 \int_\R f''(Q)(1+A_{1,0})^2 A_{1,0} +3a_{1,0}^2 \int_\R A_{1,0}Q'' \nonumber\\ 
& & -a_{1,0}\int_\R (3A_{1,0}'' +f'(Q)(1+A_{1,0}))A_{1,0}\label{FA}
\eea
Collecting (\ref{FQ}) and (\ref{FA}), we get 
\bee
b_{2,0} & = &  \frac 14 \int_\R f''(Q)(1+A_{1,0})^3 -\frac 34 a_{1,0}\int_\R f'(Q)(1+A_{1,0})A_{1,0} + \frac 94 a_{1,0} \int_\R A_{1,0}'^2  \\
& & -\frac 12 a_{1,0}^2 \int_\R f'(Q)Q (1+A_{1,0}) + 3a_{1,0}^2 \int_\R A_{1,0}Q'' -\frac 34 a_{1,0}^3 \int_\R Q'^2 .
\eee 
as desired.

\smallskip

Let us treat now the quadratic case. The procedure is similar, but more involved. Now we assume that $\ga_{2,0}= -\frac 12 b_{1,0}^2$ as in Proposition \ref{lem:b20} (i), and consider (\ref{a20})-(\ref{b20}). We get
$$
b_{2,0} = \frac 12\Big[ -\frac 12 b_{1,0}^2\int_\R (\bar P -a_{1,0}P)  +\int_\R F_{2,0}\int_0^x (a_{1,0}P-\bar P) +\int_\R G_{2,0} -a_{1,0}\int_\R G_{2,0}Q \Big].
$$
Now several remarks. Note that from (\ref{A}) and the definition of $P$ in (\ref{P}) we have $\bar P -a_{1,0} P = A_{1,0} + a_{1,0}Q$, and from (\ref{balt}),
$$
\int_\R (\bar P -a_{1,0}P) = 2b_{1,0}.
$$
Second, note that from $b_{1,0}= \pm \lim_{\pm \infty} B_{1,0}$ and $\lim_{\pm \infty } f''(Q) = 2$, 
$$
\int_\R G_{2,0} = \frac 12 \int_\R (f''(Q) -2) -\int_\R A_{1,0} -4b_{1,0}+\frac 12 \int_\R f''(Q) (2A_{1,0}+ A_{1,0}^2).
$$
On the other hand, from (\ref{F202}), $F_{2,0} = \tilde F_{2,0}' -f'(Q)B_{1,0} $, where $\tilde F_{2,0}\in \mathcal{Y}$ and is given by
\bee
\tilde F_{2,0} & : = &  -(3A_{1,0} +3B_{1,0}') +\frac 12 f''(Q)(2A_{1,0} +A_{1,0}^2) \nonumber \\
& & \quad  - a_{1,0} ( 3A_{1,0}''  -Q +f'(Q)(1+A_{1,0})  ) + 3a_{1,0}^2 Q'' + \frac 12 (f''(Q)-2).
\eee
Thus, 
\bee
\int_\R F_{2,0}\int_0^x (a_{1,0}P-\bar P) & = & \int_\R \tilde F_{2,0}(\bar P -a_{1,0}P) +\int_\R f'(Q)B_{1,0}\int_0^x (\bar P -a_{1,0}P)\\
& =&  \int_\R \tilde F_{2,0}(A_{1,0} + a_{1,0}Q) +\int_\R f'(Q)B_{1,0}\int_0^x (A_{1,0} +a_{1,0} Q).
\eee

Repeating the same computation for the cubic case, we obtain
\bee
\int_\R Q(\tilde F_{2,0}-G_{2,0}) %& = &% - 2\int_\R Q A_{1,0} -\frac 32 a_{1,0}^2 \int_\R Q'^2 + a_{1,0}\int_\R Q(\frac 32 A_{1,0}'' +Q -f'(Q)(1+A_{1,0}) )\\
%& & -\frac 12 \int_\R (\mathcal L B_{1,0})'A_{1,0}\\
& = &  - 2\int_\R Q A_{1,0} -\frac 32 a_{1,0}^2 \int_\R Q'^2 + a_{1,0}\int_\R Q(\frac 32 A_{1,0}'' +Q -f'(Q)(1+A_{1,0}) )\\
& & -\frac 12 \int_\R [f'(Q)(1+A_{1,0}) +3A_{1,0}'' -3a_{1,0}Q'' ]A_{1,0},
\eee 
and 
\bee
\int_\R \tilde F_{2,0}A_{1,0} & = & -\int_\R A_{1,0}(3A_{1,0} +3B_{1,0}') +\frac 12\int_\R f''(Q)(2A_{1,0} +A_{1,0}^2)A_{1,0}+ 3a_{1,0}^2 \int_\R Q''A_{1,0} \nonumber \\
& & \quad  - a_{1,0}\int_\R ( 3A_{1,0}''  -Q +f'(Q)(1+A_{1,0})  ) A_{1,0} + \frac 12 \int_\R (f''(Q)-2)A_{1,0}.
\eee
Collecting the above identities and after several simplifications we get
\bee
b_{2,0}& = & -\frac 12 b_{1,0}^3 + \frac 14 \int_\R (f''(Q) -2)(1+A_{1,0})^3   - 2b_{1,0}  +  \frac 12 \int_\R A_{1,0}(1+A_{1,0}^2) \\
& &  - \frac 12 a_{1,0}\int_\R Q A_{1,0} -\frac 34 a_{1,0}^3 \int_\R Q'^2 + \frac 12 a_{1,0}^2\int_\R Q(Q -f'(Q)(1+A_{1,0}) )\\
& & -\frac 34a_{1,0} \int_\R [ f'(Q)(1+A_{1,0}) + 3A_{1,0}''  ]A_{1,0}   +\frac 12 \int_\R B_{1,0}[ 3A_{1,0}' + f'(Q)\int_0^x (A_{1,0}  + a_{1,0}Q) ] \\
& &  + 3 a_{1,0}^2 \int_\R Q''A_{1,0} .
\eee
The proof is now complete.
\end{proof}

The objective is now to give the first order terms for the coefficient $b_{2,0}$. For this, we consider separate cases. 
%%%%%%%%%%%%%%%%%%%%%%
%\begin{defn}
%Let $a=a(\ve)$ be a real number and $f =f_\ve$ a real valued function, both depending smoothly on $\ve$. We define
%$$
%a^0 := \lim_{\ve\to 0} a, \quad f^0 := \lim_{\ve\to 0} f_\ve,
%$$ 
%where the limits are understood as pointwise. We also denote
%$$
%a^1 := \lim_{\ve\to 0} \frac 1\ve (a-a^0), \quad f^1 := \lim_{\ve\to 0} \frac 1\ve (f_\ve-f^0).
%$$ 
%\end{defn}
It turns out that computations in the cubic case are easy to carry out. We first deal with this case.

\subsection{Cubic case}

The objective of this paragraph is to prove the following

\begin{lem}[Asymptotic expansions, case $m=3$]\label{L14} We have
$$
b_{2,0} = b_{2,0}^1 \ve + o(\ve).
$$
where 
\be\label{b201}
b_{2,0}^1 =: c_{3,p} = -\Big[ \frac{(p-1)(p-3)(p^2-3p+8) }{8(p-2)(p+1)}\Big]\int_\R (Q^0)^p. 
\ee
In particular, for any $p\geq 4$, $b_{2,0}(\ve) \neq 0$ provided $0<\abs{\ve} \leq \ve_0$ for $\ve_0$ small. 
\end{lem}

First of all we start with an auxiliary
\begin{Cl}[Asymptotic expansions, basic functions]\label{exp6}

Suppose $f$ as in (\ref{gKdV0b}), $p\geq 4$. The following asymptotic expansions hold. 
\begin{enumerate}
\item The soliton solution $Q$ can be expanded as
\be\label{ffp1}
Q= Q^0 + \ve Q^1 +o(\ve), \quad o(\ve) \in \mathcal{Y},
\ee
where $Q^0$ and $Q^1$ satisfy the equations
\be\label{ffp2}
-(Q^0)'' + Q^0 - (Q^0)^3 =0, \quad  \mathcal L^0 Q^1 :=-(Q^1) '' + Q^1 -  3(Q^0)^2  Q^1= (Q^0)^p.
\ee
Finally, 
\be\label{ffp3}\begin{cases}
& f(Q) = (Q^0)^3 + \ve(3(Q^0)^2 Q^1 + (Q^0)^p ) + o(\ve), \\
& f'(Q) = 3(Q^0)^2 +\ve(6Q^0Q^1 + p(Q^0)^{p-1}) + o(\ve),\\
& f''(Q) = 6Q^0 + \ve( 6Q^1 + p(p-1) (Q^0)^{p-2}) + o(\ve),
\end{cases}
\ee
where every term $o(\ve)\in \mathcal{Y}$ uniformly in $\ve<\ve_0$.
\item The operator $\mathcal L$ satisfies
\be\label{ffp4}
\mathcal L = \mathcal L^0  - \ve [ 6Q^0Q^1 + p(Q^0)^{p-1}] + o(\ve), \quad  \mathcal{L}^{0} = -\partial_x^2 + 1 - 3(Q^0)^2.
\ee
\item From (\ref{P}), (\ref{barp}) and (\ref{hatp}), the test functions $P, \bar P$ and $\hat P$ satisfy the following relations
\be\label{Ps3}
\begin{cases}
\Lambda Q = \Lambda Q^0 + \ve \Lambda Q^1 +o(\ve) \in \mathcal Y\;   \hbox{ where } \Lambda Q^0 :=  \frac 12(  x(Q^0)' + Q^0),  \quad \mathcal{L}^{0} \Lambda Q^0 = -Q^0;\\
\hbox{ and } \mathcal{L}^{0} \Lambda Q^1 := (6Q^0 \Lambda Q^0 -1) Q^1 + p(Q^0)^{p-1} \Lambda Q^0.
\end{cases}
\ee
Moreover, the following identities hold
\bea
\int_\R \Lambda Q^0 &= &0,  \nonumber  \\
\int_\R \Lambda Q^1 &=&  \int_\R[ -1+ (Q^0)^2 +6Q^0\Lambda Q^0 - 6(Q^0)^3 \Lambda Q^0 ] Q^1 \nonumber \\
& & + p\int_{\R} (Q^0)^{p-1} \Lambda Q^0 (1-(Q^0)^2).\label{canc}
\eea
\item \emph{Integrals}. 
For any $p\geq 1$,
\be\label{Int23}
\int_\R (Q^0)^{p+2} = \frac{2p}{1+ p}\int_{\R} (Q^0)^p, \quad \int_\R (Q^0)^2 = 4.
\ee
\item Let  $D(\ve) = \int_\R \Lambda Q Q$. Then
\be\label{int23}
D(\ve) = 1+ O(\ve).
\ee
\item \emph{Inverse functions}. The following identities hold
\be\label{inv1}
\mathcal{L}^0 (-\frac 94 x(Q^0)' -\frac{15}4 Q^0 + \frac 32 (Q^0)^3) = \frac 92 Q^0(1-(Q^0)^2)^2, 
\ee
\be\label{inv2}
\mathcal{L}^0 (xQ^0 (Q^0)') = -4(Q^0)^2 +3(Q^0)^4 -3xQ^0(Q^0)'(1-(Q^0)^2), 
\ee
\be\label{inv3}
\mathcal{L}^0 ((Q^0)^4) =  -15(Q^0)^4 + 7(Q^0)^6.
\ee
\end{enumerate}
\end{Cl}

\begin{proof}
First of all, (\ref{ffp1})-(\ref{ffp4}) follow by Taylor expansion in $\ve$. Concerning (\ref{Ps3}), it follows from (\ref{ffp1})-(\ref{ffp4}).
Let us see (\ref{canc}). From the definition of $\mathcal{L}^0 \Lambda Q^1$ and the identity $\mathcal{L}^0 (Q^0)^2 = -3(Q^0)^2$, we have
$$
\int_\R \mathcal{L}^0 \Lambda Q^1 = \int_\R \Lambda Q^1 + \int_\R (Q^0)^2 \mathcal{L}^0 \Lambda Q^1,
$$
thus
$$
\int_\R \Lambda Q^1 =  \int_\R[ 1- (Q^0)^2] \mathcal{L}^0 \Lambda Q^1= \int_\R[ 1- (Q^0)^2] [ (6Q^0 \Lambda Q^0 -1) Q^1 + p(Q^0)^{p-1} \Lambda Q^0],
$$
where we obtain (\ref{canc}). 

To obtain (\ref{Int23}) we use integration by parts and the explicit function $Q^0(x):= \frac{\sqrt{2}}{\cosh x}$.

We prove (\ref{int23}). It follows from the fact that
$$
\int_\R Q^0\Lambda Q^0 = \frac 12 \int_\R (\frac 12 x ((Q^0)^2)' + (Q^0)^2) = \frac 14 \int_\R (Q^0)^2 =1.
$$
Finally, (\ref{inv1})-(\ref{inv3}) are obtained by simple differentiation. We left the proof to the reader. %The proof is now complete.
\end{proof}

\begin{Cl}[Asymptotic expansions, case $m=3$]\label{exp3}
The following expansion hold. 
\be\label{a10new}\begin{cases}
a_{1,0} =  a_{1,0}^0 + \ve a_{1,0}^1 +o(\ve),  \quad a_{1,0}^0 =0,\\
A_{1,0} = A_{1,0}^0 + \ve A_{1,0}^1 + o(\ve), \quad o(\ve)\in \mathcal{Y}, \quad A_{1,0}^0=  -(Q^0)^2, \\
B_{1,0} =  B_{1,0}^0 + \ve B_{1,0}^1 +o(\ve), \quad B_{1,0}^0 = -2\varphi^0 -\frac 34\sqrt{2} \pi (Q^0)' .
\end{cases}
\ee
Here $a_{1,0}^1 := \int_{\R} \Lambda Q^1$ and $A_{1,0}^1, B_{1,0}^1$ satisfy the following linear system
\be\label{AB13}\begin{cases}
(\mathcal L^{0} A_{1,0}^1)' + a_{1,0}^1(3Q^0 - 2(Q_0)^3)' = ((6Q^0Q^1+p(Q^0)^{p-1}) (1+ A_{1,0}^0))',\\
(\mathcal L^{0} B_{1,0}^1)' +3 a_{1,0}^1(Q^0)'' -3(A_{1,0}^1)'' - 3(Q^0)^2 A_{1,0}^1 =  (6Q^0Q^1 + p(Q^0)^{p-1})(1+ A_{1,0}^0).
\end{cases}
\ee
\end{Cl}

\begin{proof}
We start with the zeroth order system. From (\ref{A10})-(\ref{B10}) and using Claim \ref{exp6} we get
$$\begin{cases}
(\mathcal L^{0} A_{1,0}^0)' + a_{1,0}^0(3Q^0 - 2(Q^0)^3)' = (3(Q^0)^2)'.\\
(\mathcal L^{0} B_{1,0}^0)' +3a_{1,0}^0(Q^0)'' -3(A_{1,0}^0)'' - 3(Q^0)^2 A_{1,0}^0 = 3(Q^0)^2.
\end{cases}
$$
It is easy to verify that $a_{1,0}^0 =0$, $A_{1,0}^0 = -(Q^0)^2 \in \mathcal{Y}$ and $B^0_{1,0}= -2\varphi^0 -\frac 34\sqrt{2} \pi (Q^0)' $ satisfy this system with the required properties. In particular,
$$
\int_\R B^0_{1,0} (Q^0)' =0.
$$
Concerning the system (\ref{AB13}), it follows directly from (\ref{eq:A20})-(\ref{eq:B20}) and using Claim \ref{exp6}.
We will not solve this system explicitly, but we only compute the constant $a_{1,0}^1$.   

Indeed, from (\ref{eq:solOmega10}) and Claim \ref{exp6}, we have $a_{1,0} = a_{1,0}^0 + a_{1,0}^1\ve + o(\ve), $
where
$$
a_{1,0}^0 := \frac{\int_\R \Lambda Q^0 }{\int_\R \Lambda Q^0 Q^0} = 0, 
\quad \hbox{ and } \quad 
a_{1,0}^1 :=  \frac{\int_\R \Lambda Q^1 }{\int_\R \Lambda Q^0 Q^0} = \int_\R \Lambda Q^1.
$$
This finishes the proof.
\end{proof}
We finally prove Lemma \ref{L14}.

\begin{proof}[Proof of Lemma \ref{L14}]
From (\ref{b20new}) and (\ref{a10new}) we have $b_{2,0} = b_{2,0}^0 + \ve b_{2,0}^1 + o(\ve)$, where
$$b_{2,0}^0 = \frac 14 \int_\R 6Q^0 (1+A_{1,0}^0)^3 =  \frac 14 \int_\R 6Q^0 (1-(Q^0)^2)^3 =0, $$
and 
\bee
b_{2,0}^1&  = &  \frac 14 \int_\R (6Q^1 +p(p-1)(Q^0)^{p-2})(1-(Q^0)^2)^3 +  \frac {9}2 \int_\R Q^0 (1-(Q^0)^2)^2A_{1,0}^1 \\
& & + \frac 94 a_{1,0}^1 \Big[ \int_\R (Q^0)^2 (1-(Q^0)^2)(Q^0)^2 +  \int_\R 4(Q^0)^2((Q^0)^2 - \frac 12 (Q^0)^4)\Big].
\eee
From (\ref{inv1}), the selfadjointness of the operator $\mathcal{L}^0$ and by using (\ref{AB13}), we get
\bee
 \frac 92\int_\R Q^0(1-(Q^0)^2))^2 A_{1,0}^1&  = & \frac 34 \int_\R (-3 x(Q^0)' -5 Q^0 + 2 (Q^0)^3) \mathcal{L}^0 A_{1,0}^1 \\
 &  = & \frac 34  a_{1,0}^1 \int_\R  (3 x(Q^0)' + 5 Q^0 - 2 (Q^0)^3)(3Q^0 -2(Q^0)^3)\\
 & &  - \frac 34  \int_\R (3 x(Q^0)' + 5 Q^0 - 2 (Q^0)^3)  (6Q^0Q^1 + p(Q^0)^{p-1})(1-(Q^0)^2).
\eee
Therefore,
\bee
b_{2,0}^1&  = & \frac 32 \int_\R Q^1 (1-(Q^0)^2)[1-2(Q^0)^2 + (Q^0)^4 -3Q^0(3 x(Q^0)' + 5 Q^0 - 2 (Q^0)^3) ] \\
& & +  \frac p4 \int_\R (Q^0)^{p-2}(1-(Q^0)^2)[(p-1)(1-2(Q^0)^2+(Q^0)^4) -3Q^0(3 x(Q^0)' + 5 Q^0 - 2 (Q^0)^3) ]  \\
& & + \frac 34 \int_\R \Lambda Q^1\Big[  \int_\R (3 x(Q^0)' + 5 Q^0 - 2 (Q^0)^3)(3Q^0 -2(Q^0)^3) + 3\int_\R (Q^0)^4 (5-3(Q^0)^2 )\Big].
\eee
Note that, from (\ref{Int23}) 
$$
\int_\R (3 x(Q^0)' + 5 Q^0 - 2 (Q^0)^3)(3Q^0 -2(Q^0)^3) + 3\int_\R (Q^0)^4 (5-3(Q^0)^2 ) =2, 
$$
thus from (\ref{canc}) we get
\bee
b_{2,0}^1&  = & \frac 32 \int_\R Q^1 (1-(Q^0)^2)[-14(Q^0)^2 + 7(Q^0)^4 -6 xQ^0(Q^0)'  ] \\
& & +  \frac p4 \int_\R (Q^0)^{p-2}(1-(Q^0)^2)[(p-1)(1-2(Q^0)^2+(Q^0)^4) -3Q^0(3 x(Q^0)' + 5 Q^0 - 2 (Q^0)^3) ]  \\
& & + \frac 32  p\int_{\R} (Q^0)^{p-1} \Lambda Q^0 (1-(Q^0)^2)
\eee
Finally, from (\ref{inv2}), (\ref{inv3}) and the identity $\mathcal L^0 (Q^0)^2 = -3(Q^0)^2$, we have
$$
\mathcal L^0 [2(Q^0)^2 -(Q^0)^4 +2xQ^0(Q^0)' ] = (1-(Q^0)^2)[-14(Q^0)^2 + 7(Q^0)^4 -6 xQ^0(Q^0)' ].
$$
Using the selfadjointness of $\mathcal L^0$ and the equation for $Q^1$ in (\ref{ffp2}),  and after integrating by parts, we conclude that
\bee
b_{2,0}^1&  = & \frac 32 \int_\R (Q^0)^p (2(Q^0)^2 -(Q^0)^4 +2xQ^0(Q^0)' )+ \frac 32  p\int_{\R} (Q^0)^{p-1} \Lambda Q^0 (1-(Q^0)^2)  \\
& & +  \frac p4 \int_\R (Q^0)^{p-2}(1-(Q^0)^2)[(p-1)(1-2(Q^0)^2+(Q^0)^4) -3Q^0(3 x(Q^0)' + 5 Q^0 - 2 (Q^0)^3) ]  \\
& = & \frac p4(p-1) \int_\R (Q^0)^{p-2}  -\frac 34 (p^2 +3p-2) \int_\R (Q^0)^{p} + \frac{3}{4(p+2)}(p^3 +7p^2 +12p + 4)\int_\R (Q^0)^{p+2} \\
& &  - \frac 14 (p^2 + 5p +6)\int_\R (Q^0)^{p+4}.   
\eee
Finally, from (\ref{Int23}), and after some simplifications,
\be\label{f2}	
b_{2,0}^1 = -\Big[ \frac{(p-1)(p-3)(p^2-3p+8) }{8(p-2)(p+1)}\Big]\int_\R (Q^0)^p. %-\Big[ \frac{(p-1)(p+2)(p-3)(p+3)(p^2-3p+8) }{8(p-2)(p+1)(p+2)(p+3)}\Big]\int_\R (Q^0)^p
 %-\Big[ \frac{p^6 -2p^5-6p^4 +32p^3 -43p^2 -126p+144}{8(p-2)(p+1)(p+2)(p+3)}\Big]\int_\R (Q^0)^p
\ee
The proof is now complete.
\end{proof}

\begin{rem}
Note that even though the higher regularity needed in our results ($f\in C^{5}$ for $m=3$), we are able to take, at least formally, the limit $p\downarrow 3$ in (\ref{f2}), recovering the results from the integrable case (that is, $b_{2,0}^1=0$). This gain of regularity comes from (\ref{F203}) and (\ref{b20new}): for these identities, we only need $f\in C^3(\R)$.
\end{rem}

\medskip

\subsection{Gardner and quadratic nonlinearities}

These two nonlinearities are very similar to handle. Although computations are harder for the Gardner nonlinearity, a simple trick will allow to link both results. As a consequence, we are reduced to consider only the quadratic case.

Finally, recall the soliton $Q_{\tilde \mu, 1}$ introduced in (\ref{SolG}), well defined for $\tilde \mu<\frac 29$. 
Given $\tilde\mu, \nu \in \R$, $\tilde\mu<\frac 29$ and $\nu$ small enough, let $d_{\tilde\mu,\nu}$ be the \emph{defect} (possibly zero) associated the the nonlinearity $f_{\tilde\mu, \nu} (s):= s^2 -\tilde \mu s^3 + \nu s^p$, namely
\be\label{dmn}
d_{\tilde \mu, \nu} := b_{2,0}(f_{\tilde\mu, \nu}) + \frac 16 b_{1,0}^3(f_{\tilde\mu, \nu}).
\ee

We following reduction Lemma is the key ingredient of the proof.

\medskip

\begin{lem}[]~

Let $d(\ve)$ be the defect parameter introduced in (\ref{diff0}) for the nonlinearity $f(s)$ described in (\ref{gKdV0b}), $m=2$, and let $d_{\tilde \mu, \nu}$ be the defect introduced in (\ref{dmn}), for $\tilde \mu, \nu$ small. Then the following properties are satisfied:
\begin{enumerate}
\item For all $\tilde \mu <\frac 29$, $\nu\in \R$ small, $d_{\tilde \mu, \nu}$ is a smooth function of $\tilde \mu, \nu$ and for all $\tilde \mu <\frac 29$,
\be\label{Nodef}
d_{\tilde\mu, 0} =0.
\ee
\item Given $\ve$ small, let $\tilde \mu = \mu(\ve)$ and $\nu =\ve$. Then the following expansion holds
\be\label{K1}
d_{\mu(\ve),\ve} = -\ve\Big[ \frac{(p-3)(2p-1) (48 -46p +6 p^2  +4p^3)}{72(p^2 -1)(p-2)}\Big] \int_\R \Big[\frac{3}{2\cosh^2 (x/2)}\Big] ^p+ o(\ve), %\quad c_{2,p}= \int_\R (Q^0)^p \neq 0,
\ee
for all $\abs{\ve}<\ve_0$ and $p\geq 3$.  %Here $c_{2,p}$ is given by (\ref{c2p}).
\item The following expansion holds
\be\label{Red1}
d(\ve) = d_{\mu(\ve), \ve} + o(\ve), \quad \hbox{ as } \ve \to 0.
\ee  
\end{enumerate}
\end{lem}

\begin{proof}[Proof of (\ref{Red1})]
This is an easy consequence of the definition of $f$ in (\ref{gKdV0}), and the fact that $f(Q) = f_{\mu(\ve), \ve}(Q) + o(\ve)$, with $o(\ve) \in \mathcal Y$. In particular, the soliton $Q$ and each term $a_{1,0}, b_{1,0}, A_{1,0}, B_{1,0}$ and $b_{2,0}$ depends smoothly in $\ve$ and can be expanded in a similar way.
\end{proof}

\begin{proof}[Proof of (\ref{Nodef})]
The smoothness is a direct consequence of the formula for $b_{2,0}$ in Claim \ref{CL} and $b_{1,0}$ in (\ref{balt}).  We have to prove that for all $\tilde \mu<\frac 29$,
$$
b_{2,0} (f_{\tilde\mu, 0}) + \frac 16 b_{1,0}^3(f_{\tilde\mu, 0})  =0,
$$
In order to prove this identity, we claim the following

\begin{Cl}[Basic functions]
Let $Q^0 := Q_{\tilde \mu,1 }$ be the soliton for the Gardner equation. Then we have
\begin{enumerate}
\item The soliton solution $Q^0$ satisfies
$$
-(Q^0)'' + Q^0 - (Q^0)^2 + \tilde \mu (Q^0)^3 =0, \quad \mathcal L^0 (Q^0)' =0,
$$
where $\mathcal L^0 := -\partial_{xx} + 1 -(2Q^0 -3\tilde \mu (Q^0)^2) $.
\item From the definition of $Q^0$, we have%(\ref{P}), (\ref{barp}) and (\ref{hatp}), the test functions $P, \bar P$ and $\hat P$ satisfy the following relations
\be\label{Ps20}
\Lambda Q^0 :=  \frac 12 ( x(Q^0)' +2Q^0) + \frac{3\tilde \mu}{4\rho^2} (3Q^0 -(Q^0)^2), \quad 
\mathcal{L}^0 \Lambda Q^0 = -Q^0. 
\ee
Moreover, %$(\Lambda Q^0)^2 = \frac 14 x^2(Q^0)^2 -\frac 16 x^2 (Q^0)^3 +xQ^0(Q^0)' +(Q^0)^2$ and
$$
\int_\R \Lambda Q^0 % \frac{1}{2\rho^2} (\int_\R Q^0 -\frac 32\tilde \mu \int_\R (Q^0)^2) 
 =\frac 3{\rho^2}, \quad \int_\R \Lambda Q^0 Q^0 =\frac 9{2\rho^2}.
$$
\item \emph{Resonance functions}. Define  $\varphi^0 := -\frac{(Q^0)'}{Q^0}$. Then 
\be\label{phi12}
(\varphi^0)' =\frac 13 Q^0 -\frac{\tilde \mu}2 (Q^0)^2, \quad (\varphi^0)^2 = 1 -\frac 23 Q^0 + \frac{\tilde \mu}2 (Q^0)^2.
\ee
\item \emph{Integrals}. 
For any $p\geq 1$,
\be\label{Int22}
 \tilde \mu \int_\R (Q^0)^{p+2} = \frac{2(2p+1)}{3(1+p)}\int_{\R} (Q^0)^{p+1}-\frac{2p}{p+1} \int_{\R} (Q^0)^{p}. %\quad \int_\R Q^0=\int_\R (Q^0)^2 =6.
\ee
\end{enumerate}
\end{Cl}

\begin{proof}
A direct computation, see e.g. Claim \ref{exp6} for a similar proof.
\end{proof}

Now we proceed to give the explicit value the constants and functions related to system $(\Omega_{1,0})$, see (\ref{A10})-(\ref{B10}). 

\begin{Cl}[Resolution of $(\Omega_{1,0})$ for the Gardner equation]
Denote by $(a_{1,0}^0, A_{1,0}^0, B_{1,0}^0)$ the solution of the linear system $(\Omega_{1,0})$, for the Gardner nonlinearity. Then we have
$$
\begin{cases}
a_{1,0}^0 = \frac 23, \\
A_{1,0}^0 = -\frac 43 Q^0 + \tilde \mu (Q^0)^2, \quad B_{1,0}^0 = -2\varphi^0 + \kappa_{1,0}^0 (Q^0)',\\
b_{1,0}^0 = \lim_{+\infty} B_{1,0}^0=  -2,
\end{cases}
$$
with
\be\label{k10}
\kappa_{1,0}^0 = \frac{3\tilde \mu (\int_\R (Q^0)^2 -3\int_\R Q^0)}{ (3\tilde \mu-1) \int_\R (Q^0)^2 + \int_\R Q^0}  = -\frac {10}3 +o_{\tilde \mu}(1).
\ee
\end{Cl}
\smallskip

\begin{rem}
It is remarkable the similarity among the functions solution of the Gardner system $(\Omega_{1,0})$ and the corresponding ones for the quadratic nonlinearity (let $\tilde\mu \to 0$). 
\end{rem}

\begin{proof}
First of all, the explicit value of $(a_{1,0}^0, A_{1,0}^0, B_{1,0}^0)$ comes from a straightforward verification. More precisely, this triplet is a solution of the zeroth order system
$$\begin{cases}
(\mathcal L^{0} A_{1,0}^0)' + a_{1,0}^0(3Q^0 -2(Q^0)^2 + 2\tilde \mu (Q^0)^3) ' = (2Q^0 -3\tilde \mu (Q^0)^2)', & \\
(\mathcal L^{0} B_{1,0}^0)' +3a_{1,0}^0 (Q^0)'' -3(A_{1,0}^0)'' - (2Q^0-3\tilde \mu (Q^0)^2) A_{1,0}^0 = 2Q^0 -3\tilde \mu (Q^0)^2, &
\end{cases}
$$
which comes from (\ref{A10})-(\ref{B10}). In particular, we choose $\kappa_{1,0}$ such that
$
\int_\R B^0_{1,0} (Q^0)' =0.
$
The value of $b_{1,0}^0$ comes from the fact that $b_{1,0}^0 =-2\lim_{+\infty} \varphi^0 =-2$. On the other hand, from (\ref{balt}), one has
$$
b_{1,0}^0 = -\frac 13\int_\R Q^0 + \frac 12 \tilde \mu \int_\R (Q^0)^2 =-2.
$$
\end{proof}

Now we are able to prove (\ref{Nodef}). (Note that this is also a consequence of the integrability of the Gardner equation.) First, we claim that
$$
3 (A_{1,0}^0)'  + (2 Q^0 -3\tilde \mu (Q^0)^2) \int_0^x ( A_{1,0}^0 + a_{1,0}^0 Q^0) = 3 (A_{1,0}^0)'  -2 (2 Q^0 -3\tilde \mu (Q^0)^2)\varphi^0 = 0.
$$
This is an easy consequence of (\ref{phi12}) and the values of of $A_{1,0}^0$ and $a_{1,0}^0$.
Consider now the expression for $b_{2,0}$, $m=2$ in Claim \ref{CL}. Note that the term containing $B_{1,0}^0$ disappears. Replacing the values of $a_{1,0}^0$ and $A_{1,0}^0$, and using the recursive formula (\ref{Int22}), we have
\bee%\label{b202newbeta}
d_{\tilde \mu,0} & = &  -\frac 13 (b_{1,0}^0)^3 - \frac 32 \tilde \mu \int_\R Q^0(1+A_{1,0}^0)^3   -2b_{1,0}^0  +  \frac 12\int_\R A_{1,0}^0(1+ (A_{1,0}^0)^2)- \frac 13\int_\R Q^0 A_{1,0}^0  \nonumber \\
& & \quad  -\frac 29  \int_\R (Q^0)'^2 + \frac 29 \int_\R (Q^0)^2[1 - (2 -3\tilde \mu Q^0 ) (1+A_{1,0}^0) ] + \frac 43 \int_\R (Q^0)''A_{1,0}^0 \nonumber \\
& & \quad -\frac 12 \int_\R [ (2Q^0 -3\tilde \mu (Q^0)^2)(1+A_{1,0}^0) + 3(A_{1,0}^0)''  ]A_{1,0}^0  \nonumber   \\
& = &   -\frac 13 (b_{1,0}^0)^3  -\frac 32\tilde \mu \int_\R Q^0 + (\frac{11}2 \tilde \mu -\frac{28}9)\int_\R (Q^0)^2 + (-\frac 92\tilde \mu^2 +\frac 23 \tilde \mu +\frac{20}9) \int_\R (Q^0)^3\nonumber   \\
& & \quad +\tilde \mu (\frac{15}2\tilde \mu -\frac{13}{3}) \int_\R (Q^0)^4 -\frac 92\tilde \mu^3  \int_\R (Q^0)^5  +\frac 72\tilde \mu^3 \int_\R (Q^0)^6 -\frac 32 \tilde \mu^4 \int_\R (Q^0)^7 \nonumber   \\
& =&   -\frac 13 (b_{1,0}^0)^3 -\frac 49 \int_\R Q^0 + \frac 23 \tilde \mu \int_\R (Q^0)^2\\
&= &   \frac 13 b_{1,0}^0 (4- (b_{1,0}^0)^2) =0.
\eee
This proves (\ref{Nodef}).
\end{proof}

\begin{proof}[Proof of (\ref{K1})]
The proof of (\ref{K1}) is a consequence of (\ref{Nodef}) and the simple relationship
$$
d_{\tilde \mu, \nu} = d_{\tilde \mu, 0} + \nu \partial_\nu d_{\tilde \mu, \nu} + o_{\nu}(\nu) = \nu (\partial_\nu d_{0, \nu} + o_{\tilde \mu}(1)) + o_{\nu}(\nu),
$$
valid for any $\tilde \mu, \nu\in \R$ small enough. This result says that, in order to prove (\ref{K1}), we only need to compute the defect for first order expansion in $\nu$ of the quadratic nonlinearity $f(s) = s^2 + \nu s^p$, $p\geq 3$. Then we use the fact that $\tilde \mu = \mu(\ve) \sim \ve^{1/(p-2)}$ and $\nu = \ve$ to conclude.
Consequently, in what follows we are reduced to prove that $\partial_\nu d_{0, \nu} \neq 0 $ for all $\nu$ small enough.
\end{proof}

\begin{Cl}[Asymptotic expansions, case $m=2$, basic functions]\label{exp12}

Suppose now $f(s)=f_{0,\nu}(s) =s^2 + \nu s^p$. Let $Q^0(x) =\frac{3}{2\cosh^2 (x/2)}$ be the soliton solution for the quadratic case. Then the following asymptotic expansions hold. 
\begin{enumerate}
\item The soliton solution $Q$ for the nonlinearity $f$ can be expanded as
\be\label{ffq12}
Q= Q^0 + \nu Q^1 +o(\nu), \quad o(\nu) \in \mathcal{Y},
\ee
where $Q^1$ satisfies the equation $ \mathcal L^0 Q^1 :=-(Q^1) '' + Q^1 -  2Q^0Q^1= (Q^0)^p.$
We also have 
\be\label{ffp23}\begin{cases}
& f(Q) = (Q^0)^2 + \nu ((Q^0)^p +2Q^0Q^1) + o(\nu), \\
& f'(Q) = 2Q^0 +\nu (2Q^1+ p(Q^0)^{p-1}) + o(\nu),\\
& f''(Q) = 2 + \nu p(p-1) (Q^0)^{p-2} + o(\nu),
\end{cases}
\ee
where every term $o(\nu)\in \mathcal{Y}$ uniformly in $\nu<\nu_0$ small.
\item The operator $\mathcal L$ satisfies
$$
\mathcal L = \mathcal L^0  - \nu [ 2Q^1+ p(Q^0)^{p-1}] + o(\nu).
$$
\item From the definition of $Q$, we have%(\ref{P}), (\ref{barp}) and (\ref{hatp}), the test functions $P, \bar P$ and $\hat P$ satisfy the following relations
\be\label{Ps2}
\Lambda Q = \Lambda Q^0 + \nu \Lambda Q^1 +o(\nu),  \quad  \mathcal{L}^0 \Lambda Q^1 = -Q^1 +(2Q^1+p(Q^0)^{p-1})\Lambda Q^0. 
\ee
\item Let  $D(\nu):= \int_\R \Lambda Q Q$. Then
\be\label{int22}
D(\nu) =\frac 9{2} %\frac{3}{4\rho^2} (\int_\R Q^0 -\frac 32\tilde \mu \int_\R (Q^0)^2) 
+ \nu ( \int_\R \Lambda Q^1 Q^0 + \int_\R \Lambda Q^0 Q^1) + o(\nu).
%D(\ve) = \frac 92 + \ve( \int_\R \Lambda Q^1 Q^0 + \int_\R \Lambda Q^0 Q^1) + o(\ve).
\ee
\item \emph{Inverse functions}. We have
\be\label{inv0c2}
\mathcal{L}^0 \big[ 1 - \frac 4{3}\Lambda Q^0 \big] = 1-\frac 23 Q^0.
\ee
\be\label{inv02}
\mathcal{L}^0 \big[ (1-Q^0) (1+\frac 13x^2 Q^0)-Q^0 \big] = 1-\frac 83 \Lambda Q^0 + \frac 83(\Lambda Q^0)^2.
\ee
\be\label{inv0a2}
\mathcal{L}^0 \big[ -5  + \frac {68}9 Q^0 -6\Lambda Q^0 \big]  = -5 +16 Q^0  -\frac {68}9(Q^0)^2.
\ee
\be\label{inv0b2}
\mathcal{L}^0 \big[ 2 + \frac{20}3 \Lambda Q^0 - \frac{170}{27}Q^0 \big] = 2 -\frac{32}3 Q^0 +\frac{170}{27}(Q^0)^2.
\ee
\end{enumerate}
\end{Cl}

\begin{proof}
First of all, (\ref{ffq12})-(\ref{Ps2}) are a direct consequence of a Taylor expansion of the considered functions. The expression for $(\Lambda Q^0)^2$ comes from a simple computation. 

The expansion of $D(\nu)$ in (\ref{int22}) follows from a Taylor expansion and the fact that 
$$
\int_\R \Lambda Q^0 Q^0 = \frac 34 \int_\R (Q^0)^2 = \frac 92 .
$$
Finally we prove (\ref{inv02}), (\ref{inv0a2}) and (\ref{inv0b2}). These follow 
from the identities
$\mathcal L^0 1 = 1-2Q^0, \quad \mathcal L^0 Q^0 = -(Q^0)^2$, $\mathcal L^0 \Lambda Q^0 =-Q^0$, $\mathcal{L}^0 (x^2 Q^0)  =-2Q^0 -4x(Q^0)' -x^2 (Q^0)^2$, and $$ \mathcal{L}^0 (x^2 (Q^0)^2) = -2(Q^0)^2  -8xQ^0 (Q^0)' -3x^2(Q^0)^2 +\frac 43 x^2 (Q^0)^3.$$  
This finishes the proof.
\end{proof}

\smallskip

Now we proceed to give an asymptotic expansion of the constants and functions related to system $(\Omega_{1,0})$, see (\ref{A10})-(\ref{B10}). 

\begin{Cl}[Asymptotic expansions II, case $m=2$]\label{exp23}
There exists $\nu_0$ small enough such that for all $|\nu|\leq \nu_0$, the following holds. Let $f(s) = s^2 + \nu s^p $, then the corresponding solution to the system $(\Omega_{1,0})$ for this case can be expanded as follows:
\be\label{Mae}
a_{1,0} = \frac 23 + \nu a_{1,0}^1 + o(\nu), \quad A_{1,0} = -\frac 43 Q^0 + \nu A_{1,0}^1 + o(\nu) \in \mathcal Y,  \quad b_{1,0} = -2 + b_{1,0}^1\nu + o(\nu),
\ee
where $\nu^{-1}o(\nu)\to 0$ as $\nu \to 0$ and $A_{1,0}^1\in \mathcal{Y}$ is a solution of the following linear equation
\be\label{O101}
 (\mathcal L^{0} A_{1,0}^1)' + a_{1,0}^1(3Q^0 - 2(Q^0)^2)'  =  [p(Q^0)^{p-1} -\frac 43 (p-1)(Q^0)^p ]', 
\ee
Finally, the following two expressions are satisfied
\be\label{ab22}
a_{1,0}^1 =  - \frac 19\Big[ \frac{(p-3)(2p-1)}{p+1} \Big] \int_\R (Q^0)^p,
\ee
and
\be\label{ba22}
b_{1,0}^1 = \frac 12 \int_\R A_{1,0}^1+ \frac 12 a_{1,0}^1\int_\R Q^0 + \frac 13 \int_\R Q^1.
\ee
\end{Cl}
\smallskip

\begin{proof}
The proof of (\ref{Mae}) and  (\ref{O101}) is direct from Claim \ref{exp12} and (\ref{A10})-(\ref{B10}). To prove (\ref{ab22}), first note that from Claim \ref{exp12}
\bee
a_{1,0}^1 & =  & \frac 29 \Big[ \int_\R(1 - \frac 2{3} Q^0 )\Lambda Q^1  -  \frac 2{3} \int_\R Q^1 \Lambda Q^0 \Big] \nonumber \\
& =  & \frac 29 \int_\R [1 - \frac 4{3} \Lambda Q^0][-Q^1 +(2Q^1+p(Q^0)^{p-1})\Lambda Q^0] -  \frac 4{27} \int_\R Q^1 \Lambda Q^0 \\
& =& - \frac 29 \int_\R Q^1 [ 1-\frac 83\Lambda Q^0 + \frac 83 (\Lambda Q^0)^2 ]  + \frac{2 p}{9}\int_\R  (Q^0)^{p-1}\Lambda Q^0[1- \frac 43 \Lambda Q^0].
\eee
Thus from (\ref{inv02}) and Claim \ref{exp12} (i), we get after integration by parts
\bee
a_{1,0}^1&  = & \frac 29\int_\R (Q^0)^{p-1}\big[ \, p\Lambda Q^0 (1-\frac 43 \Lambda Q^0) -Q^0(1-Q^0)(1+\frac 13 x^2 Q^0) +(Q^0)^2 \big]  \nonumber \\
& = &  \frac 29\Big[ \int_\R (Q^0)^{p-1}\big[ (p-1)Q^0 +(2-\frac 43 p)(Q^0)^2 +\frac p2 x(Q^0)' -\frac 43 pxQ^0(Q^0)'  \big]  \nonumber \\
& & \quad \quad  -\frac 1{3} \int_\R x^2(Q^0)^{p+1} \big[ (p+1)   - (1+\frac 23p)Q^0  \big]   \Big]  \nonumber \\
& = &  \frac 29\Big[ (p-\frac 32) \int_\R (Q^0)^p +  (\frac {10}3 -\frac 43p - \frac 43 \frac 1{p+1}) \int_\R (Q^0)^{p+1}  \\
& & \qquad -\frac 1{3} \int_\R x^2(Q^0)^{p+1} \big[ (p+1)   - (1+\frac 23p)Q^0  \big] \Big] . 
\eee
Now, recall that from the equation satisfied by $Q^0$, $[(Q^0)^{p+1}]'' = (p+1)(Q^0)^{p+1}[ (p+1) - (1+\frac 23 p)Q^0] , $
so that 
$$
\int_\R x^2 (Q^0)^{p+1}[ (p+1) - (1+\frac 23 p)Q^0 ] = \frac{1}{p+1} \int_\R x^2 [(Q^0)^{p+1}]''  = \frac{2}{p+1}\int_\R (Q^0)^{p+1}.
$$
In conclusion, from (\ref{Int22}),
$$
a_{1,0}^1  =  \frac 29\Big[ (p-\frac 32) \int_\R (Q^0)^p +  (\frac {10}3 -\frac 43p -  \frac {2}{p+1} ) \int_\R (Q^0)^{p+1} \Big]   =  - \frac 19\Big[ \frac{(p-3)(2p-1)}{p+1} \Big] \int_\R (Q^0)^p,
$$
as desired. Finally, from (\ref{balt}) we obtain (\ref{ba22}).
\end{proof}

Now we deal with the second order system $(\Omega_{2,0})$ written in (\ref{eq:A20}), (\ref{eq:B20}), (\ref{F202}) and (\ref{G202}). 

\begin{Cl}[Asymptotic expansions III, case $m=2$]\label{exp223}
The following identity holds
$$
\partial_\nu d_{0, \nu}|_{\nu=0}  = - \Big[ \frac{(p-3)(2p-1) (24 -23p +3 p^2  +2p^3)}{36(p^2 -1)(p-2)}\Big] \int_\R (Q^0)^p, %\quad c_{2,p}= \int_\R (Q^0)^p \neq 0,
$$
for all $p\geq 3$.  %Here $c_{2,p}$ is given by (\ref{c2p}).
\end{Cl}
\begin{proof} 
The proof of the above result is equivalent to prove that for the nonlinearity $f_{0,\nu}(s) = s^2 +\mu s^p$, $p\geq 3$ and $\nu$ small, we have
$$
d_{0,\nu} =  b_{2,0}(f_{0,\nu}) + \frac 16 b_{1,0}^3(f_{0,\nu}) = - \nu \Big[ \frac{(p-3)(2p-1) (24 -23p +3 p^2  +2p^3)}{36(p^2 -1)(p-2)}\Big] \int_\R (Q^0)^p + o (\nu).
$$
First of all, note that we can expand $b_{2,0} = b_{2,0}^0 + \nu b_{2,0}^1 + o(\nu)$, with $b_{2,0}^0 = \frac 43$ (cf. \cite{MMcol1}, Lemma 3.1). By considering (\ref{b202new}) in Claim \ref{CL}, Claim \ref{exp12} and expanding at first order in $\nu$, we get
\bea\label{B1nul} 
b_{2,0}^1 & = &   -8b_{1,0}^1 + \frac 14 p(p-1) \int_\R (Q^0)^{p-2}(1+A_{1,0}^0)^3  +  \frac 12\int_\R A_{1,0}^1(1+3(A_{1,0}^0)^2)  - \frac 12 a_{1,0}^1\int_\R Q^0 A_{1,0}^0  \nonumber \\
& & - \frac 13 \int_\R Q^1 A_{1,0}^0- \frac 13 \int_\R Q^0 A_{1,0}^1  - a_{1,0}^1 \int_\R (Q^0)'^2 - \frac 49 \int_\R (Q^0)'(Q^1)'   \nonumber \\
& & + \frac 23 a_{1,0}^1 \int_\R (Q^0)^2(1 - 2(1+A_{1,0}^0) )  + \frac 29 \int_\R Q^0 Q^1(1 -2(1+A_{1,0}^0) ) \nonumber \\
& & + \frac 29 \int_\R Q^0[Q^1- 2Q^0A_{1,0}^1 -(2Q^1+ p(Q^0)^{p-1})(1+A_{1,0}^0) ] \nonumber  \\
& & -\frac 34a_{1,0}^1 \int_\R [ 2Q^0(1+A_{1,0}^0) + 3(A_{1,0}^0)''  ] A_{1,0}^0  -\frac 12 \int_\R [ 2Q^0(1+A_{1,0}^0) + 3(A_{1,0}^0)''  ]A_{1,0}^1   \nonumber \\
& & -\frac 12 \int_\R [ (2Q^1+ p(Q^0)^{p-1})(1+A_{1,0}^0) + 2Q^0 A_{1,0}^1  + 3(A_{1,0}^1)''  ]A_{1,0}^0  \nonumber \\
%& &  + \frac 12 \int_\R B_{1,0}^1 [ 3(A_{1,0}^0)' + 2Q^0\int_0^x (A_{1,0}^0  + a_{1,0}^0Q^0) ]   \\
& &  +\frac 12 \int_\R B_{1,0}^0( 2Q^1+ p (Q^0)^{p-1})\int_0^x (A_{1,0}^0  + \frac 23Q^0)  +\frac 12  \int_\R B_{1,0}^0 [ 3(A_{1,0}^1)'    \nonumber \\
& & + 2Q^0 \int_0^x (A_{1,0}^1  + a_{1,0}^1Q^0 + \frac 23 Q^1) ] +4 a_{1,0}^1  \int_\R (Q^0)''A_{1,0}^0 + \frac 43 \int_\R (Q^1)''A_{1,0}^0 + \frac 43 \int_\R (Q^0)''A_{1,0}^1. \nonumber %\\
\eea
%Now we replace the values of $b_{1,0}^1$ and $B_{1,0}^0$ given by (\ref{ba22}) and Claim \ref{exp23}. We arrange the above expression according to $a_{1,0}^1$, $A_{1,0}^1$, $Q^1$ and the rest terms. We obtain,
Now we arrange the above expression according to $a_{1,0}^1$, $A_{1,0}^1$, $Q^1$, $b_{1,0}^1$, $B^{1}_{1,0}$ and the rest terms. We obtain,
\bee
b_{2,0}^1 & = &  -8 b_{1,0}^1+ \frac 14 p(p-1) \int_\R (Q^0)^{p-2}(1+A_{1,0}^0)^3  -\frac p3 \int_\R (Q^0)^{p-1}(1+A_{1,0}^0)(\frac 32 A_{1,0}^0 + \frac 23 Q^0)\\
& &  + p\int_\R B_{1,0}^0 (Q^0)^{p-2} (Q^0)'   + \int_\R B_{1,0}^0 Q^0 \int_0^x A_{1,0}^1 + \frac 12 \int_\R A_{1,0}^1 F_A  +  \delta  a_{1,0}^1  +  \int_\R  Q^1 F_Q \\
& &  + \frac 23 \int_\R B_{1,0}^0 Q^0 \int_0^x Q^1,
\eee
where 
\bee
F_A & :=&  1+3(A_{1,0}^0)^2 - \frac 23 Q^0 - \frac 89 (Q^0)^2 - (2Q^0 +2Q^0A_{1,0}^0 +3(A_{1,0}^0)'' )  - ( 2Q^0 A_{1,0}^0 +3(A_{1,0}^0)'') \\
& & \qquad -3(B_{1,0}^0)' + \frac 83 (Q^0)'' ; 
\eee
\bee
\delta & := & -\frac 12 \int_\R Q^0 A_{1,0}^0 + \int_\R Q^0 (Q^0)'' - \frac 23 \int_\R (Q^0)^2 (1+2A_{1,0}^0)\\
& & \qquad  -\frac 34 \int_\R (2Q^0 +2Q^0A_{1,0}^0 +3(A_{1,0}^0)'')A_{1,0}^0 + 4 \int_\R (Q^0)'' A_{1,0}^0 ;
\eee
and
\bee
F_Q & := & -\frac 13 A_{1,0}^0 +\frac 49 (Q^0)'' - \frac 49 Q^0(1+2A_{1,0}^0) -(1+A_{1,0}^0)A_{1,0}^0-2\varphi^0 B_{1,0}^0 + \frac 43 (A_{1,0}^0)'' .
\eee

Note that we have used that 
\be\label{unique}
\int_\R B_{1,0}^0 (Q^0)' =0, \quad \int_0^x (A_{1,0}^0 + \frac 23 Q^0) = -2\varphi^0.
\ee
Now we use the expressions (\ref{ab22}) and (\ref{ba22}) in Claim \ref{exp23} to have
\bea
  b_{2,0}^1 + 2 b_{1,0}^1 & =  &  -  \frac{12} 5 a_{1,0}^1 + \frac 14 p(p-1) \int_\R (Q^0)^{p-2}(1-\frac 43Q^0)^3 + \frac {4p}9 \int_\R (Q^0)^{p}(1-\frac 43 Q^0) \nonumber \\ 
& & + p\int_\R (Q^0)^{p-1}\big[ 2-\frac{14}3 Q^0 +\frac{20}9 (Q^0)^2\big] + \frac 12 \int_\R A_{1,0}^1\Big[ -5 +16Q^0 -\frac{68}9 (Q^0)^2  \Big] \nonumber\\ % + 6 p\int_\R (1-\frac 73 Q^0 +\frac{10}3 (Q^0)^2) (Q^0)^{p-1}
& &     +  \int_\R  Q^1\Big[  2 -\frac{32}3 Q^0 +\frac{170}{27}(Q^0)^2  \Big]. \nonumber 
\eea
Note that we have also made use of (\ref{Int22}) with $\tilde \mu =0$ to obtain
\bee
& & - 3 \int_\R Q^0 -\frac 12 \int_\R Q^0 A_{1,0}^0 +  \int_\R Q^0 (Q^0)'' - \frac 23 \int_\R (Q^0)^2 (1+2A_{1,0}^0)  \\
& & \quad \quad \quad    -\frac 34 \int_\R (2Q^0 +2Q^0A_{1,0}^0 +3(A_{1,0}^0)'')A_{1,0}^0 + 4 \int_\R (Q^0)'' A_{1,0}^0 = -\frac{12}5.
\eee
Using (\ref{inv0a2}), (\ref{inv0b2}) and (\ref{O101}), we have
\bee
& &  b_{2,0}^1 + 2 b_{1,0}^1 =     - \frac{12} 5 a_{1,0}^1  + \frac 14 p(p-1)  \int_\R (Q^0)^{p-2}(1-\frac 43Q^0)^3 + \frac {4p}9 \int_\R (Q^0)^{p}(1-\frac 43 Q^0)  \nonumber\\ % + 6 p\int_\R (1-\frac 73 Q^0 +\frac{10}3 (Q^0)^2) (Q^0)^{p-1}
& & \quad + \frac 12 \int_\R\Big[-a^1_{1,0}(3Q^0 -2(Q^0)^2) + p(Q^0)^{p-1} -\frac 43 (p-1)(Q^0)^p \Big] \Big[ -5 -6\Lambda Q^0 +\frac {68}9 Q^0 \Big]  \\
& & \quad   +  \int_\R  (Q^0)^p\Big[ 2 + \frac{20}3 \Lambda Q^0 - \frac{170}{27}Q^0 \Big] +  p\int_\R (Q^0)^{p-1}\big[ 2-\frac{14}3 Q^0 +\frac{20}9 (Q^0)^2\big] .
\eee
\noindent
A simple computation using (\ref{Ps2}) and (\ref{Int22}) with $\tilde \mu =0$ shows that 
$$
\int_\R(3Q^0 - 2(Q^0)^2)( -5 - 6\Lambda Q^0 + \frac{68}9 Q^0)= -\frac{59}5.
$$
Thus, replacing the value of $a_{1,0}^1$ given by (\ref{ab22}),
\bee
& &  b_{2,0}^1 + 2 b_{1,0}^1 =  - \frac 7{18}\Big[ \frac{(p-3)(2p-1)}{p+1} \Big] \int_\R (Q^0)^p + \frac {4p}9 \int_\R (Q^0)^{p}(1-\frac 43 Q^0)  \\
& & \quad + \frac 14 p(p-1) \int_\R (Q^0)^{p-2}(1-\frac 43Q^0)^3 +  p\int_\R (Q^0)^{p-1}\big[ 2-\frac{14}3 Q^0 +\frac{20}9 (Q^0)^2\big]  \nonumber\\ % + 6 p\int_\R (1-\frac 73 Q^0 +\frac{10}3 (Q^0)^2) (Q^0)^{p-1}
& & \quad + \frac 12 \int_\R\Big[p(Q^0)^{p-1} -\frac 43 (p-1)(Q^0)^p \Big] \Big[ -5 -6\Lambda Q^0 +\frac {68}9 Q^0 \Big]  \\
& & \quad  + 2 \int_\R  (Q^0)^p  + \frac {10}3 \frac{2p+1}{p+1}\int_\R  (Q^0)^{p+1}  - \frac{170}{27} \int_\R  (Q^0)^{p+1}.   
\eee
Simplifiying, we get
\bee
&& b_{2,0}^1 + 2 b_{1,0}^1  = - \frac 7{18} \frac{(p-3)(2p-1)}{p+1} \int_\R (Q^0)^p +  \frac 14p(p-1)\int_\R (Q^0)^{p-2} -  p(p-\frac 12) \int_\R (Q^0)^{p-1} \\
& & \quad + (\frac 16 - \frac {13}9p + \frac 43 p^2 ) \int_\R (Q^0)^p + (- \frac {16}{27} + \frac{32}{27} p  + \frac {2}3\frac 1{p+1} -\frac{16}{27} p^2 )\int_\R (Q^0)^{p+1}.
\eee
Using (\ref{Int22}) with $\tilde \mu =0$ and the fact that $p\geq 3$, we finally obtain
\bee
b_{2,0}^1 + 2 b_{1,0}^1&  =&     \Big[ - \frac 7{18}\frac{(p-3)(2p-1)}{p+1}  +  \frac 1{36}p%(p-1) 
\frac{(2p-1)(2p-3)}{%(p-1)
p-2} -  \frac 13 p(p- \frac 12) \frac{(2p-1)}{p-1} +\\
& & \quad   + (\frac 16- \frac {13}{9}p + \frac 43 p^2)  + \frac{3p}{1+2p}(- \frac {16}{27} + \frac{32}{27} p  + \frac 23\frac {1}{p+1} -\frac{16}{27} p^2  ) \Big]\int_\R (Q^0)^p\\
%& = & -\Big[ \frac{72 -237 p +218 p^2  -61 p^3  - 8 p^4  +4 p^5}{36(p^2 -1)(p-2)}\Big] \int_\R (Q^0)^p.
& =& -\frac{(p-3)(2p-1) (24 - 23p +3 p^2  + 2p^3)}{36(p^2 -1)(p-2)}\int_\R (Q^0)^p.
%& = & \Big[ \frac{88p^6 -180p^5 +574p^4 -2199p^3 +2092p^2 +453p -648}{36(p^2 -1)(p-2) (1+2p)}\Big] \int_\R (Q^0)^p.
\eee
Let us define, for $p$ real, $ f(p):= 24 - 23p + 3 p^2  + 2p^3.$
Then we have
\be\label{pos1}
f(p)\geq 36 \quad \hbox{ for all } p\geq 3.
\ee
It is clear that this last affirmation allows us to conclude the proof. Let us prove (\ref{pos1}). Note that $f(3) = 36$ and $f'(p)$ is given by $f'(p)  = 6 p^2 + 6p - 23 >0$ for all $p\geq 3$. This implies (\ref{pos1}). The proof is complete. 
\end{proof}

\begin{rem}\label{2inf}
First of all, note that in the above expression we recover the integrability condition of the Gardner equation $(p=3)$.
Furthermore, note that this term is divergent when we formally take the limit $p\downarrow 2$ and the equation approaches the integrable case. This can be explained by the higher regularity needed in our results ($f\in C^{4}$ for $m=2$), to justify the asymptotics. Indeed, from (\ref{F202}) and (\ref{b202new}) we need at least $f\in C^3(\R)$, and $f(s):=s^2 +\ve s^p$ is not $C^3$ at zero as $p\downarrow 2$, $p>2$. In addition, the terms in (\ref{F202}), (\ref{G202})
$$
\frac 12(f''(Q) -2)', \; \frac 12 (f''(Q) -2),
$$
vanish in the integrable case $m=p=2$. For the computation in the quadratic case, see Proposition 2.1 and Lemma 3.1 in \cite{MMcol1}.
\end{rem}

%%%%%%%%%%%%%%%%%%%%%%%%%%%%%%%%%%%%%%
%%%%%%%%%%%%%%%%%%%%%%%%%%%%%%%%%%%%%%
%%%%%%%%%%%%%%%%%%%%%%%%%%%%%%%%%%%%%%
%%%%%%%%%%%%%%%%%%%%%

{\bf Acknowledgments}. The author would like to thank Yvan Martel and Frank Merle for their continuous encouragement in the elaboration of this work. Part of this work was written at DIM, Universidad de Chile. The author thanks the DIM members for their warm hospitality. Finally, the author wishes to thank J. Colliander for pointing him out the case of the Gardner nonlinearity.

\end{document}